\documentclass[11pt]{article} 
\usepackage{authblk}
\usepackage[utf8]{inputenc}
\usepackage{graphicx}
\usepackage{placeins}
\usepackage{pgfplots}
\usepgfplotslibrary{groupplots}
\usepackage{float} 
\usetikzlibrary{arrows.meta}
\usepackage{bbm}
\usepackage{appendix}
\usepackage{chngcntr}
\usepackage{apptools}
\AtAppendix{\counterwithin{lemma}{section}}
\usepackage{aliascnt}
\usepackage{hyperref}
\hypersetup{
	colorlinks   = true, 
	urlcolor     = blue, 
	linkcolor    = blue!90!black, 
	citecolor   = black 
}

\definecolor{lavender}{rgb}{0.9, 0.9, 0.98}
\definecolor{gainsboro}{rgb}{0.86, 0.86, 0.86}

\usepackage[backend=bibtex, style=numeric, maxbibnames=99]{biblatex}
\DeclareFieldFormat{prefixnumber}{\mkbibbold{#1}}
\DeclareFieldFormat{labelnumber}{\mkbibbold{#1}}

\AtBeginBibliography{%
	\DeclareFieldFormat{prefixnumber}{#1}%
	\DeclareFieldFormat{labelnumber}{#1}}

\usepackage{tikz}
\usetikzlibrary{shapes,arrows,shadows,patterns,pgfplots.groupplots,decorations.markings,chains,calc}
\usetikzlibrary{datavisualization}
\usetikzlibrary{datavisualization.formats.functions}
\usepackage{amsgen,amssymb,amsfonts,amsbsy,amsmath,amscd,amsthm}
\usepackage{mathrsfs}
\usepackage{comment} 
\usepackage{geometry}
\newtheorem{theorem}{Theorem}
\newaliascnt{lemma}{theorem}
\newtheorem{definition}{Definition}
\newtheorem{lemma}[lemma]{Lemma}
\newtheorem{corollary}{Corollary}[theorem]

\newtheorem{remark}{Remark}
\aliascntresetthe{lemma}

\newtheoremstyle{named}{}{}{\itshape}{}{\bfseries}{.}{.5em}{\thmnote{#3}}
\theoremstyle{named}

\title{A block-Toeplitz approach to the spectrum of the Neumann-Poincaré operator on self-similar chains}
\author{Matias Ruiz\thanks{\noindent E-mail address: \texttt{mr447@leicester.ac.uk}\\The author acknowledges support from \textit{The Royal Society of Edinburgh}.
}}
\date{}
\affil{\textit{\small{School of Computing and Mathematical Sciences,
			University of Leicester,
			Leicester, UK}}}

\begin{document}
\maketitle
\begin{abstract}
We characterise the spectrum of the Neumann-Poincar\'e operator acting on a mean-zero subspace of the Sobolev-Slobodeckij space $W^{s,p}(\Gamma)$, where $\Gamma$ is a non-Lipschitz domain, consisting of an infinite self-similar chain of disjoint smooth domains accumulating at a limit point. We exploit the discrete geometric scaling of the chain to obtain an exact block-Toeplitz representation of the Neumann-Poincaré operator, as well as the single- and double-layer operators. The corresponding operator-valued symbol depends entirely on the single scaling parameter $\alpha=(d-1)/p-s$ and belongs to the Wiener algebra whenever $0<\alpha<d$. This allows us to exploit block-Toeplitz operator theory to characterise the essential spectrum and Fredholm regions within this regime. In the energy space $H^{-1/2}_0(\Gamma)$, corresponding to $\alpha = d/2$, we prove the spectrum is real and characterise any isolated eigenvalues outside the essential spectrum through an operator-valued Wiener-Hopf factorisation. Furthermore, we apply an operator-valued Szeg\H{o} limit theorem to derive the asymptotic spectral distribution for large finite truncations and establish an eigenvalue counting formula based on the operator symbol. We illustrate these results through explicit analytical computations for a chain of concentric annuli and numerical approximations for a chain of disks.
\end{abstract}
	\smallskip



\section{Introduction}\label{sec:introduction}

\subsection{Problem formulation}\label{subsec:problem fomulation}

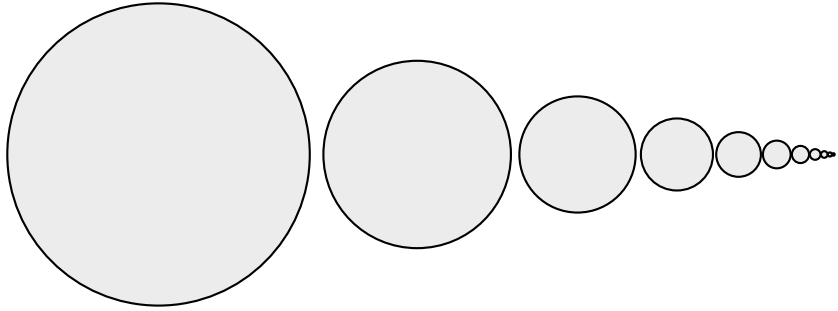
\begin{figure}[!htbp]
	\centering
	\begin{tikzpicture}[scale=2]
		
		\def\r{0.62}      
		\def\a{4.5}       
		\def\nDisks{10}    
		
		\foreach \j in {0,...,\nDisks} {
			\pgfmathsetmacro{\rad}{(\r)^(\j)}
			\pgfmathsetmacro{\cx}{-\a*(\r)^(\j)}
			\draw[thick, fill=gray!15] (\cx,0) circle (\rad);
		}
		
		\foreach \j in {0,...,\nDisks} {
			\pgfmathsetmacro{\rad}{(\r)^(\j)}
			\pgfmathsetmacro{\cx}{-\a*(\r)^(\j)}
		}

	\end{tikzpicture}
	\caption{Geometry of a self-similar chain of disjoint disks. The reference domain $\Gamma_0$ is scaled successively by a geometric factor $r \in (0,1)$ such that $\Gamma_j = r^j\Gamma_0$, generating an infinite sequence of boundaries accumulating at the origin.}
	\label{fig:disks}
\end{figure}

\paragraph{The Neumann-Poincaré operator.} Let $\Omega \subset \mathbb{R}^d$, $d = 2,\,3$ be a bounded open domain. For a given density $\varphi$ defined on a boundary surface $\Gamma: =\partial \Omega$, the \textit{Neumann-Poincaré operator} $K^*_{\Gamma}$ is given by the principal value integral
$$ K^*_{\Gamma} \varphi(x) := \mathrm{p.v.} \int_\Gamma \frac{\partial G}{\partial \nu_x}(x, y) \varphi(y) \, d\sigma(y), \quad x \in \Gamma, $$
where $\nu_x$ is the outward-pointing unit normal vector at $x \in \Gamma$ and $G(x,y)$ is the fundamental solution to the Laplace equation in $\mathbb{R}^d$:
$$
G(x,y) \;=\;
\begin{cases}
	\dfrac{1}{2\pi}\,\log|x-y|, & d = 2,\\[2mm]
	-\dfrac{1}{4\pi}\,|x-y|^{2-d}, & d = 3.
\end{cases}
$$
\paragraph{Self-similar chains.} In this work, we consider the specific case where $\Gamma$ is formed by an infinite sequence of smooth domains staggered in  a self-similar way (e.g. Fig~\ref{fig:disks}). Let $\Omega_0 \subset \mathbb{R}^d$ be a bounded and connected domain with $C^{\infty}$ boundary $\Gamma_0 := \partial \Omega_0$, and assume that the origin $0 \notin \overline{\Omega}_0$. For a fixed geometric scaling factor $r \in (0, 1)$, we define the scaled domains $\Omega_j := r^j \Omega_0$. We assume that $r$ and $\Omega_0$ are chosen such that these components are mutually disjoint, meaning $\overline{\Omega}_j \cap \overline{\Omega}_k = \emptyset$ for all $j \neq k$. The self-similar chain is then defined as
$$
\Omega  := \bigcup_{j\geq 0} \Omega_j.
$$
Thus, the infinite self-similar boundary $\partial \Omega = \Gamma\cup\{0\}$, where 
$$
\Gamma := \bigcup_{j=0}^\infty \Gamma_j,\quad \Gamma_j:= r^j\Gamma_0.
$$
Note that $\Gamma$ is not Lipschitz since every neighbourhood of the origin meets infinitely many components of the chain, so no neighbourhood of that point is the graph of a Lipschitz function.

We shall denote by $\Gamma^{(N)} := \bigcup_{j=0}^{N-1} \Gamma_j$ the finite truncated boundary of $N$ components.

\paragraph{Functional setting.} Let $W^{s,p}(\Gamma)$ denote the usual Sobolev-Slobodeckij space on the boundary, with $p \in (1, \infty)$ and $s \in (-1,1)$. For $s=0$, this space coincides with the standard Lebesgue space $L^p(\Gamma)$. For $s \in (0, 1)$, it is equipped with the standard norm
$$
\|\varphi\|_{W^{s,p}(\Gamma)}^p := \|\varphi\|_{L^p(\Gamma)}^p + |\varphi|_{W^{s,p}(\Gamma)},
$$
where $|\varphi|_{W^{s,p}(\Gamma)}$ denotes the Gagliardo seminorm:
\begin{equation}\label{eq:Gagliadrdo}
	|\varphi|_{W^{s,p}(\Gamma)} := \int_{\Gamma} \int_{\Gamma} \frac{|\varphi(x) - \varphi(y)|^p}{|x - y|^{d-1+sp}} \, d\sigma(x) \, d\sigma(y).
\end{equation}
For $s \in (-1, 0)$, the space is defined via $L^2$-duality with $W^{-s, p'}(\Gamma)$, with $1/p + 1/p' =1$.

In this paper we focus on the Neumann-Poincaré operator acting on the component-wise mean-zero subspace of $W^{s,p}(\Gamma)$, defined by
$$
W^{s,p}_0(\Gamma) := \left\{ \varphi \in W^{s,p}(\Gamma) \;:\; \int_{\Gamma_j} \varphi \, d\sigma = 0 \quad \text{for all } j \ge 0 \right\},
$$
As will be established, $W^{s,p}_0(\Gamma)$ is invariant under $K^*_\Gamma$.

Restricting the analysis to this mean-zero subspace is a natural rather than technical choice. Physically, the zero-mean condition expresses the electrostatic charge neutrality of each isolated component, which is a fundamental constraint in the physical contexts this operator models. This also corresponds to the componentwise compatibility condition for the Laplace boundary value problems with Neumann data.

We are in particular interested in the cases $W_0^{0,p}(\Gamma)=L_0^p(\Gamma)$, relevant to boundary value problems, and $W_0^{s,2}(\Gamma)=H_0^s(\Gamma)$, relevant to transmission problems. However, working in the general space $W_0^{s,p}(\Gamma)$ treats both settings simultaneously while making explicit the interplay between the parameters $s$ and $p$.

\paragraph{Main problem.}  Our main objective is to study
$$
\sigma(\widetilde K^*_{\Gamma}),
$$
i.e. the spectrum of $\tilde K^*_{\Gamma} := K^*_{\Gamma}|_{ W_0^{s,p}(\Gamma)}: W_0^{s,p}(\Gamma)\to W_0^{s,p}(\Gamma)$, for the specific values of $s$ and $p$ where $K^*_{\Gamma}$ is bounded.

\subsection{Background and motivation}

Boundary integral operators are fundamental tools in the mathematical analysis and numerical solution of boundary value and transmission problems for elliptic partial differential equations \cite{McLean2000, hsiao2008boundary}. Among these, the Neumann-Poincar\'e operator occupies a particularly prominent place. Historically, it arose as a central object in the study of classical Neumann and Dirichlet boundary value problems in electrostatics, potential theory, and continuum mechanics \cite{kellogg1929foundations, ammari2009layer}. In recent years, the Neumann-Poincaré operator has attracted substantial renewed interest due to its foundational role in mathematical wave physics, particularly in the study of subwavelength metamaterials and plasmonics \cite{Grieser:14, Ammari_book:18, bonnet2012t}. Consequently, a comprehensive understanding of its spectrum is essential both for establishing the solvability of boundary transmission problems and for characterizing subwavelength resonances in wave physics.

For bounded domains with smooth boundaries, the spectral theory of the Neumann-Poincar\'e operator is classical \cite{khavinson2007poincare}. The operator is compact and admits a symmetrisation on the physically relevant Sobolev space $H^{-1/2}(\Sigma)$ (with $\Sigma$ denoting the domain boundary), thus yielding a strictly real and purely discrete spectrum that accumulates only at zero.

When the boundary loses its smoothness, however, the compactness of the Neumann-Poincaré operator is lost, and the spectral picture becomes significantly richer. On Lipschitz domains with corners, the symmetrisation on $H^{-1/2}(\Sigma)$ persists, ensuring that the spectrum remains strictly real. Yet, the loss of compactness generates a nontrivial essential spectrum consisting of a connected interval determined entirely by the extremal angles of the boundary \cite{BonnetBenDhiaChesnelClaeys2013, perfekt2014spectral, PerfektPutinar2017, bonnetier2019characterization, perfekt2021plasmonic}. Subsequent investigations have also revealed the presence of eigenvalues embedded within this continuous band \cite{helsing2017classification, li2019embedded, LiPerfektShipman2022, bonnet2021complex, faria2025complex}. Similar essential spectrum phenomena have since been rigorously characterised for three-dimensional singularities, including rotationally symmetric conical points \cite{helsing2018spectra, perfekt2019transmission, bonnet2021maxwell}, and domains with edges \cite{de2022quasi}.

Remarkably, because the Neumann-Poincaré operator is generally non-self-adjoint, its spectral properties are highly sensitive to the choice of the ambient function space. Transitioning the analysis from $H^{-1/2}(\Sigma)$ to $L^2(\Sigma)$ or general $L^p(\Sigma)$ spaces causes the essential spectrum to detach from the real line and expand into complex curves or regions \cite{perfekt2021plasmonic}. More recently, attention has expanded to Lipschitz domains with highly oscillatory components respecting a local dilation-invariant property \cite{ChandlerWildeHaggerPerfektVirtanen2023}. Utilizing Floquet transform techniques (a discrete analogue to the continuous Mellin transform standardly used for isolated singularities like corners) a characterization of the essential spectrum on $L^2(\Sigma)$ was achieved.

This paper aims to contribute to the spectral study of the Neumann-Poincaré operator by investigating a \textit{non-Lipschitz} domain $\Gamma$ formed by an infinite sequence of components accumulating at a single point. The motivation is twofold. From a physical perspective, such self-similar chains have been proposed in the plasmonics literature as a mechanism for field enhancement and energy concentration \cite{LiStockmanBergman2003, DingDengYanCabriniZuckermannBokor2010, EssoneMezemeBrosseau2012, HoppenerLapinBharadwajNovotny2012}. From a mathematical perspective, we aim to give what is, to our knowledge, the first spectral analysis of the Neumann--Poincar\'e operator on a non-Lipschitz boundary. Doing so, we put forward a novel framework for studying the spectral properties of layer potentials in dilation-invariant geometries. Indeed, while $\Gamma$ is dilation-invariant (and therefore typically studied via Mellin and Floquet techniques), the fact that all of its components are exact scaled copies of one another allows us to represent the Neumann-Poincar\'e operator exactly in block-Toeplitz form, rather than modulo compact errors. This in turn allows us to apply well-established results in block-Toeplitz operator theory cleanly and directly, yielding not only the essential spectrum and the Fredholm index but also an operator-valued Wiener--Hopf factorisation \cite{GohbergLeiterer1973, GohbergLeiterer2009} and a Szeg\H{o}-type trace asymptotic \cite{MirandaTilli2000, BottcherSilbermann2006}. As such, this approach may provide a conceptual template for obtaining analogous results in continuously dilation-invariant geometries, such as corners or conical vertices, where logarithmic radial coordinates naturally yield Wiener--Hopf operators on the half-line (the continuous equivalent of Toeplitz operators).

Our block-Toeplitz framework relies fundamentally on the scale invariance of the underlying Green's function. Since fundamental solutions in elastostatics and hydrostatics share analogous scaling properties, and their associated Neumann-Poincar\'e-type operators display similar function-space dependence \cite{Mitrea1999, mitrea2002spectra}, we expect these techniques to carry over to the study of the layer potentials on self-similar chains in those systems, and more broadly in any second-order elliptic system with comparable scale invariance properties.

\subsection{Our approach and contributions}
The central strategy of this work is to establish a similarity transformation that identifies the Neumann--Poincar\'e operator on $W_0^{s,p}(\Gamma)$ with a block-Toeplitz operator acting on the sequence space $\ell^p\big(\mathbb{N}, W_0^{s,p}(\Gamma_0)\big)$:
\begin{equation*}
	\widetilde K^*_{\Gamma} \sim 
	\begin{pmatrix} 
		a_0 & a_{-1} & a_{-2} & \cdots \\ 
		a_1 & a_0 & a_{-1} & \cdots \\ 
		a_2 & a_1 & a_0 & \cdots \\ 
		\vdots & \vdots & \vdots & \ddots 
	\end{pmatrix}.
\end{equation*}
Here, the matrix entries $a_m$ are compact operators acting on the reference space $W_0^{s,p}(\Gamma_0)$. These operator-valued coefficients (which will be made precise later) depend explicitly on the geometric and functional parameter
\begin{equation*}
	\alpha := \frac{d-1}{p} - s.
\end{equation*}
This representation allows us to deduce the spectral properties of the global operator directly from its operator-valued symbol:
\begin{equation*}
	\kappa_\alpha(\theta) = \sum_{m \in \mathbb{Z}} a_m e^{im\theta}.
\end{equation*} 
Exploiting this framework, the present paper delivers three primary contributions:

\paragraph{Main result 1: Essential spectrum, Fredholm index and the spectrum.} For the scaling regime $0 < \alpha < d$, we establish that the essential spectrum of $\widetilde K^*_{\Gamma}$ is given exactly by the union of the spectra of the operator-valued symbol across the Floquet phase parameter:
$$
\sigma_{\mathrm{ess}}(\widetilde K^*_{\Gamma}) = \bigcup_{\theta \in [-\pi, \pi]} \sigma(\kappa_\alpha(\theta)).
$$
Outside these continuous bands we provide a formula for the Fredholm index via the winding number of finite-dimensional symbol truncations, paving the way for direct numerical evaluation. See Theorem~\ref{thm:essential-spectrum}, Theorem~\ref{thm:fredholm-index} and Corollary~\ref{cor:index-winding}.

\paragraph{Main result 2: Physical spectrum and discrete edge states.} In the energy space $H^{-1/2}_0(\Gamma)$, corresponding to $\alpha = d/2$, we prove that the projected block-Toeplitz operator is boundedly similar to a self-adjoint operator, guaranteeing a real spectrum. We then characterise the emergence of isolated eigenvalues outside the essential spectrum through an operator-valued Wiener--Hopf factorisation, linking the exact dimension of these highly localised discrete edge states to the non-canonical partial indices of the symbol. See Theorem~\ref{thm:physical-spectrum} and Theorem~\ref{thm:wiener-hopf}.

\paragraph{Main result 3: Large finite chains and the density of states.} We link the infinite-chain theory back to physical, finite systems by analysing truncations consisting of $N$ components. Writing $\{\lambda^{(N)}_j\}_{j\ge1}$ for the non-zero eigenvalues of  $\widetilde K^*_{\Gamma^{(N)}}$ and $\{\mu_j(\theta)\}_{j\ge1}$ for those of the symbol $\kappa_{d/2}(\theta)$, operator-valued Szeg\H{o} limit theorems give
$$
\lim_{N\to\infty}\frac{1}{N}\sum_{j\ge1} f\bigl(\lambda^{(N)}_j\bigr)
=\frac{1}{2\pi}\int_{-\pi}^{\pi}\sum_{j\ge1} f\bigl(\mu_j(\theta)\bigr)\,d\theta
$$
for every continuous $f$ vanishing near the origin. Taking $f$ to approximate the indicator of a spectral window yields an eigenvalue counting formula directly from the fibre spectra of $\kappa_{d/2}(\theta)$. See Theorem~\ref{thm:szego-chain} and Corollary~\ref{cor:density-of-states}.

\medskip 
We illustrate these results with an explicit, computable example of a self-similar chain of concentric two-dimensional annuli.

\subsection{Outline of the paper}

The remainder of this paper is structured as follows. Section~\ref{sec:preliminaries} briefly reviews classical results of layer potentials on Lipschitz domains that will be useful for our analysis. In Section~\ref{sec:layer potentials infinite chain}, we introduce the layer potentials on the infinite self-similar chain and define the projected layer potentials on the mean-zero subspace $W_0^{s,p}(\Gamma)$.  Section~\ref{sec:sequence spaces} establishes an isomorphisms between $W_0^{s,p}(\Gamma)$ and the sequence spaces $\ell^p(\mathbb N,W_0^{s,p}(\Gamma_0))$; this is the key step of our construction. Section~\ref{sec:block-Toeplitz} builds on these isomorphisms to derive the block-Toeplitz representations of the layer potentials. Section~\ref{sec:spectrum} builds on the block-Toeplitz representations to carry out the spectral analysis of the infinite chain, where we characterise the essential spectrum, the Fredholm index, and the emergence of discrete edge states via Wiener--Hopf factorisation. In Section~\ref{sec:asymptotics}, we apply operator-valued Szeg\H{o} limit theorems to derive the asymptotic eigenvalue distribution and density of states for large finite truncations. Section~\ref{sec:annuli} presents an explicit, computable example of the theory using a self-similar chain of concentric two-dimensional annuli. Section~\ref{sec:disks} provides further numerical illustrations by discretising the Neumann-Poincaré operator defined on the chain of disks of Figure~\ref{fig:disks}. Finally, Section~\ref{sec:conclusion} offers concluding remarks. 

The paper is complemented with three appendices, which collect, in an abstract setting, the operator-valued Toeplitz theory used throughout. Appendix~\ref{app:toeplitz} covers the Wiener algebra and the Fredholm criterion, Appendix~\ref{app:wiener-hopf} the Wiener--Hopf factorisation and the partial indices, and Appendix~\ref{app:szego} the Szeg\H{o}-type trace asymptotic, which we prove in the form required here.

\subsection{Notation and conventions}

Throughout this paper, we adopt the following standard mathematical notation and conventions. 

\begin{itemize}
	\item \textbf{Operators and spaces}. For Banach spaces $X$ and $Y$, $\mathcal{L}(X, Y)$ denotes the space of bounded linear operators from $X$ to $Y$, and we write $\mathcal{L}(X) := \mathcal{L}(X, X)$. The set of compact operators on $X$ is denoted by $\mathcal{K}(X)$. When $X$ is a Hilbert space, $\mathfrak{S}_q(X)$ denotes the standard Schatten class of order $q \in [1, \infty)$. 
	\item\textbf{ Spectrum and Fredholm operators}. For an operator $A \in \mathcal{L}(X)$, its spectrum is denoted by $\sigma(A)$. The essential spectrum, denoted by $\sigma_{\mathrm{ess}}(A)$, is defined as the set of $\lambda \in \mathbb{C}$ for which the operator $\lambda I - A$ is not Fredholm.
	\item \textbf{Inner products, duality, and adjoint.} The notation $\langle \cdot, \cdot \rangle$ stands for the standard $L^2$ inner product. When applied to functions in fractional Sobolev spaces, it is naturally extended by continuity to denote the duality pairing between $W^{s,p}$ and its dual space $W^{-s,p'}$, where $1/p + 1/p' = 1$.
	$A^*$ denotes the $L^2$ adjoint of an operator $A$.
	\item \textbf{Inequalities.} We use the notation $a \lesssim b$ to indicate that there exists a universal constant $C > 0$, independent of the critical parameters (such as the geometric scaling factor, the component index, or the matrix size), such that $a \le C b$. We write $a \asymp b$ to indicate two-sided bounds, i.e., $a \lesssim b$ and $b \lesssim a$.
	\item \textbf{Indicator function.} We denote by $\chi_{E}$ the indicator function of a subset $E\subset \mathbb R^d$.
\end{itemize}

\section{Preliminaries: layer potentials on Lipschitz domains}\label{sec:preliminaries}

We recall some properties of classical layer potentials on Lipschitz domains \cite{McLean2000, ammari2009layer}. Throughout this section $\Omega\subset\mathbb R^d$, $d=2,3$, is a bounded open set with Lipschitz boundary $\Sigma:=\partial\Omega$, not assumed connected or simply connected; the bounded connected components of $\mathbb R^d\setminus\overline\Omega$ are called the \emph{holes} of $\Omega$. Note that every truncation $\Gamma^{(N)}=\bigcup_{j=0}^{N-1}\Gamma_j$ of the self-similar chain is a finite disjoint union of $C^\infty$ closed surfaces and is therefore covered by all of the following.

\paragraph{Definition and mapping properties.}
For a density $\varphi$ on $\Sigma$ the single layer operator, the Neumann--Poincar\'e operator and the double layer operator are
\begin{align*}
S_\Sigma\varphi(x)&:=\int_\Sigma G(x,y)\varphi(y)\,d\sigma(y),\\
K^*_\Sigma\varphi(x)&:=\mathrm{p.v.}\!\int_\Sigma\frac{\partial G}{\partial\nu_x}(x,y)\varphi(y)\,d\sigma(y),\\
K_\Sigma\varphi(x)&:=\mathrm{p.v.}\!\int_\Sigma\frac{\partial G}{\partial\nu_y}(x,y)\varphi(y)\,d\sigma(y),
\end{align*}
with $\nu$ the outward normal. Thus $K_\Sigma$ is the formal $L^2(\Sigma)$ adjoint of $K^*_\Sigma$, and $S_\Sigma$ is self-adjoint in $L^2(\Sigma)$. For Lipschitz $\Sigma$, $K_\Sigma$ and $K^*_\Sigma$ are bounded on $L^p(\Sigma)$, $p\in(1,\infty)$, and on $H^{-1/2}(\Sigma)$, while $S_\Sigma:H^{-1/2}(\Sigma)\to H^{1/2}(\Sigma)$ \cite{ref:cmm,ref:verchota, Costabel1988}. If $\Sigma$ is $C^\infty$, then $S_\Sigma$ and $K^*_\Sigma$ are pseudodifferential operators of order $-1$ \cite{MiyanishiRozenblum2019}, so that for all $p\in(1,\infty)$ and $s\in\mathbb{R}$
\begin{equation}\label{eq:smoothing}
	S_\Sigma:W^{s,p}(\Sigma)\to W^{s+1,p}(\Sigma),
	\qquad
	K_\Sigma,\ K^*_\Sigma:W^{s,p}(\Sigma)\to W^{s+1,p}(\Sigma).
\end{equation}
(On a $C^\infty$ surface the scale is understood for every $s\in\mathbb{R}$, the
restriction $|s|<1$ being needed only on $\Gamma$, where we use the Gagliardo
semi-norm for the definition of the space.)

\paragraph{Jump relations and the energy identity.}
For $\varphi\in H^{-1/2}(\Sigma)$ the potential $u:=S_\Sigma\varphi$ is harmonic off $\Sigma$, continuous across it, and satisfies the jump relation $\partial_\nu u|_\pm=(\pm\frac12I+K^*_\Sigma)\varphi$. If $d=3$, or $d=2$ and $\int_\Sigma\varphi\,d\sigma=0$, then $u$ decays at infinity and
\begin{equation}\label{eq:energy}
	\big\langle\varphi,-S_\Sigma\varphi\big\rangle=\int_{\mathbb R^d}|\nabla u|^2\,dx .
\end{equation}

\paragraph{Definiteness in the mean-zero space.}
The operator $-S_\Sigma$ is positive definite on $H^{-1/2}(\Sigma)$ for $d=3$, whereas for $d=2$ its definiteness on the full space depends on the logarithmic capacity of $\Omega$. On the mean-zero subspace $H^{-1/2}_0(\Sigma)$, consisting of the densities with vanishing average on the boundary of each connected component of $\Omega$, this dichotomy disappears: there is $c>0$ with
\begin{equation}\label{eq:coercivity}
	\big\langle\varphi,-S_\Sigma\varphi\big\rangle\ \ge\ c\,\|\varphi\|^2_{H^{-1/2}(\Sigma)},
	\qquad\varphi\in H^{-1/2}_0(\Sigma),
\end{equation}
in both dimensions. By \eqref{eq:energy} the two sides of \eqref{eq:coercivity} carry the same homogeneity under dilations, so $c$ is unchanged if $\Sigma$ is replaced by a dilate of $\Sigma$ and norms on the dilated copy are taken by pullback. This dilation invariance will make this estimate usable uniformly along the self-similar chain.

\paragraph{Symmetrisation.}
The Plemelj symmetrisation principle, or Calder\'on identity,
\begin{equation*}\label{eq:plemelj}
	S_\Sigma K^*_\Sigma=K_\Sigma S_\Sigma ,
\end{equation*}
holds as an identity of bounded operators on $W^{s,p}(\Sigma)$. Together with \eqref{eq:coercivity} it makes $\langle\varphi,\psi\rangle_*:=\langle\varphi,-S_\Sigma\psi\rangle$ an inner product on $H^{-1/2}_0(\Sigma)$, equivalent to the $H^{-1/2}$ one, for which $K^*_\Sigma$ is self-adjoint. Hence the spectrum of $K^*_\Sigma$ on the mean-zero energy space is real and contained in $[-\frac12,\frac12)$.

\paragraph{The endpoints $\pm\frac12$.}
For a connected component $\Omega'$ of $\Omega$ (which could be $\Omega$ itself), let $\chi_{\partial\Omega'}$
denote the indicator of the whole of $\partial\Omega'$, outer and inner
boundaries alike. Whatever the boundary regularity,
\begin{equation*}
	K_\Sigma \chi_{\partial\Omega'} = \tfrac{1}{2}\chi_{\partial\Omega'}
	\ \text{ on } \partial\Omega',
	\qquad
	K_\Sigma \chi_{\partial\Omega'} = 0
	\ \text{ on } \Sigma\setminus\partial\Omega' .
\end{equation*}
Hence $\lambda = \frac{1}{2}$ is an eigenvalue of $K_\Sigma$ of multiplicity
the number of connected components of $\Omega$, and the same count holds for
$K^*_\Sigma$, whose eigenfunctions are the equilibrium densities of those
components. These carry non-zero charge. At the other endpoint,
$\lambda = -\frac{1}{2}$ is an eigenvalue of $K^*_\Sigma$ of multiplicity the
number of holes of $\Omega$, the eigenfunctions being the densities whose
single layer potential is constant on each hole and vanishes at infinity, and
these have vanishing average on the boundary of every component. The endpoint
$-\frac{1}{2}$ therefore persists on the mean-zero space, while $\frac{1}{2}$
does not.

\paragraph{Spectral properties.}
For $C^\infty$ boundaries $K^*_\Sigma$ is compact on $W^{s,p}(\Sigma)$ by \eqref{eq:smoothing}, so its spectrum is discrete and accumulates only at $0$. Bootstrapping \eqref{eq:smoothing} shows that eigenfunctions associated with non-zero eigenvalues are smooth, whence the non-zero spectrum and the eigenspaces are the same in every $W^{s,p}(\Sigma)$. Combined with the previous two paragraphs, they are real and lie in $[-\frac12,\frac12)$ on the mean-zero space. This applies verbatim to every truncated chain $\Gamma^{(N)}$. When the boundary is merely Lipschitz, or has infinitely many components, compactness is lost and a non-trivial essential spectrum appears.

\section{Layer potentials on the infinite chain}\label{sec:layer potentials infinite chain}
Because the boundary $\Gamma = \bigcup_{j=0}^\infty \Gamma_j$ consists of infinitely many disconnected domains, the global layer potentials act on $\varphi$ through an infinite sum of local self-interactions and cross-component interactions. Therefore, for $x \in \Gamma$, we define the single layer potential $S_\Gamma$, the Neumann--Poincar\'e operator 
$K^*_\Gamma$, and the double layer potential $K_\Gamma$ by the infinite series
$$
S_\Gamma \varphi(x) := \sum_{j=0}^\infty \int_{\Gamma_j} G(x, y) \varphi(y) \, d\sigma(y),
$$
$$
K^*_\Gamma \varphi(x) := \sum_{j=0}^\infty \int_{\Gamma_j} \frac{\partial G}{\partial \nu_x}(x, y) \varphi(y) \, d\sigma(y),
$$
$$
K_\Gamma \varphi(x) := \sum_{j=0}^\infty \int_{\Gamma_j} \frac{\partial G}{\partial \nu_y}(x, y) \varphi(y) \, d\sigma(y),
$$
Note that because the reference boundary $\Gamma_0$ is smooth, we do not need to read the self-interaction terms (when $x \in \Gamma_k$) in a principal value sense.

\subsection{Layer potentials on the mean-zero subspace}\label{subsec: layer potentils in mean zero}

Our interest is on the restriction of these operators to the componentwise mean-zero
subspace
\[
W^{s,p}_0(\Gamma)
:=
\left\{
\varphi\in W^{s,p}(\Gamma)
:
\int_{\Gamma_j}\varphi\,d\sigma=0
\quad\text{for all }j\ge0
\right\}.
\]

\paragraph{Why focus on the mean-zero subspace?}
This subspace is invariant under the Neumann--Poincar\'e operator $K^*_\Gamma$. To see this, let $\varphi \in W^{s,p}_0(\Gamma)$ and let $\chi_{\Gamma_j}$ be the indicator function of the $j$-th component. By duality and the classical identity $K_\Gamma\chi_{\Gamma_j}=\frac12\chi_{\Gamma_j}$, we have
\[
\bigl\langle K^*_\Gamma\varphi,\chi_{\Gamma_j}\bigr\rangle
= \bigl\langle\varphi, K_\Gamma\chi_{\Gamma_j}\bigr\rangle
= \frac12\bigl\langle\varphi,\chi_{\Gamma_j}\bigr\rangle = 0,
\qquad \text{for all } j\ge0 .
\]
Thus, $K^*_\Gamma\varphi$ also has zero mean on every component.

To understand the remainder of the space, let $Y:=\{\varphi:\ S_\Gamma\varphi \text{ is constant on each } \Omega_j\}$, which is precisely the eigenspace associated with the eigenvalue $\lambda = 1/2$, meaning $K^*_\Gamma$ acts identically as $\frac12 I$ on this subspace. Assuming the decomposition $W^{s,p}(\Gamma)=W^{s,p}_0(\Gamma)\oplus Y$ (not proved in general), the operator admits a block-diagonal representation:
\[
K^*_\Gamma=
\begin{pmatrix}
	K^*_\Gamma|_{W^{s,p}_0(\Gamma)} & 0\\[2pt]
	0 & \frac12 I
\end{pmatrix}.
\]
Consequently, the spectrum splits entirely as
\[
\sigma\bigl(K^*_\Gamma\bigr)
=\sigma\bigl(K^*_\Gamma|_{W^{s,p}_0(\Gamma)}\bigr)\cup\Bigl\{\frac12\Bigr\} .
\]
The second block simply contributes the eigenvalue $1/2$, which remains exactly the same for every chain, completely independent of the geometric scaling $r$, the shape of the reference cell $\Gamma_0$, and the functional parameters $s$ and $p$. Therefore, all the non-trivial spectral information encoded by the self-similar geometry is concentrated within the mean-zero subspace.

\paragraph{The projected operators.} The above motivates introducing the
\emph{projected layer potentials}, which will be the main focus of our study:
\[
\widetilde S_\Gamma:=P_0S_\Gamma P_0,
\qquad
\widetilde K_\Gamma^*:=P_0K_\Gamma^*P_0,
\qquad
\widetilde K_\Gamma:=P_0K_\Gamma P_0,
\]
where, writing $\varphi_j:=\varphi|_{\Gamma_j}$,
\[
P_0\varphi:=\sum_{j=0}^{\infty}\bigl(P_{0,j}\varphi_j\bigr)\chi_{\Gamma_j},
\qquad
P_{0,j}\varphi_j:=\varphi_j-\frac{1}{|\Gamma_j|}\int_{\Gamma_j}\varphi_j\,d\sigma,
\qquad
\varphi\in W^{s,p}(\Gamma),
\]
and $|\Gamma_j|$ denotes the surface measure of $\Gamma_j$, the integral being
understood as the duality pairing against $\chi_{\Gamma_j}$ when $s<0$.

The local projections $P_{0,j}:W^{s,p}(\Gamma_j)\to W^{s,p}_0(\Gamma_j)$ are
uniformly bounded in $j$. The global projection $P_0$ is also bounded on
$W^{s,p}(\Gamma)$ for every $p\in(1,\infty)$ and $s\in(-1,1)$; this follows from
Lemma~\ref{lem:sequence-isomorphism}, and is stated explicitly in
Corollary~\ref{cor:P0-bounded}. Since each $P_{0,j}$ is the
$L^2(\Gamma_j)$-orthogonal projection onto the mean-zero densities and $P_0$ acts
componentwise, we have that $P_0$ is symmetric for the $L^2(\Gamma)$ pairing; in particular,
when $s=0$ and $p=2$ it is the orthogonal projection onto $L^2_0(\Gamma)$.

\begin{remark}
	Since
	\[
	\widetilde K_\Gamma^*\varphi
	=
	K_\Gamma^*\varphi,
	\qquad
	\varphi\in W^{s,p}_0(\Gamma),
	\]
	$\widetilde K_\Gamma^*$ may equivalently be regarded as the
	restriction of $K_\Gamma^*$ to $W^{s,p}_0(\Gamma)$.
\end{remark}

\subsection{Some properties of the projected operators}

The following properties will be used in later sections and may be skipped on a
first reading. Their proofs rely only on Lemma~\ref{lem:sequence-isomorphism} and Lemma~\ref{lem:scaling-isomorphism}, stated in Section~\ref{sec:sequence spaces},
and are independent of this section.


\begin{lemma}[Projected symmetrisation]\label{lem:projected-symmetrization}
	Let $p \in (1, \infty)$ and $s \in (-1, 1)$ be such that
	$\widetilde{S}_{\Gamma}$, $\widetilde{K}^*_{\Gamma}$ and
	$\widetilde{K}_{\Gamma}$ are bounded on $W^{s,p}_0(\Gamma)$. Then
	$$
	\widetilde{S}_{\Gamma} \widetilde{K}^*_{\Gamma}
	= \widetilde{K}_{\Gamma} \widetilde{S}_{\Gamma}.
	$$
\end{lemma}

\begin{proof}
	Let $Q_N$ denote the projection onto the first $N$ components of $\Gamma$.
	Since $Q_N$ and $P_0$ both act componentwise they commute, and since the
	kernels of the layer potentials do not see the discarded components,
	\[
	Q_N\tilde S_\Gamma Q_N=\tilde S_{\Gamma^{(N)}},\qquad
	Q_N\tilde K^*_\Gamma Q_N=\tilde K^*_{\Gamma^{(N)}},\qquad
	Q_N\tilde K_\Gamma Q_N=\tilde K_{\Gamma^{(N)}} .
	\]
	
	The boundary $\Gamma^{(N)}$ is a finite union of disjoint smooth components,
	so the classical symmetrisation principle applies to it. Write
	$P_0^c:=I-P_0$, whose range consists of the componentwise constants. We have that
	\[
	P_0K_{\Gamma^{(N)}}P_0^c=\tfrac12P_0P_0^c=0,
	\qquad
	P_0^cK^*_{\Gamma^{(N)}}P_0=0;
	\]
	the first because $K_{\Gamma^{(N)}}$ acts as $\tfrac12 I$ on the
	componentwise constants, the second because $K^*_{\Gamma^{(N)}}$ preserves
	the componentwise mean-zero condition, as shown above. Using the above and $P_0 = I-P_0^c$ we get
	\[
	\tilde S_{\Gamma^{(N)}}\tilde K^*_{\Gamma^{(N)}}
	=P_0S_{\Gamma^{(N)}}P_0K^*_{\Gamma^{(N)}}P_0
	=P_0S_{\Gamma^{(N)}}K^*_{\Gamma^{(N)}}P_0 ,
	\]
	\[
	\tilde K_{\Gamma^{(N)}}\tilde S_{\Gamma^{(N)}}
	=P_0K_{\Gamma^{(N)}}P_0S_{\Gamma^{(N)}}P_0
	=P_0K_{\Gamma^{(N)}}S_{\Gamma^{(N)}}P_0.
	\]
	The two right-hand sides agree by the classical identity. Hence
	$\tilde S_{\Gamma^{(N)}}\tilde K^*_{\Gamma^{(N)}}
	=\tilde K_{\Gamma^{(N)}}\tilde S_{\Gamma^{(N)}}$ for every $N$.
	
	By Lemma~\ref{lem:sequence-isomorphism} the space
	$W^{s,p}_0(\Gamma)$ is isomorphic to an $\ell^p$ direct sum with $p<\infty$,
	so the $Q_N$ are uniformly bounded and $Q_N\to I$ strongly. The compressions
	above therefore converge strongly to the corresponding infinite-chain
	operators and are uniformly bounded, the latter by the hypothesis of the
	lemma, so their products converge strongly as well. Passing to the limit
	gives $\tilde S_\Gamma\tilde K^*_\Gamma=\tilde K_\Gamma\tilde S_\Gamma$.
\end{proof}

\begin{lemma}[Coercivity on the infinite chain]\label{lem:uniform-coercivity}
	There exists a constant $c>0$ such that
	$$
	\big\langle \varphi,-\widetilde{S}_{\Gamma}\varphi \big\rangle \ge c\,\|\varphi\|_{H^{-1/2}(\Gamma)}^2, \qquad \varphi\in H^{-1/2}_0(\Gamma).
	$$
\end{lemma}

\begin{proof}
	Each truncation $\Gamma^{(N)}$ is a Lipschitz boundary, so the classical
	coercivity of the single layer on the mean-zero space gives a constant
	$c_N>0$ with
	\[
	\bigl\langle\varphi,-\tilde S_\Gamma\varphi\bigr\rangle
	\ \ge\ c_N\|\varphi\|^2_{H^{-1/2}(\Gamma)}
	\]
	for every $\varphi$ supported on the first $N$ components, the operator and
	the norm being unchanged by the truncation on such densities. Since
	$H^{-1/2}_0(\Gamma)$ is equivalent to an $\ell^2$ direct sum over the components by
	Lemma~\ref{lem:sequence-isomorphism}, these densities are dense, and both sides are
	continuous, the left-hand side because $\tilde S_\Gamma$ is bounded by
	Theorem~\ref{thm:boundedness} (which is independent of this result). The lemma therefore holds with
	$c=\inf_N c_N$, and it remains only to check that this infimum is positive.
	
	To do that,  let $\varphi\in H^{-1/2}_0(\Gamma)$ be supported on finitely many components and
		put $u:=S_\Gamma\varphi$. Since the layer potential does not see the components
		where $\varphi$ vanishes, the energy identity \eqref{eq:energy} on the smooth
		boundary $\Gamma^{(N)}$ gives
		$\langle\varphi,-\tilde S_\Gamma\varphi\rangle=\int_{\mathbb R^d}|\nabla u|^2$,
		and $\varphi_j$ is the jump $[\partial_\nu u]$ across $\Gamma_j$.
		
	Fix a collar $U_0$ of $\Gamma_0$ whose dilates $U_j:=r^jU_0$ are pairwise
disjoint. We claim that
\begin{equation}\label{eq:collar}
	\|[\partial_\nu v]\|_{H^{-1/2}(\Gamma_0)}\le C_0\|\nabla v\|_{L^2(U_0)}
\end{equation}
for every $v$ harmonic in $U_0\setminus\Gamma_0$ with $\nabla v\in L^2(U_0)$.
Indeed, given $\xi\in H^{1/2}(\Gamma_0)$, the trace theorem provides
$\chi\in H^1_0(U_0)$ with $\chi|_{\Gamma_0}=\xi$ and
$\|\nabla\chi\|_{L^2(U_0)}\lesssim\|\xi\|_{H^{1/2}(\Gamma_0)}$. Green's identity
on each side of $\Gamma_0$ gives
$\langle[\partial_\nu v],\xi\rangle=\int_{U_0}\nabla v\cdot\nabla\chi$, the
boundary terms on $\partial U_0$ vanishing since $\chi$ has compact support, and
taking the supremum over $\xi$ with $\|\xi\|_{H^{1/2}(\Gamma_0)}\le1$ gives
\eqref{eq:collar}.
		
		Apply this to $w(\tilde x):=u(r^j\tilde x)$, harmonic in $U_0\setminus\Gamma_0$
		because $r^jU_0=U_j$. 	On the left $[\partial_\nu w](\tilde x)=r^j\varphi_j(r^j\tilde x)$, and
		Lemma~\ref{lem:scaling-isomorphism} at $(s,p)=(-1/2,2)$, for which $\alpha=d/2$,
		gives
		$\|\varphi_j\|_{H^{-1/2}(\Gamma_j)}\asymp
		r^{jd/2}\|\varphi_j(r^j\,\cdot\,)\|_{H^{-1/2}(\Gamma_0)}$. On the right
		$\|\nabla w\|_{L^2(U_0)}=r^{j(1-d/2)}\|\nabla u\|_{L^2(U_j)}$. Both sides carry
		the factor $r^{j(1-d/2)}$, which cancels, so
		\[
		\|\varphi_j\|_{H^{-1/2}(\Gamma_j)}\le C\,\|\nabla u\|_{L^2(U_j)}
		\]
		with the same constant on every component.
		
		Squaring and summing over $j$, the collars being disjoint,
		\[
		\sum_{j\ge0}\|\varphi_j\|^2_{H^{-1/2}(\Gamma_j)}
		\le C^2\int_{\mathbb R^d}|\nabla u|^2
		= C^2\bigl\langle\varphi,-\tilde S_\Gamma\varphi\bigr\rangle ,
		\]
		and Lemma~\ref{lem:sequence-isomorphism} identifies the left-hand side with
		$\|\varphi\|^2_{H^{-1/2}(\Gamma)}$ up to a constant. Such densities are dense in
		$H^{-1/2}_0(\Gamma)$ by that same lemma, and the quadratic form is continuous
		because $\tilde S_\Gamma$ is bounded by Theorem~\ref{thm:boundedness}, so the
		estimate passes to the whole space.
\end{proof}

\section{Equivalence of $W_0^{s,p}(\Gamma)$ with a sequence space}\label{sec:sequence spaces}
In Section~\ref{sec:block-Toeplitz} we aim to represent $\tilde K_\Gamma^*: W_0^{s,p}(\Gamma)\to W_0^{s,p}(\Gamma)$ as a block-Toeplitz operator acting on the space
$$
\ell^p\big(\mathbb{N}, W_0^{s,p}(\Gamma_0)\big) := \left\{ \begin{pmatrix} u_0 \\ u_1 \\ \vdots \end{pmatrix} \in \bigoplus_{j=0}^{\infty} W_0^{s,p}(\Gamma_0) \;:\; \sum_{j=0}^\infty \|u_j\|^p_{W^{s,p}(\Gamma_0)} < \infty \right\},
$$
for suitable $p \in (1, \infty)$ and $s \in (-1,1)$. To that end, in this section we show that $W_0^{s,p}(\Gamma)$ is isomorphic to $\ell^p\big(\mathbb{N}, W_0^{s,p}(\Gamma_0)\big)$, through a composition of two maps: 
\begin{enumerate}
	\item  a restriction map $\mathcal T$ which isomorphically maps $W_0^{s,p}(\Gamma)$ to the $\ell^p$ sum of the local component norms 
	$$
	V_0^{s,p}(\Gamma) := \left\{ \begin{pmatrix} \varphi_0 \\ \varphi_1 \\ \vdots \end{pmatrix} \in \bigoplus_{j=0}^{\infty} W_0^{s,p}(\Gamma_j) \;:\; \sum_{j=0}^\infty \|\varphi_j\|^p_{W^{s,p}(\Gamma_j)} < \infty \right\};
	$$
	\item a composition of a pullback map $\mathcal P$ and a scaling map $D^\alpha$, which together isomorphically map $V_0^{s,p}(\Gamma)$ to $\ell^p\big(\mathbb{N}, W_0^{s,p}(\Gamma_0)\big)$. Here $\alpha$ is a parameter that depends on $s$ and $p$ that will be properly introduced below.
\end{enumerate}
Together, they define the isomorphic map $\mathcal M_\alpha:W_0^{s,p}(\Gamma)\to \ell^p\big(\mathbb{N}, W_0^{s,p}(\Gamma_0)\big)$, defined as
$$
\mathcal M_\alpha : = D^\alpha \mathcal{P}\mathcal T,
$$
which will later allows us to construct the block-Toeplitz representation of the projected layer potential operators.

\subsection{The restriction map}
Let $\mathcal T$ be the map that identifies a global density
$\varphi\in W^{s,p}_0(\Gamma)$, $p\in(1,\infty)$ and $s\in(-1,1)$, with the
vertical block vector of its localised densities $\varphi_j$:
\[
\mathcal T\varphi=
\begin{pmatrix}\varphi_0\\ \varphi_1\\ \vdots\end{pmatrix}
\in\bigoplus_{j=0}^{\infty}W^{s,p}_0(\Gamma_j).
\]
For $s\ge 0$ this is pointwise restriction, for $s<0$ it is defined by duality.

\begin{lemma}\label{lem:sequence-isomorphism}
	Let $p\in(1,\infty)$ and $s\in(-1,1)$. The map
	$
	\mathcal T:W^{s,p}_0(\Gamma)\longrightarrow V^{s,p}_0(\Gamma),
	$
	is an isomorphism. 
\end{lemma}

\begin{proof}
	The proof divides into three regularity regimes, the negative range being
	obtained from the positive one by duality.
	
	\noindent \textbf{Case $s\in(0,1)$:} Writing $\varphi_j:=\varphi|_{\Gamma_j}$, the global norm decomposes into local contributions and cross-component interactions:
	\begin{equation}
		\|\varphi\|_{W^{s,p}(\Gamma)}^p = \sum_{j\ge0}\|\varphi_j\|_{W^{s,p}(\Gamma_j)}^p + \sum_{j\neq k} I_{j,k},
		\label{eq:global-decomposition}
	\end{equation}
	where
	\[
	I_{j,k} := \int_{\Gamma_j}\int_{\Gamma_k} \frac{|\varphi_j(x)-\varphi_k(y)|^p}{|x-y|^{d-1+sp}} \,d\sigma(y)\,d\sigma(x).
	\]
	Since the cross-interaction terms are non-negative, dropping them immediately yields the lower bound:
	\begin{equation}
		\sum_{j\ge0}\|\varphi_j\|_{W^{s,p}(\Gamma_j)}^p \le \|\varphi\|_{W^{s,p}(\Gamma)}^p.
		\label{eq:local-below-global}
	\end{equation}
	
	To prove the reverse inequality, we bound the cross-interactions. Using $|a-b|^p \le 2^{p-1}(|a|^p+|b|^p)$ and the symmetry of the double sum, we obtain:
	\[
	\sum_{j\neq k}I_{j,k} \lesssim \sum_{j\ge0} \int_{\Gamma_j}|\varphi_j(x)|^p \left( \sum_{k\neq j} \int_{\Gamma_k} \frac{d\sigma(y)}{|x-y|^{d-1+sp}} \right) d\sigma(x).
	\]
	By the self-similarity and strict separation of the components, we have $\operatorname{dist}(\Gamma_j,\Gamma_k) \gtrsim r^{\min\{j,k\}}$ and $|\Gamma_k|=r^{k(d-1)}|\Gamma_0|$. Splitting the sum into $k<j$ and $k>j$ and summing the resulting geometric series bounds the inner bracket by $\mathcal{O}(r^{-jsp})$. Hence,
	\[
	\sum_{j\neq k}I_{j,k} \lesssim \sum_{j\ge0} r^{-jsp}\|\varphi_j\|_{L^p(\Gamma_j)}^p.
	\]
	Crucially, because each $\varphi_j$ has zero mean on $\Gamma_j$, the fractional Poincar\'e--Friedrichs inequality on $\Gamma_0$ \cite{DiNezzaPalatucciValdinoci2012}, transported to $\Gamma_j$ via the dilation $x=r^j\widetilde x$, yields $\|\varphi_j\|_{L^p(\Gamma_j)} \lesssim r^{js}|\varphi_j|_{W^{s,p}(\Gamma_j)}$. Consequently, the cross-interactions are controlled entirely by the local semi-norms:
	\[
	\sum_{j\neq k}I_{j,k} \lesssim \sum_{j\ge0} \|\varphi_j\|_{W^{s,p}(\Gamma_j)}^p.
	\]
	Combining this bound with \eqref{eq:global-decomposition} and \eqref{eq:local-below-global} establishes the norm equivalence:
	\begin{equation}
		\|\varphi\|_{W^{s,p}(\Gamma)} \asymp \left( \sum_{j\ge0} \|\varphi_j\|_{W^{s,p}(\Gamma_j)}^p \right)^{1/p}.
		\label{eq:sequence-equivalence-positive}
	\end{equation}
	Conversely, if $(\varphi_j)_{j\ge0}\in V^{s,p}_0(\Gamma)$, applying this same estimate to the sequence tails shows that the partial sums $\sum_{j<N}\varphi_j\chi_{\Gamma_j}$ (extended by zero) form a Cauchy sequence in $W^{s,p}_0(\Gamma)$. Thus, they converge to a unique element $\varphi\in W^{s,p}_0(\Gamma)$ satisfying $\mathcal T\varphi=(\varphi_j)_{j\ge0}$, proving that $\mathcal T$ is surjective.
	
	We record a consequence used below. Applying
	\eqref{eq:sequence-equivalence-positive} to the componentwise mean-zero
	density $P_0\varphi$, then the uniform boundedness of the local projections
	$P_{0,j}$, and finally \eqref{eq:local-below-global}, which holds for every
	$\varphi\in W^{s,p}(\Gamma)$, gives
	\begin{equation}
		\|P_0\varphi\|_{W^{s,p}(\Gamma)}
		\lesssim \Bigl( \sum_{j\ge0}\|P_{0,j}\varphi_j\|_{W^{s,p}(\Gamma_j)}^p \Bigr)^{1/p}
		\lesssim \Bigl( \sum_{j\ge0}\|\varphi_j\|_{W^{s,p}(\Gamma_j)}^p \Bigr)^{1/p}
		\lesssim \|\varphi\|_{W^{s,p}(\Gamma)}.
		\label{eq:P0-positive}
	\end{equation}
	
	\noindent \textbf{Case $s=0$:} The claim follows trivially from the disjointness of the components, yielding the exact identity $\|\varphi\|^p_{L^p(\Gamma)}=\sum_{j\ge0}\|\varphi_j\|^p_{L^p(\Gamma_j)}$, and \eqref{eq:P0-positive} holds for the same reason.
	
	\noindent \textbf{Case $s\in(-1,0)$:} Let $t=-s\in(0,1)$, and let $p'$ be the H\"older conjugate exponent, so that $W^{s,p}(\Gamma)=\bigl(W^{t,p'}(\Gamma)\bigr)^*$ and $W^{s,p}(\Gamma_j)=\bigl(W^{t,p'}(\Gamma_j)\bigr)^*$, where $(\,\cdot\,)^*$ denotes the $L^2$ conjugate. 
	
	We first note that
	$\bigl(W^{t,p'}_0(\Gamma)\bigr)^{*}=W^{s,p}_0(\Gamma)$ and
	$\bigl(W^{t,p'}_0(\Gamma_j)\bigr)^{*}=W^{s,p}_0(\Gamma_j)$.  Indeed, by \eqref{eq:P0-positive}, $P_0$ is a bounded projection so
	$W^{t,p'}(\Gamma)=P_0W^{t,p'}(\Gamma)\oplus\ker P_0$. Therefore,
	\[
	\bigl(W^{t,p'}_0(\Gamma)\bigr)^{*}
	=\bigl(P_0W^{t,p'}(\Gamma)\bigr)^{*}
	=P_0^*\bigl(W^{t,p'}(\Gamma)\bigr)^{*}
	=P_0W^{s,p}(\Gamma)
	=W^{s,p}_0(\Gamma),
	\]
	with the dual
	norm equivalent to that of $W^{s,p}(\Gamma)$. The same argument on a single
	component gives $\bigl(W^{t,p'}_0(\Gamma_j)\bigr)^{*}=W^{s,p}_0(\Gamma_j)$,
	uniformly in $j$.
	
Next, by the positive-regularity case at $(t,p')$, the map
\[
\mathcal T_{t,p'}:W^{t,p'}_0(\Gamma)\longrightarrow
\Bigl(\bigoplus_{j\ge0}W^{t,p'}_0(\Gamma_j)\Bigr)_{\ell^{p'}}
\]
is an isomorphism. Since $p'<\infty$, the dual of an $\ell^{p'}$
direct sum is the $\ell^p$ direct sum of the duals, we have that $\mathcal T_{t,p'}^{*}$
maps $V^{s,p}_0(\Gamma)$ onto $W^{s,p}_0(\Gamma)$.

We want to identify $\mathcal T$ with $\bigl(\mathcal T_{t,p'}^{*}\bigr)^{-1}$.
To that end, let $\varphi\in W^{s,p}_0(\Gamma)$ and $\psi\in W^{t,p'}_0(\Gamma)$,
and write $\psi^{(N)}:=\sum_{j<N}\psi_j\chi_{\Gamma_j}$. By the
positive-regularity case $\psi^{(N)}\to\psi$ in $W^{t,p'}_0(\Gamma)$, so
\[
\langle\varphi,\psi\rangle
=\lim_{N\to\infty}\bigl\langle\varphi,\psi^{(N)}\bigr\rangle
=\lim_{N\to\infty}\sum_{j<N}\bigl\langle\varphi,\psi_j\chi_{\Gamma_j}\bigr\rangle
=\sum_{j\ge0}\langle\varphi_j,\psi_j\rangle .
\]
The right-hand side is the duality pairing between $V^{s,p}_0(\Gamma)$ and the
$\ell^{p'}$ direct sum, evaluated at $\mathcal T\varphi$ and
$\mathcal T_{t,p'}\psi$, so $\langle\varphi,\psi\rangle=\langle\mathcal
T\varphi,\mathcal T_{t,p'}\psi\rangle$ for all $\psi$. This identifies
$\mathcal T$ with $\bigl(\mathcal T_{t,p'}^{*}\bigr)^{-1}$. Since the adjoint of
an isomorphism between Banach spaces is again one, $\mathcal T$ is a
isomorphism for negative Sobolev indices as well.
\end{proof}

Before carrying on, we complete the proof of the boundedness of the
mean-zero projection $P_0$, introduced in
Section~\ref{subsec: layer potentils in mean zero}.

\begin{corollary}\label{cor:P0-bounded}
	For every $p\in(1,\infty)$ and $s\in(-1,1)$, the componentwise
	mean-zero projection
	\[
	P_0:W^{s,p}(\Gamma)\longrightarrow W^{s,p}_0(\Gamma)
	\]
	is bounded.
\end{corollary}

\begin{proof}
	For $s\in[0,1)$ this is \eqref{eq:P0-positive}. For $s\in(-1,0)$ set
	$t=-s$. We have that $P_0$ is symmetric
	for the $L^2(\Gamma)$ pairing. Hence, for $\psi\in W^{t,p'}(\Gamma)$,
	\[
	|\langle P_0\varphi,\psi\rangle|
	=|\langle \varphi,P_0\psi\rangle|
	\lesssim \|\varphi\|_{W^{s,p}(\Gamma)}\|\psi\|_{W^{t,p'}(\Gamma)},
	\]
	by the case already proved. Taking the supremum over $\psi$ gives the claim.
\end{proof}

\subsection{The pullback and scaling maps}
Let $\eta_j : \Gamma_0 \to \Gamma_j$ denote the bijective similarity mapping that constructs the $j$-th sub-component, satisfying $\eta_j(\tilde{x}) = r^j \tilde{x}$ for all $\tilde{x} \in \Gamma_0$, where $r \in (0,1)$ is the geometric scaling factor. Because the surface measure scales uniformly, pulling back a mean-zero density preserves the zero-mean property. 

We define the canonical spatial pullback operator $\mathcal{P} : V_0^{s,p}(\Gamma) \to \ell^p\big(\mathbb{N}, W_0^{s,p}(\Gamma_0)\big)$, acting component-wise as:
$$ (\mathcal{P}\varphi)_j := \varphi_j \circ \eta_j. $$

While $\mathcal{P}$ standardizes the domain, it does not account for the geometric scaling of the surface measure or the regularity index of the functional space. To correct this, we introduce a family of diagonal scaling operators $D^\alpha$, parametrised by an exponent $\alpha \in \mathbb{R}$, where $D^\alpha$ acts on the sequence space as the infinite diagonal block matrix:
$$ D^\alpha = \begin{pmatrix} I & 0 & 0 & \cdots \\ 0 & r^\alpha I & 0 & \cdots \\ 0 & 0 & r^{2\alpha} I & \cdots \\ \vdots & \vdots & \vdots & \ddots \end{pmatrix}. $$

The composite scaled pullback mapping is then $D^\alpha \mathcal{P}$. The specific value of $\alpha$ entirely dictates the functional space for which this mapping acts as an isomorphism.

\begin{lemma}\label{lem:scaling-isomorphism}
	For $p \in (1, \infty)$ and $s \in (-1,1)$, let 
	$$
	\alpha := \frac{d-1}{p} - s.
	$$ 
	The mapping $D^\alpha \mathcal{P}$ defines an isomorphism from $V_0^{s,p}(\Gamma)$ onto $\ell^p\big(\mathbb{N}, W_0^{s,p}(\Gamma_0)\big)$.
\end{lemma}

\begin{proof}
	The norm on $V_0^{s,p}(\Gamma)$ is the $\ell^p$ sum of the local norms
	$\|\varphi_j\|_{W^{s,p}(\Gamma_j)}$, and $D^\alpha\mathcal P$ acts
	componentwise by $\varphi_j\mapsto r^{j\alpha}u_j$ with
	$u_j:=\varphi_j\circ\eta_j$. Each $\eta_j$ is a bijection and $D^\alpha$ is
	diagonal with non-zero entries, so the map is a componentwise bijection and
	it suffices to prove
	\begin{equation}
		\|\varphi_j\|_{W^{s,p}(\Gamma_j)}
		\asymp r^{j\alpha}\|u_j\|_{W^{s,p}(\Gamma_0)},
		\label{eq:component-scaling}
	\end{equation}
	with constants independent of $j$.
	
	\noindent\textbf{Case $s\in[0,1)$:} The substitution $x=\eta_j(\widetilde x)$
gives $d\sigma(x)=r^{j(d-1)}d\sigma_0(\widetilde x)$ and
$|x-y|=r^{j}|\widetilde x-\widetilde y|$, whence
\[
\|\varphi_j\|^p_{L^p(\Gamma_j)}=r^{j(d-1)}\|u_j\|^p_{L^p(\Gamma_0)} .
\]
For $s=0$ this is \eqref{eq:component-scaling} with $\alpha=\frac{d-1}{p}$, and
$D^\alpha\mathcal P$ is an exact isometry. For $s\in(0,1)$ the norm carries in
addition the Gagliardo seminorm (denoted $|\,\cdot\,|_{W^{s,p}}$, defined by the double
integral in \eqref{eq:Gagliadrdo}). There, the two surface measures
contribute $r^{2j(d-1)}$ and the kernel $r^{-j(d-1+sp)}$, so that
\[
|\varphi_j|^p_{W^{s,p}(\Gamma_j)}=r^{j(d-1-sp)}|u_j|^p_{W^{s,p}(\Gamma_0)} .
\]
The two summands therefore scale with different powers of $r^j$, and since
$r^{j(d-1)}\le r^{j(d-1-sp)}$,
\[
r^{j(d-1-sp)}|u_j|^p_{W^{s,p}(\Gamma_0)}
\;\le\;\|\varphi_j\|^p_{W^{s,p}(\Gamma_j)}
\;\le\; r^{j(d-1-sp)}\Bigl(\|u_j\|^p_{L^p(\Gamma_0)}+|u_j|^p_{W^{s,p}(\Gamma_0)}\Bigr).
\]
Since $u_j$ has zero mean, the fractional Poincar\'e--Friedrichs inequality on
$\Gamma_0$ \cite{DiNezzaPalatucciValdinoci2012} bounds the bracket by
$|u_j|^p_{W^{s,p}(\Gamma_0)}$ up to a constant, and both outer terms are
$\asymp r^{j(d-1-sp)}\|u_j\|^p_{W^{s,p}(\Gamma_0)}$. Taking $p$-th roots gives
\eqref{eq:component-scaling} with $\alpha=\frac{d-1}{p}-s$. The inequality is
used only on $\Gamma_0$, so the constants do not depend on $j$.
	
	\noindent\textbf{Case $s\in(-1,0)$:} We argue by duality. Let $t=-s\in(0,1)$, let $p'$ be the conjugate exponent, and let $\alpha':=\frac{d-1}{p'}-t$ be the exponent associated with $(t,p')$, so that
	\[
	\alpha+\alpha'=(d-1)\Bigl(\tfrac1p+\tfrac1{p'}\Bigr)-s-t=d-1 .
	\]
	On each component, $W^{s,p}_0(\Gamma_j)$ is the dual of $W^{t,p'}_0(\Gamma_j)$, the mean being a bounded functional on $W^{t,p'}(\Gamma_j)$. Writing $v_j:=\psi_j\circ\eta_j$ we have $\langle\varphi_j,\psi_j\rangle_{L^2(\Gamma_j)}=r^{j(d-1)}\langle u_j,v_j\rangle_{L^2(\Gamma_0)}$, while $\|\psi_j\|_{W^{t,p'}(\Gamma_j)}\asymp r^{j\alpha'}\|v_j\|_{W^{t,p'}(\Gamma_0)}$ by the case already proved. Taking the supremum over $\psi_j\in W^{t,p'}_0(\Gamma_j)$,
	\[
	\|\varphi_j\|_{W^{s,p}(\Gamma_j)}
	\asymp r^{j(d-1-\alpha')}\|u_j\|_{W^{s,p}(\Gamma_0)}
	= r^{j\alpha}\|u_j\|_{W^{s,p}(\Gamma_0)},
	\]
	which is \eqref{eq:component-scaling}.
\end{proof}

\section{Block-Toeplitz structure}\label{sec:block-Toeplitz}
The goal of this section is to derive the block-Toeplitz representation of the projected layer potentials $\widetilde{S}_{\Gamma}$,  $\widetilde{K}^*_{\Gamma}$, and $\widetilde{K}_{\Gamma}$. Doing so will also reveal for which values of $s$ and $p$ these operators are bounded on $W_0^{s,p}(\Gamma)$. 

Consider the canonical isomorphism $\mathcal M_\alpha:W_0^{s,p}(\Gamma)\to \ell^p\big(\mathbb{N}, W_0^{s,p}(\Gamma_0)\big)$ established in Section~\ref{sec:sequence spaces}. We introduce the following block-matrix scaled operators, acting on $\ell^p\big(\mathbb{N}, W_0^{s,p}(\Gamma_0)\big)$:
\begin{equation*}\label{eq: S Toeplitz}
\widetilde{\mathbb{S}}_\alpha
:=
\mathcal M_\alpha\widetilde{S}_\Gamma
\mathcal M_\alpha^{-1},
\end{equation*}
\begin{equation*}\label{eq: K^* Toeplitz}
\widetilde{\mathbb{K}}^*_\alpha
:=
\mathcal M_\alpha\widetilde{K}^*_\Gamma
\mathcal M_\alpha^{-1},
\end{equation*}
and
\begin{equation*}\label{eq: K Toeplitz}
\widetilde{\mathbb{K}}_\alpha
:=
\mathcal M_\alpha\widetilde{K}_\Gamma
\mathcal M_\alpha^{-1}.
\end{equation*}
Note that the above are block-Matrix operators that are similar to the projected operators $\widetilde{S}_{\Gamma}$,  $\widetilde{K}^*_{\Gamma}$, and $\widetilde{K}_{\Gamma}$. 

Recall that $\eta_j : \Gamma_0 \to \Gamma_j$, defined by $\eta_j(\tilde{x}) = r^j \tilde{x}$, is the bijective similarity mapping constructing the $j$-th component, and that $P_{0,j} : W^{s,p}(\Gamma_j) \to W_0^{s,p}(\Gamma_j)$ is the local zero-mean projection on the domain $\Gamma_j$. 

\subsection{Block-Toeplitz structure of $\widetilde{\mathbb{K}}^*_\alpha$ and $\widetilde{\mathbb{K}}_\alpha$}
We only detail the computations for $\widetilde{\mathbb{K}}^*_\alpha$. Similar results can be obtained for $\widetilde{\mathbb{K}}_\alpha$ by taking $L^2$ transpose of $\widetilde{\mathbb{K}}^*_\alpha$.

\begin{lemma}[Structure of the scaled NP operator]\label{lem: Toeplitz structure of NP operator}
	The scaled Neumann-Poincaré operator $\widetilde{\mathbb{K}}^*_\alpha$ admits the explicit block-Toeplitz representation:
	\begin{equation*}
		\widetilde{\mathbb{K}}^*_\alpha = \begin{pmatrix} \tilde{K}^*_0 & r^{-\alpha} \tilde{K}^*_{-1} & r^{-2\alpha} \tilde{K}^*_{-2} & \cdots \\ r^\alpha \tilde{K}^*_1 & \tilde{K}^*_0 & r^{-\alpha} \tilde{K}^*_{-1} & \cdots \\ r^{2\alpha} \tilde{K}^*_2 & r^\alpha \tilde{K}^*_1 & \tilde{K}^*_0 & \cdots \\ \vdots & \vdots & \vdots & \ddots \end{pmatrix}. 
	\end{equation*}
	The blocks are defined by $(\widetilde{\mathbb{K}}^*_\alpha)_{j,k} = r^{(j-k)\alpha} \tilde{K}^*_{j-k}$ for $j, k \ge 0$. The reference operator $\tilde{K}^*_m : W_0^{s,p}(\Gamma_0) \to W_0^{s,p}(\Gamma_0)$ is given by:
	\begin{equation*}
		(\tilde{K}^*_m \psi)(\tilde{x}) = P_{0,0} \left( \mathrm{} \int_{\Gamma_0} \frac{\partial G}{\partial \nu_{\tilde{x}}}(\eta_m(\tilde{x}), \tilde{y}) \psi(\tilde{y}) \, d\sigma_0(\tilde{y}) \right).
	\end{equation*}
\end{lemma}

\begin{proof}
	Recall that $\mathcal M_\alpha : = D^\alpha \mathcal{P}\mathcal T$, so 
	$$
	\widetilde{\mathbb{K}}^*_\alpha
	:=
	D^\alpha \mathcal{P}\mathcal T\widetilde{K}^*_\Gamma
	\mathcal \mathcal \mathcal T^{-1}\mathcal{P}^{-1} D^{-\alpha}.
	$$
	It follows directly from the definition of $\mathcal T$ that
	\[
	\tilde{\mathsf{K}}^*
	:=
	\mathcal{T}\widetilde{K}^*_\Gamma\mathcal{T}^{-1}
	=
	\begin{pmatrix}
		\tilde{\mathsf{K}}^*_{0,0} & \tilde{\mathsf{K}}^*_{0,1} & \cdots \\
		\tilde{\mathsf{K}}^*_{1,0} & \tilde{\mathsf{K}}^*_{1,1} & \cdots \\
		\vdots & \vdots & \ddots
	\end{pmatrix},
	\]
	where 
	\[
	(\tilde{\mathsf{K}}^*_{j,k}\varphi_k)(x)
	:=
	P_{0,j}\left(
	\int_{\Gamma_k}
	\frac{\partial G}{\partial \nu_{\cdot}}(\cdot,y)
	\varphi_k(y)\,d\sigma(y)
	\right)(x),
	\qquad x\in\Gamma_j.
	\]
	
    We next evaluate the purely spatial pullback of the $(j,k)$-th block. For $\psi \in W_0^{s,p}(\Gamma_0)$,
	\begin{align*}
		(\mathcal{P} \tilde{\mathsf{K}}^* \mathcal{P}^{-1})_{j,k} \psi(\tilde{x}) &= P_{0,0} \left( \mathrm{} \int_{\Gamma_k} \frac{\partial G}{\partial \nu_x}(\eta_j(\tilde{x}), y) \psi(\eta_k^{-1}(y)) \, d\sigma(y) \right) \\
		&= P_{0,0} \left( \mathrm{} \int_{\Gamma_0} \frac{\partial G}{\partial \nu_{\tilde{x}}}(\eta_j(\tilde{x}), \eta_k(\tilde{y})) \psi(\tilde{y}) r^{k(d-1)} \, d\sigma_0(\tilde{y}) \right).
	\end{align*}
	By the geometry of this radial scaling, these mappings naturally satisfy the algebraic identity $\eta_j(\tilde{x}) - \eta_k(\tilde{y}) = r^k (\eta_{j-k}(\tilde{x}) - \tilde{y})$ for all $\tilde{x}, \tilde{y} \in \Gamma_0$. Furthermore, the normal vectors are translation- and dilation-invariant, such that $\nu_{\eta_j(\tilde{x})} = \nu_{\eta_{j-k}(\tilde{x})} = \nu_{\tilde{x}}$. Thus,
	\begin{equation*}
		\frac{\partial G}{\partial \nu_{\tilde{x}}}(\eta_j(\tilde{x}), \eta_k(\tilde{y})) = r^{k(1-d)} \frac{\partial G}{\partial \nu_{\tilde{x}}}(\eta_{j-k}(\tilde{x}), \tilde{y}).
	\end{equation*}
	Substituting this back into the integral, the geometric scaling factor of the kernel $r^{k(1-d)}$ perfectly cancels the Jacobian of the surface measure $r^{k(d-1)}$. Thus, the spatial pullback is Toeplitz: $(\mathcal{P} \tilde{\mathsf{K}}^* \mathcal{P}^{-1})_{j,k} \psi = \tilde{K}^*_{j-k} \psi$. 
	
	Conjugating by the diagonal scaling matrices $ D^{\alpha}$ and $ D^{-\alpha}$ introduces the weight $r^{j\alpha} r^{-k\alpha} = r^{(j-k)\alpha}$, yielding $(\widetilde{\mathbb{K}}^*_\alpha)_{j,k} = r^{(j-k)\alpha} \tilde{K}^*_{j-k}$.
	\end{proof}

\begin{lemma}\label{lem: block decay of NP operator}
	For $m \neq 0$, the blocks in Lemma~\ref{lem: Toeplitz structure of NP operator} satisfy the exponential decay bounds:
	\begin{equation*}
		\| r^{m\alpha} \tilde{K}^*_m \|_{\mathcal{L}(W_0^{s,p}(\Gamma_0))} \le \begin{cases} 
			C r^{m\alpha}, & m > 0, \\ 
			C r^{|m|(d-\alpha)}, & m < 0, 
		\end{cases}
	\end{equation*}
	for some constant $C > 0$ independent of $m$. Moreover, $\tilde{K}^*_m \in \mathcal{K}(W_0^{s,p}(\Gamma_0))$ for all $m \in \mathbb{Z}$. 
\end{lemma}

\begin{proof}
	The diagonal block $\tilde{K}^*_0$ is the classical Neumann--Poincar\'e operator
	on a smooth boundary, which is compact by \eqref{eq:smoothing}. For $m\neq0$ the surfaces
	$\eta_m(\Gamma_0)$ and $\Gamma_0$ are disjoint, so
	$\tilde{K}^*_m = P_{0,0}A_{g_m}$, where $A_{g_m}\psi := \int_{\Gamma_0}g_m(\cdot,\tilde y)\psi(\tilde y)\,d\sigma_0(\tilde y)$
	and
	\[
	g_m(\tilde x,\tilde y) := \frac{\partial G}{\partial\nu_{\tilde x}}\bigl(\eta_m(\tilde x),\tilde y\bigr)
	\in C^\infty(\Gamma_0\times\Gamma_0).
	\]
	We note that for any $g\in C^1(\Gamma_0\times\Gamma_0)$,
	\begin{equation}\label{eq:kernel-bound}
		\|A_{g}\|_{\mathcal{L}(W^{s,p}_0(\Gamma_0))} \le C\,\|g\|_{C^1(\Gamma_0\times\Gamma_0)},
	\end{equation}
	with $C$ depending only on $\Gamma_0$, $s$, and $p$. Indeed, for $s\in[0,1)$
	compactness of $\Gamma_0$ gives $\|\psi\|_{L^1}\lesssim\|\psi\|_{L^p}$, whence
	$\|A_{g}\psi\|_{L^p}\lesssim\|g\|_{C^0}\|\psi\|_{L^p}$, while
	$|A_{g}\psi(\tilde x)-A_{g}\psi(\tilde x')|\le\|\nabla_{\tilde x}g\|_{C^0}
	|\tilde x-\tilde x'|\,\|\psi\|_{L^1}$ bounds the Gagliardo semi-norm, the
	resulting integral converging because $s<1$. For $s\in(-1,0)$, the kernel
	$(\tilde x,\tilde y)\mapsto g(\tilde y,\tilde x)$ is equally smooth and $A_{g}$
	is the adjoint on $W^{-s,p'}_0(\Gamma_0)$ of the operator it generates, so the
	bound follows from the case already proved. Compactness follows the same
	pattern: for $s\in[0,1)$ the operator $A_{g}$ maps bounded subsets of
	$W^{s,p}_0(\Gamma_0)$ into bounded subsets of $C^1(\Gamma_0)$, which embeds
	compactly into $W^{s,p}_0(\Gamma_0)$; for $s\in(-1,0)$, $A_{g}$ is the adjoint of
	the compact operator with kernel $g(\tilde y,\tilde x)$ and is therefore compact.
	Since $P_{0,0}$ is bounded, $\tilde{K}^*_m$ is compact for every $m$.
	
	By \eqref{eq:kernel-bound}, all $m$-dependence now sits in the kernel norms. Write
	$c_0 := \operatorname{dist}(0,\Gamma_0) > 0$ and $R_0 := \max_{\Gamma_0}|\tilde y|$,
	both finite and positive since $0\notin\overline{\Omega_0}$ and $\Gamma_0$ is
	compact.
	
	\noindent\emph{Case $m>0$.} Here $\eta_m(\tilde x)=r^m\tilde x$ lies within
	distance $r^mR_0$ of the origin, so $|\eta_m(\tilde x)-\tilde y|\ge c_0-r^mR_0$, which is
	bounded below uniformly once $r^mR_0\le c_0/2$; for the finitely many remaining
	$m$ the separation is positive by disjointness. Hence
	$\|g_m\|_{C^n}\le C_n$ for every $n\in\mathbb N$, with $C_n$ independent of $m$
	and no gain in $m$, and the decay of the block comes entirely from the sequence
	weight: $\|r^{m\alpha}\tilde{K}^*_m\| \le Cr^{m\alpha}$.
	
	\noindent\emph{Case $m<0$.} Let $l:=-m>0$. Since
	$\partial G/\partial\nu_x(x,y)$ is homogeneous of degree $1-d$ in $x-y$, writing
	$r^{-l}\tilde x-\tilde y=r^{-l}(\tilde x-r^l\tilde y)$ gives
	\[
	\frac{\partial G}{\partial\nu_{\tilde x}}\bigl(\eta_{-l}(\tilde x),\tilde y\bigr)
	=r^{l(d-1)}\frac{\partial G}{\partial\nu_{\tilde x}}\bigl(\tilde x,\eta_l(\tilde y)\bigr).
	\]
	Because $\partial G/\partial\nu_{\tilde x}(\tilde x,0)$ does not depend on
	$\tilde y$ and $\psi\in W^{s,p}_0(\Gamma_0)$ has zero mean, this monopole
	contribution may be subtracted, so $\tilde{K}^*_{-l}=P_{0,0}A_{\tilde g_l}$ with
	\[
	\tilde g_l(\tilde x,\tilde y) := r^{l(d-1)}\left[
	\frac{\partial G}{\partial\nu_{\tilde x}}\bigl(\tilde x,\eta_l(\tilde y)\bigr)
	-\frac{\partial G}{\partial\nu_{\tilde x}}(\tilde x,0)\right].
	\]
	The source points satisfy $|\eta_l(\tilde y)|\le r^lR_0$ while $|\tilde x|\ge c_0$,
	so for $r^lR_0\le c_0/2$ the segment joining $\eta_l(\tilde y)$ to the origin
	stays at distance at least $c_0/2$ from $\tilde x$, where the kernel and its
	$\tilde x$-derivatives have bounded gradients in the second variable. The mean
	value theorem then gives $\|\tilde g_l\|_{C^n}\le C_nr^{l(d-1)}r^l=C_nr^{ld}$
	for every $n\in\mathbb N$, and the same bound holds for the finitely many
	remaining $l$ by disjointness. Multiplying by the Toeplitz weight $r^{-l\alpha}$
	yields $Cr^{l(d-\alpha)}$.
\end{proof}

\begin{lemma}[Schatten class of the blocks]\label{lem:schatten-blocks}
	Let $q>d-1$. Then $\tilde K^*_m\in\mathfrak S_q\bigl(H^{-1/2}_0(\Gamma_0)\bigr)$ for every $m\in\mathbb Z$, with
	\[
	\|\tilde K^*_m\|_{\mathfrak S_q}\le
	\begin{cases}
		C, & m>0,\\
		Cr^{|m|d}, & m<0,
	\end{cases}
	\]
	and $C$ independent of $m$. For $m\neq0$ the same bounds hold for every $q\ge1$.
\end{lemma}

\begin{proof}
	The diagonal block $\tilde K^*_0$ is the mean-zero projection of the
	Neumann--Poincar\'e operator of the reference boundary $\Gamma_0$. On a smooth
	closed surface this operator is pseudodifferential of order $-1$
	\cite{MiyanishiRozenblum2019}, hence maps $H^{-1/2}(\Gamma_0)$ into
	$H^{1/2}(\Gamma_0)$. Since the embedding
	$H^{1/2}(\Gamma_0)\hookrightarrow H^{-1/2}(\Gamma_0)$ has singular values
	decaying at the rate $-1/(d-1)$, it belongs to $\mathfrak S_q$ exactly for
	$q>d-1$, and the ideal property transfers this to $\tilde K^*_0$; the
	projection $P_{0,0}$ is bounded and changes nothing. This is consistent with the
	sharp asymptotics of the Neumann--Poincar\'e eigenvalues on smooth boundaries
	\cite{MiyanishiSuzuki2017, MiyanishiRozenblum2019, Miyanishi2022}, which show that the exponent $d-1$
	cannot be improved.
	
	For $m\neq0$ the surfaces $\eta_m(\Gamma_0)$ and $\Gamma_0$ are disjoint, so the
	kernel $g_m$ of $\tilde K^*_m$ is smooth on $\Gamma_0\times\Gamma_0$ and the
	operator is smoothing of infinite order: for every $t$ it maps
	$H^{-1/2}_0(\Gamma_0)$ boundedly into $H^{t}(\Gamma_0)$, with norm controlled by
	finitely many derivatives of $g_m$. Fix $q\ge1$ and choose $t$ with
	$t+1/2>(d-1)/q$. Write $\tilde K^*_m=\iota_t A_m$, where $A_m$ denotes
	$\tilde K^*_m$ regarded as a bounded map
	$H^{-1/2}_0(\Gamma_0)\to H^{t}(\Gamma_0)$ and $\iota_t$ the inclusion
	$H^{t}(\Gamma_0)\hookrightarrow H^{-1/2}(\Gamma_0)$. The singular values of
	$\iota_t$ decay at the rate $-(t+1/2)/(d-1)$, so $\iota_t\in\mathfrak S_q$ and
	the ideal property~\eqref{eq:app-ideal} gives
	$\|\tilde K^*_m\|_{\mathfrak S_q}\le\|\iota_t\|_{\mathfrak S_q}\|A_m\|$. The
	proof of Lemma~\ref{lem: block decay of NP operator} bounds finitely many
	derivatives of $g_m$ by $C_n$ for $m>0$ and by $C_nr^{|m|d}$ for $m<0$, uniformly
	in $m$, so $\|A_m\|$ obeys the same bounds and hence so does
	$\|\tilde K^*_m\|_{\mathfrak S_q}$, for every $q\ge1$.
\end{proof}

\subsection{Block-Toeplitz structure of $\widetilde{\mathbb{S}}_\alpha$}
\begin{lemma}[Structure of the scaled Single Layer operator]\label{lem: Toeplitz structure Single layer}
	For $d \in \{2, 3\}$, the scaled single-layer potential $\widetilde{\mathbb{S}}_\alpha$ admits the exact factorisation:
	\begin{equation*}
		\widetilde{\mathbb{S}}_\alpha = D^\alpha \widetilde{\mathbb{S}}_{\mathrm{Toep}} D^{1-\alpha},
	\end{equation*}
	where $\widetilde{\mathbb{S}}_{\mathrm{Toep}}$ is the block-Toeplitz operator:
	\begin{equation*}
		\widetilde{\mathbb{S}}_{\mathrm{Toep}} := \begin{pmatrix} \tilde{S}_0 & \tilde{S}_{-1} & \tilde{S}_{-2} & \cdots \\ \tilde{S}_1 & \tilde{S}_0 & \tilde{S}_{-1} & \cdots \\ \tilde{S}_2 & \tilde{S}_1 & \tilde{S}_0 & \cdots \\ \vdots & \vdots & \vdots & \ddots \end{pmatrix}.
	\end{equation*}
	The block entries are given by $(\widetilde{\mathbb{S}}_{\mathrm{Toep}})_{j,k} = \tilde{S}_{j-k}$ for $j, k \ge 0$, where $\tilde{S}_m : W_0^{s,p}(\Gamma_0) \to W_0^{s,p}(\Gamma_0)$ is defined by:
	\begin{equation*}
		(\tilde{S}_m \psi)(\tilde{x}) = P_{0,0} \left( \int_{\Gamma_0} G(\eta_m(\tilde{x}), \tilde{y}) \psi(\tilde{y}) \, d\sigma_0(\tilde{y}) \right).
	\end{equation*}
\end{lemma}

\begin{proof}
		Recall that
	$$
	\widetilde{\mathbb{S}}_\alpha
	:=
	D^\alpha \mathcal{P}\mathcal T\widetilde{S}_\Gamma
	\mathcal \mathcal \mathcal T^{-1}\mathcal{P}^{-1} D^{-\alpha}.
	$$
	It follows directly from the definition of $\mathcal T$ that
	\[
\tilde{\mathsf{S}}
:=
\mathcal{T}\widetilde{S}_\Gamma\mathcal{T}^{-1}
=
\begin{pmatrix}
	\tilde{\mathsf{S}}_{0,0} & \tilde{\mathsf{S}}_{0,1} & \cdots \\
	\tilde{\mathsf{S}}_{1,0} & \tilde{\mathsf{S}}_{1,1} & \cdots \\
	\vdots & \vdots & \ddots
\end{pmatrix},
\]
	where 
	\[
(\tilde{\mathsf{S}}_{j,k}\varphi_k)(x)
:=
P_{0,j}\left(
\int_{\Gamma_k}
G(\cdot,y)\varphi_k(y)\,d\sigma(y)
\right)(x),
\qquad x\in\Gamma_j.
\]
We next evaluate the purely spatial pullback of the $(j,k)$-th block. For $\psi \in W_0^{s,p}(\Gamma_0)$,
	\begin{align*}
		(\mathcal{P} \tilde{\mathsf{S}} \mathcal{P}^{-1})_{j,k} \psi(\tilde{x}) &= P_{0,0} \left( \int_{\Gamma_k} G(\eta_j(\tilde{x}), y) \psi(\eta_k^{-1}(y)) \, d\sigma(y) \right) \\
		&= P_{0,0} \left( \int_{\Gamma_0} G\big(r^k(\eta_{j-k}(\tilde{x}) - \tilde{y})\big) \psi(\tilde{y}) r^{k(d-1)} \, d\sigma_0(\tilde{y}) \right).
	\end{align*}
	By the homogeneity of the fundamental solution, $G(r^k z) = r^{k(2-d)} G(z) + c_k$, with $c_k = (2\pi)^{-1}\log r^k$ for $d = 2$ and $c_k = 0$ for $d = 3$. Because $\psi \in W_0^{s,p}(\Gamma_0)$ has exactly zero mean ($\int_{\Gamma_0} \psi \, d\sigma_0 = 0$), the scalar shift $c_k$ integrates to zero identically, allowing the fractional homogeneity relation to hold within the mean-zero subspace. Factoring out the scaling yields:
	\begin{align*}
		(\mathcal{P} \tilde{\mathsf{S}} \mathcal{P}^{-1})_{j,k} \psi(\tilde{x}) &= r^{k(2-d)} r^{k(d-1)} P_{0,0} \left( \int_{\Gamma_0} G(\eta_{j-k}(\tilde{x}), \tilde{y}) \psi(\tilde{y}) \, d\sigma_0(\tilde{y}) \right) \\
		&= r^k (\tilde{S}_{j-k} \psi)(\tilde{x}).
	\end{align*}
	Conjugating by the diagonal scaling weights $D^\alpha$ and $D^{-\alpha}$ gives the algebraic factorisation:
	\begin{equation*}
		(\widetilde{\mathbb{S}}_\alpha)_{j,k} = r^{j\alpha} \big(r^k \tilde{S}_{j-k}\big) r^{-k\alpha} = r^{j\alpha} \tilde{S}_{j-k} r^{k(1-\alpha)}.
	\end{equation*}
	This structurally verifies the factorisation $\widetilde{\mathbb{S}}_\alpha = D^\alpha \widetilde{\mathbb{S}}_{\mathrm{Toep}} D^{1-\alpha}$.

\end{proof}

\begin{lemma}\label{lem:decay-SL-blocks}
	For $m \neq 0$, the Toeplitz blocks in Lemma~\ref{lem: Toeplitz structure Single layer} satisfy the exponential decay bounds:
	\begin{equation*}
		\| \tilde{S}_m \|_{\mathcal{L}(W_0^{s,p}(\Gamma_0))} \le \begin{cases} 
			C r^{m}, & m > 0, \\ 
			C r^{|m|(d-1)}, & m < 0, 
		\end{cases}
	\end{equation*}
	for some constant $C > 0$ independent of $m$.
\end{lemma}
\begin{proof}
	As in the proof of Lemma~\ref{lem: block decay of NP operator}, for $m\neq0$ the
	source and target boundaries are strictly separated, so the kernels are
	$C^\infty$ on $\Gamma_0\times\Gamma_0$ and \eqref{eq:kernel-bound} reduces the
	claim to bounding their $C^1$ norms, with a constant depending only on
	$\Gamma_0$, $s$, and $p$. We keep the notation
	$c_0=\operatorname{dist}(0,\Gamma_0)>0$ and $R_0=\max_{\Gamma_0}|\tilde y|$.
	Throughout, the homogeneity $G(\rho z)=\rho^{2-d}G(z)$ holds exactly for $d=3$
	and modulo an additive constant for $d=2$; since every density occurring below
	lies in a mean-zero space, such constants integrate to zero and are omitted.
	
	\noindent\emph{Case $m<0$.} Set $l=-m>0$. Homogeneity of $G$ of degree $2-d$
	gives $G(\eta_{-l}(\tilde x),\tilde y)=r^{l(d-2)}G(\tilde x,\eta_l(\tilde y))$.
	Subtracting the monopole contribution at the origin,
	\begin{equation*}
		(\tilde S_{-l}\psi)(\tilde x)=P_{0,0}\left(r^{l(d-2)}\int_{\Gamma_0}
		\bigl[G(\tilde x,\eta_l(\tilde y))-G(\tilde x,0)\bigr]\psi(\tilde y)\,d\sigma_0(\tilde y)\right).
	\end{equation*}
	Since $|\eta_l(\tilde y)|\le r^l R_0$ while $|\tilde x|\ge c_0$, for
	$r^l R_0\le c_0/2$ the segment joining $\eta_l(\tilde y)$ to the origin stays at
	distance at least $c_0/2$ from $\tilde x$, where $G(\tilde x,\cdot)$ and its
	$\tilde x$-derivatives have bounded gradients in the second variable. The mean
	value theorem gives $\mathcal O(r^l)$ for the bracket and its derivatives,
	whence $\|\tilde S_{-l}\|_{\mathcal L}\le Cr^{l(d-2)}r^l=Cr^{l(d-1)}$; the
	finitely many remaining $l$ are covered by disjointness of the components.
	
	\noindent\emph{Case $m>0$.} We use duality. For $\phi,\psi$ in the relevant
	mean-zero spaces, the symmetry of $G$ give
	\begin{equation*}
		\langle\tilde S_m\psi,\phi\rangle_{L^2(\Gamma_0)}
		=\int_{\Gamma_0}\psi(\tilde y)\left(\int_{\Gamma_0}
		G(\eta_m(\tilde x),\tilde y)\phi(\tilde x)\,d\sigma_0(\tilde x)\right)d\sigma_0(\tilde y),
	\end{equation*}
	and the scaling identity
	$G(\eta_m(\tilde x),\tilde y)=r^{-m(d-2)}G(\eta_{-m}(\tilde y),\tilde x)$
	evaluates the inner integral as $r^{-m(d-2)}(\tilde S_{-m}\phi)(\tilde y)$.
	Hence $\tilde S_m^{*}=r^{-m(d-2)}\tilde S_{-m}$, the adjoint being taken with
	respect to the $L^2(\Gamma_0)$ pairing, for which $W^{-s,p'}_0(\Gamma_0)$ is the
	dual of $W^{s,p}_0(\Gamma_0)$ as in the proof of
	Lemma~\ref{lem:sequence-isomorphism}, so that
	$\|\tilde S_m\|_{\mathcal L(W^{s,p}_0)}
	=\|\tilde S_m^{*}\|_{\mathcal L(W^{-s,p'}_0)}$. The case $m<0$ was proved for
	every admissible pair of parameters, in particular for $(-s,p')$, so
	\begin{equation*}
		\|\tilde S_m\|_{\mathcal L(W^{s,p}_0)}
		=r^{-m(d-2)}\|\tilde S_{-m}\|_{\mathcal L(W^{-s,p'}_0)}
		\le r^{-m(d-2)}\cdot Cr^{m(d-1)}=Cr^{m}.
	\end{equation*}
\end{proof}

\subsection{Boundedness and symmetrisation}
\begin{theorem}[Boundedness on the infinite chain]\label{thm:boundedness}
	Let $p\in(1,\infty)$ and $s\in(-1,1)$, and define the geometric scaling exponent
	\[
	\alpha=\frac{d-1}{p}-s.
	\]
	If $0<\alpha<d$, then the rescaled sequence-space operators $\widetilde{\mathbb{S}}_\alpha$, $\widetilde{\mathbb{K}}^*_\alpha$, and $\widetilde{\mathbb{K}}_\alpha$ are bounded linear operators on $\ell^p\bigl(\mathbb{N},W_0^{s,p}(\Gamma_0)\bigr)$.
	
	Consequently, through the isomorphisms established above, the layer potentials on the infinite chain are bounded operators on the global mean-zero space:
	\[
	\widetilde{S}_\Gamma,\,\widetilde{K}^*_\Gamma,\,\widetilde{K}_\Gamma
	\in
	\mathcal{L}\bigl(W_0^{s,p}(\Gamma)\bigr),
	\]
	where $s$ and $p$ are determined by $0<\alpha<d$.
\end{theorem}

\begin{proof}
	By virtue of the isomorphism between $W_0^{s,p}(\Gamma)$ and $V_0^{s,p}(\Gamma)$, it suffices to establish boundedness of the rescaled operators on the sequence space.
	
	To do so, we rely on the following Schur test for operator-valued matrices. Let $M = [M_{j,k}]$ have entries $M_{j,k} \in \mathcal{L}(X)$ with
	\[
	\sup_j \sum_k \|M_{j,k}\|_{\mathcal{L}(X)} \le C,
	\qquad
	\sup_k \sum_j \|M_{j,k}\|_{\mathcal{L}(X)} \le C .
	\]
	Then $M$ is bounded on $\ell^p(\mathbb{N},X)$ for every $p \in (1,\infty)$, with $\|M\| \le C$. Indeed, writing $v_k := \|u_k\|_X$ and applying H\"older with the weights $\|M_{j,k}\|^{1/p'}$ and $\|M_{j,k}\|^{1/p} v_k$,
	\[
	\|(Mu)_j\|_X^p \le \Bigl(\sum_k \|M_{j,k}\| \Bigr)^{p/p'} \sum_k \|M_{j,k}\| \, v_k^p
	\le C^{p/p'} \sum_k \|M_{j,k}\| \, v_k^p ,
	\]
	and summing over $j$ and using the column bound gives $\|Mu\|_p^p \le C^p \|u\|_p^p$.
	
	The operator $\widetilde{\mathbb{K}}^*_\alpha$ is exactly block-Toeplitz, with its $(j,k)$-th block given by $(\widetilde{\mathbb{K}}^*_\alpha)_{j,k} = r^{(j-k)\alpha} \widetilde{K}^*_{j-k}$. By the exponential decay bounds of Lemma~\ref{lem: block decay of NP operator}, 
	\[
	\bigl\|(\widetilde{\mathbb{K}}^*_\alpha)_{j,k}\bigr\|_{\mathcal{L}(W_0^{s,p}(\Gamma_0))} \lesssim
	\begin{cases}
		r^{(j-k)\alpha}, & j \ge k, \\
		r^{(k-j)(d-\alpha)}, & j < k.
	\end{cases}
	\]
	Since $r\in(0,1)$ and $0<\alpha<d$, these operator norms are bounded by a scalar sequence $a_{j-k}$ with $\sum_{m\in\mathbb{Z}} a_m < \infty$. The row and column sums are therefore both bounded by $\sum_{m} a_m$, and the Schur test yields boundedness on $\ell^p$. Boundedness of $\widetilde{\mathbb{K}}_\alpha$ follows identically by transposition.
	
	For the single-layer operator, the factorisation $\widetilde{\mathbb{S}}_\alpha = D^\alpha \widetilde{\mathbb{S}}_{\mathrm{Toep}} D^{1-\alpha}$ yields the block entries $(\widetilde{\mathbb{S}}_\alpha)_{j,k} = r^{j\alpha}\widetilde{S}_{j-k}r^{k(1-\alpha)}$ for $j, k \ge 0$. Substituting the decay bounds for $\widetilde{S}_m$ from Lemma~\ref{lem:decay-SL-blocks} yields:
	\[
	\bigl\|(\widetilde{\mathbb{S}}_\alpha)_{j,k}\bigr\|_{\mathcal{L}(W_0^{s,p}(\Gamma_0))}
	\lesssim
	\begin{cases}
		r^{j\alpha} r^{j-k} r^{k(1-\alpha)} = r^k r^{(j-k)(\alpha+1)} \le r^{(j-k)(\alpha+1)}, & j \ge k, \\[1mm]
		r^{j\alpha} r^{(k-j)(d-1)} r^{k(1-\alpha)} = r^j r^{(k-j)(d-\alpha)} \le r^{(k-j)(d-\alpha)}, & j < k.
	\end{cases}
	\]
	Because $\alpha+1 > 0$ and $d-\alpha > 0$, the operator matrix $\widetilde{\mathbb{S}}_\alpha$ is again entry-wise dominated in norm by a scalar Toeplitz matrix whose entries belong to $\ell^1(\mathbb{Z})$, so the row and column sums are uniformly bounded and the Schur test applies.
\end{proof}

\begin{remark}
	Within the range $p \in (1,\infty)$, $s \in (-1,1)$, the exponent
	$\alpha = (d-1)/p - s$ always satisfies $\alpha < d$, so only the lower bound $\alpha>0$ is
	a genuine restriction. Note that, by the bounds of
	Lemma~\ref{lem:decay-SL-blocks}, the single layer requires only $\alpha > -1$, which
	likewise holds automatically. For $\widetilde{\mathbb{K}}^*_\alpha$ and
	$\widetilde{\mathbb{K}}_\alpha$ the constraint $\alpha > 0$ cannot be relaxed,
	as the chain of annuli shows in Section~\ref{sec:annuli}.
\end{remark}

\begin{lemma}
	\label{lem: slefadjointness of S_alpha}
	The scaled single layer operator $\widetilde{\mathbb{S}}_{\frac{d-1}{2}}$ is self-adjoint, negative definite, and injective in $\ell^2\big(\mathbb{N},L_0^2(\Gamma_0)\big)$. Furthermore, it extends to a coercive isomorphism from the weighted sequence space $D^{-1/2}\ell^2\big(\mathbb{N},H_0^{-1/2}(\Gamma_0)\big)$ onto its dual space $D^{1/2}\ell^2\big(\mathbb{N},H_0^{1/2}(\Gamma_0)\big)$, the duality being that of the $\ell^2$ pairing, for which the two weights are conjugate.
\end{lemma}
\begin{proof}
	Throughout we take $\alpha=\frac{d-1}{2}$, which lies in $(0,d)$, so that $\widetilde{\mathbb{S}}_{\frac{d-1}{2}}$ is bounded by Theorem~\ref{thm:boundedness}. For $\mathbf{u},\mathbf{v}\in \ell^2\big(\mathbb{N},L_0^2(\Gamma_0)\big)$ define the corresponding zero-mean physical densities $f_0=\mathcal{P}^{-1}D^{-\alpha}\mathbf{u}$ and $g_0=\mathcal{P}^{-1}D^{-\alpha}\mathbf{v}$. A change of variables via $x=\eta_j(\tilde{x})$ shows that the local $L^2$ inner products scale geometrically as $d\sigma(x)=r^{j(d-1)}d\sigma_0(\tilde{x})$, so that the global physical pairing satisfies
	$$
	\langle \widetilde{S}_\Gamma f_0,g_0\rangle_{L^2(\Gamma)} = \left\langle D^{d-1}\mathcal{P}\widetilde{S}_\Gamma \mathcal{P}^{-1}D^{-\alpha}\mathbf{u}, D^{-\alpha}\mathbf{v} \right\rangle_{\ell^2} = \left\langle D^{d-1-\alpha}\mathcal{P}\widetilde{S}_\Gamma \mathcal{P}^{-1}D^{-\alpha}\mathbf{u}, \mathbf{v} \right\rangle_{\ell^2}.
	$$
	Our choice of $\alpha$ balances the weights, $d-1-\alpha=\alpha$, so the operator in the right-hand inner product is exactly $\widetilde{\mathbb{S}}_{\frac{d-1}{2}}$. Since the physical operator $\widetilde{S}_\Gamma$ is symmetric, this congruence makes $\widetilde{\mathbb{S}}_{\frac{d-1}{2}}$ symmetric, hence self-adjoint, on the sequence space. Taking $\mathbf{v}=\mathbf{u}$ and using Lemma~\ref{lem:uniform-coercivity},
	$$
	\left\langle -\widetilde{\mathbb{S}}_{\frac{d-1}{2}}\mathbf{u}, \mathbf{u} \right\rangle_{\ell^2}
	=\langle -\widetilde{S}_\Gamma f_0,f_0\rangle_{L^2(\Gamma)}
	\ \ge\ c\,\|f_0\|_{H^{-1/2}(\Gamma)}^2\ >\ 0
	\qquad (\mathbf{u}\neq0),
	$$
	so the operator is negative definite and in particular injective.
	
	To establish the isomorphism, we use the same coercivity bound together with the geometric scaling of the norms. The global Sobolev norm at $s=-1/2$ corresponds by Lemma~\ref{lem:scaling-isomorphism} to the weighted sequence norm with exponent $\frac{d-1}{2}+\frac12=\frac{d}{2}$, that is $\|f_0\|_{H^{-1/2}(\Gamma)} \asymp \|D^{d/2}\mathcal{P}f_0\|_{\ell^2(\mathbb{N},H_0^{-1/2}(\Gamma_0))}$. Substituting $f_0 = \mathcal{P}^{-1} D^{-\frac{d-1}{2}}\mathbf{u}$ shifts the weight by $D^{-(d-1)/2}$ and yields
	$$
	\|f_0\|_{H^{-1/2}(\Gamma)} \asymp \left\|D^{1/2}\mathbf{u}\right\|_{\ell^2(\mathbb{N},H_0^{-1/2}(\Gamma_0))}.
	$$
	Inserting this into the coercivity bound above gives
	$$
	\left\langle -\widetilde{\mathbb{S}}_{\frac{d-1}{2}}\mathbf{u}, \mathbf{u} \right\rangle_{\ell^2}
	\ \ge\ \widetilde{C} \left\|D^{1/2}\mathbf{u}\right\|_{\ell^2(\mathbb{N},H_0^{-1/2}(\Gamma_0))}^2,
	$$
	and by the Lax--Milgram lemma $\widetilde{\mathbb{S}}_{\frac{d-1}{2}}$ extends uniquely to an isomorphism between the stated weighted sequence spaces.
\end{proof}

\begin{lemma}[Energy symmetrisation]\label{lem: Symmetrisation of K toeplitz}
	Let $\mathbb J:=-D^{-1/2}\widetilde{\mathbb S}_{\frac{d-1}{2}}D^{-1/2}$. Then
	$\mathbb J$ is a coercive self-adjoint isomorphism from
	$\ell^2\bigl(\mathbb N,H^{-1/2}_0(\Gamma_0)\bigr)$ onto
	$\ell^2\bigl(\mathbb N,H^{1/2}_0(\Gamma_0)\bigr)$, and
	\begin{equation}\label{eq:shifted_symmetrization}
		\mathbb J\,\widetilde{\mathbb K}^*_{d/2}
		=\bigl(\widetilde{\mathbb K}^*_{d/2}\bigr)^{*}\,\mathbb J .
	\end{equation}
	Consequently $\langle u,v\rangle_{\mathbb J}:=\langle\mathbb Ju,v\rangle$ is an
	inner product on $\ell^2\bigl(\mathbb N,H^{-1/2}_0(\Gamma_0)\bigr)$, equivalent
	to the standard one, for which $\widetilde{\mathbb K}^*_{d/2}$ is self-adjoint.
\end{lemma}
\begin{proof}
	Conjugating the identity $\widetilde S_\Gamma\widetilde K^*_\Gamma
	=\widetilde K_\Gamma\widetilde S_\Gamma$ of
	Lemma~\ref{lem:projected-symmetrization} by the isomorphism $\mathcal M_\alpha$
	gives $\widetilde{\mathbb S}_\alpha\widetilde{\mathbb K}^*_\alpha
	=\widetilde{\mathbb K}_\alpha\widetilde{\mathbb S}_\alpha$ for every admissible
	$\alpha$. By Lemma~\ref{lem: slefadjointness of S_alpha},
	$\widetilde{\mathbb S}_{\frac{d-1}{2}}$ is a coercive self-adjoint isomorphism
	from $D^{-1/2}\ell^2(\mathbb N,H^{-1/2}_0(\Gamma_0))$ onto
	$D^{1/2}\ell^2(\mathbb N,H^{1/2}_0(\Gamma_0))$; pre- and post-multiplying by
	$D^{-1/2}$ gives the first assertion.
	
	Take $\alpha=\frac{d-1}{2}$. As in the congruence of
Lemma~\ref{lem: slefadjointness of S_alpha}, for
$\mathbf u,\mathbf v\in\ell^2(\mathbb N,L^2_0(\Gamma_0))$ with physical
densities $f_0=\mathcal M_\alpha^{-1}\mathbf u$ and
$g_0=\mathcal M_\alpha^{-1}\mathbf v$, the measure weight $r^{j(d-1)}$ splits as
$r^{j\alpha}\cdot r^{j\alpha}$, one factor for each argument, so
\[
\bigl\langle\widetilde{\mathbb K}^*_{\alpha}\mathbf u,\mathbf v\bigr\rangle_{\ell^2}
=\bigl\langle\widetilde K^*_\Gamma f_0,g_0\bigr\rangle_{L^2(\Gamma)}
=\bigl\langle f_0,\widetilde K_\Gamma g_0\bigr\rangle_{L^2(\Gamma)}
=\bigl\langle\mathbf u,\widetilde{\mathbb K}_{\alpha}\mathbf v\bigr\rangle_{\ell^2},
\]
the middle equality because $\widetilde K_\Gamma$ is by construction the
$L^2(\Gamma)$ adjoint of $\widetilde K^*_\Gamma$. Hence
$\widetilde{\mathbb K}_{\frac{d-1}{2}}
=\bigl(\widetilde{\mathbb K}^*_{\frac{d-1}{2}}\bigr)^{*}$, and the identity of
the first paragraph reads
\[
\widetilde{\mathbb S}_{\frac{d-1}{2}}\widetilde{\mathbb K}^*_{\frac{d-1}{2}}
=\bigl(\widetilde{\mathbb K}^*_{\frac{d-1}{2}}\bigr)^{*}\widetilde{\mathbb S}_{\frac{d-1}{2}} .
\]
Multiplying on both sides by $-D^{-1/2}$ and inserting
$D^{-1/2}D^{1/2}=I$ between the factors gives
\eqref{eq:shifted_symmetrization}, since
$\widetilde{\mathbb K}^*_{d/2}
=D^{1/2}\widetilde{\mathbb K}^*_{\frac{d-1}{2}}D^{-1/2}$ and $D$ is
self-adjoint.

	Finally, for $\mathbf u,\mathbf v\in\ell^2(\mathbb N,H^{-1/2}_0(\Gamma_0))$,
\eqref{eq:shifted_symmetrization} gives
\[
\bigl\langle\widetilde{\mathbb K}^*_{d/2}\mathbf u,\mathbf v\bigr\rangle_{\mathbb J}
=\bigl\langle\mathbb J\widetilde{\mathbb K}^*_{d/2}\mathbf u,\mathbf v\bigr\rangle
=\bigl\langle\bigl(\widetilde{\mathbb K}^*_{d/2}\bigr)^{*}\mathbb J\mathbf u,\mathbf v\bigr\rangle
=\bigl\langle\mathbb J\mathbf u,\widetilde{\mathbb K}^*_{d/2}\mathbf v\bigr\rangle
=\bigl\langle\mathbf u,\widetilde{\mathbb K}^*_{d/2}\mathbf v\bigr\rangle_{\mathbb J},
\]
so $\widetilde{\mathbb K}^*_{d/2}$ is symmetric, hence self-adjoint, for
$\langle\cdot,\cdot\rangle_{\mathbb J}$; that this is an inner product
equivalent to the standard one is the coercivity of $\mathbb J$.
\end{proof}

\section{The spectrum for the infinite chain (main results 1 and 2)}\label{sec:spectrum}

In this section, we provide a characterisation of the spectrum of the projected Neumann--Poincar\'e operator ${K}^*_\Gamma$ acting on the mean-zero Sobolev spaces $W_0^{s,p}(\Gamma)$, i.e. $\widetilde{K}^*_\Gamma$. To do so, we apply block-Toeplitz operator theory to the rescaled sequence-space operator $\widetilde{\mathbb{K}}^*_\alpha$ acting on $\ell^p\bigl(\mathbb{N}, W_0^{s,p}(\Gamma_0)\bigr)$, with 
$$\alpha = \frac{d-1}{p} - s \in (0, d).$$

\subsection{The symbol of $\widetilde{\mathbb{K}}^*_\alpha$}
\begin{definition}[Floquet--Toeplitz Symbol]\label{def:symbol}
	For $\theta \in [-\pi, \pi]$, we define the operator-valued symbol
	$$
	\kappa_\alpha(\theta) := \sum_{m\in\mathbb{Z}} e^{im\theta} r^{m\alpha} \tilde{K}^*_m.
	$$
\end{definition}
The exponential decay of the Toeplitz blocks $r^{m\alpha}\tilde{K}^*_m$, established in Lemma~\ref{lem: block decay of NP operator}, implies that the Laurent series $\sum_{m\in\mathbb{Z}} z^m r^{m\alpha}\tilde{K}^*_m$ converges in the operator norm on the annulus $r^{d-\alpha}<|z|<r^{-\alpha}$, which contains the unit circle precisely because $0<\alpha<d$; the series above is its restriction to $|z|=1$. In particular $\kappa_\alpha$ belongs to the Wiener algebra $\mathcal W(\mathcal{L}(W_0^{s,p}(\Gamma_0)))$, defined in Appendix~\ref{app:toeplitz}, and the map $\theta \mapsto \kappa_\alpha(\theta)$ is continuous in the operator norm. Furthermore, since each block $\tilde{K}^*_m$ is compact (Lemma~\ref{lem: block decay of NP operator}), $\kappa_\alpha(\theta)$ is a compact operator for every $\theta \in [-\pi, \pi]$.

\subsection{Main result 1: the spectrum in $W_0^{s,p}(\Gamma)$}
\begin{theorem}[Essential spectrum]\label{thm:essential-spectrum}
	Let $\widetilde{K}^*_\Gamma$ be the projected Neumann--Poincar\'e operator acting on $W_0^{s,p}(\Gamma)$, and let $\alpha = (d - 1)/p - s \in (0, d)$. Then
	$$
	\sigma_{\mathrm{ess}}(\widetilde{K}^*_\Gamma) = \bigcup_{\theta\in[-\pi,\pi]} \sigma\bigl(\kappa_\alpha(\theta)\bigr).
	$$
\end{theorem} 

\begin{proof}
	Because $\widetilde{K}^*_\Gamma$ is similar to $\widetilde{\mathbb{K}}^*_\alpha$, we have
	$$
	\sigma_{\mathrm{ess}}(\widetilde{K}^*_\Gamma) = \sigma_{\mathrm{ess}}(\widetilde{\mathbb{K}}^*_\alpha).
	$$
	Since $\kappa_\alpha$ belongs to the operator-valued Wiener algebra and is compact-valued, Theorem~\ref{thm:app-fredholm} implies that $\lambda I - \widetilde{\mathbb{K}}^*_\alpha$ is Fredholm if and only if $\lambda I - \kappa_\alpha(\theta)$ is boundedly invertible for every $\theta \in [-\pi, \pi]$. That is, $\lambda I - \widetilde{\mathbb{K}}^*_\alpha$ is not Fredholm if and only if $\lambda \in \sigma\bigl(\kappa_\alpha(\theta)\bigr)$ for at least one $\theta \in [-\pi, \pi]$. 
\end{proof}

\begin{remark}
	Note that since each $\kappa_\alpha(\theta)$ is compact on the infinite-dimensional space $W^{s,p}_0(\Gamma_0)$, we necessarily have $0 \in \sigma(\kappa_\alpha(\theta))$ for every $\theta \in [-\pi,\pi]$. Therefore, $0 \in \sigma_{\mathrm{ess}}(\widetilde{K}^*_\Gamma)$.
\end{remark}

\begin{remark}\label{rem:alpha-only}
	We note that the essential spectrum depends on the functional setting only
	through $\alpha$. This is because $\kappa_\alpha(\theta)$ is
	smoothing of order one, since each block maps $W_0^{s,p}(\Gamma_0)$ into
	$W_0^{s+1,p}(\Gamma_0)$ with the decay bounds of
	Lemma~\ref{lem: block decay of NP operator}. Hence, any $v$ satisfying
	$v = \lambda^{-1}\kappa_\alpha(\theta) v$ with $\lambda \neq 0$ is smooth by
	bootstrapping, and therefore its non-zero spectrum is common to every
	$W_0^{s,p}(\Gamma_0)$.
\end{remark}

\paragraph{A computational approach to the Fredholm index.}  The Fredholm index is constant on each connected component of $\mathbb{C} \setminus \sigma_{\mathrm{ess}}(\widetilde{K}^*_\Gamma)$. We now show that these indices can be computed from finite-dimensional approximations of the symbol. This indicates that for a sufficiently large discretization of the symbol, the Fredholm index can be evaluated numerically by tracing the determinant along a discrete grid of $\theta \in [-\pi, \pi]$ and computing its winding number around the origin.

Let $(\mathcal Q_M)_{M\ge1}$ be a uniformly bounded family of finite-rank projections on $W_0^{s,p}(\Gamma_0)$ such that $\mathcal Q_M \to I$ and $\mathcal Q_M^* \to I$ strongly on $W_0^{s,p}(\Gamma_0)$ and its dual, respectively, and assume that $\dim \operatorname{Im}(\mathcal Q_M) = M$. We define the truncated finite-dimensional symbol
$$
\kappa^{(M)}_\alpha(\theta) := \mathcal Q_M \kappa_\alpha(\theta) \mathcal Q_M \big|_{\operatorname{Im}(\mathcal Q_M)}, \quad \theta \in [-\pi, \pi].
$$
Thus $\kappa^{(M)}_\alpha(\theta)$ may be identified with an $M \times M$ matrix-valued symbol. Since $\kappa_\alpha(\theta)$ is a norm-continuous family of compact operators,
$$
\sup_{\theta\in[-\pi,\pi]} \|\kappa_\alpha(\theta) - \mathcal Q_M\kappa_\alpha(\theta)\mathcal Q_M\| \longrightarrow 0 \quad \text{as } M \to \infty.
$$ 

\begin{theorem}[Fredholm Index]\label{thm:fredholm-index}
	Let $\lambda \notin \sigma_{\mathrm{ess}}(\widetilde{K}^*_\Gamma)$. For every $M$ sufficiently large that
	$$
	\sup_{\theta\in[-\pi,\pi]} \|\kappa_\alpha(\theta) - \mathcal Q_M\kappa_\alpha(\theta)\mathcal Q_M\| < \min_{\theta\in[-\pi,\pi]} \|(\lambda I - \kappa_\alpha(\theta))^{-1}\|^{-1},
	$$
	the Fredholm index is given by
	$$
	\operatorname{ind}(\lambda I - \widetilde{K}^*_\Gamma) = - \frac{1}{2\pi i} \int_{-\pi}^\pi \frac{d}{d\theta} \log \det \left( \lambda I_M - \kappa^{(M)}_\alpha(\theta) \right) d\theta.
	$$
\end{theorem} 

\begin{proof}
	By the similarity between $\widetilde{K}^*_\Gamma$ and $\widetilde{\mathbb{K}}^*_\alpha$, we have $\operatorname{ind}(\lambda I - \widetilde{K}^*_\Gamma) = \operatorname{ind}(\lambda I - \widetilde{\mathbb{K}}^*_\alpha)$. The remainder of the proof relies on standard homotopy arguments. Let
	$$
	\delta_\lambda := \min_{\theta\in[-\pi,\pi]} \|(\lambda I - \kappa_\alpha(\theta))^{-1}\|^{-1} > 0,
	$$
	which exists because $\lambda \notin \sigma_{\mathrm{ess}}(\widetilde{K}^*_\Gamma)$. Fix $M$ large enough that the approximation error is strictly smaller than $\delta_\lambda$, and consider the homotopy
	$$
	\kappa^{(M)}_{\alpha,t}(\theta) = (1 - t)\kappa_\alpha(\theta) + t \mathcal Q_M\kappa_\alpha(\theta)\mathcal Q_M, \quad 0 \le t \le 1.
	$$
	Since $\kappa^{(M)}_{\alpha,t}(\theta)-\kappa_\alpha(\theta)$ has norm at most $t\sup_\theta\|\kappa_\alpha(\theta)-\mathcal Q_M\kappa_\alpha(\theta)\mathcal Q_M\|<\delta_\lambda$ for every $t\in[0,1]$, the operator $\lambda I - \kappa^{(M)}_{\alpha,t}(\theta)$ remains boundedly invertible for all $(\theta,t)\in[-\pi,\pi]\times[0,1]$. The corresponding Toeplitz operators therefore form a continuous Fredholm homotopy, and hence
	$$
	\operatorname{ind}(\lambda I - \widetilde{\mathbb{K}}^*_\alpha) = \operatorname{ind} \bigl( \lambda I - T(\mathcal Q_M\kappa_\alpha \mathcal Q_M) \bigr),
	$$
	where $T(\mathcal Q_M\kappa_\alpha \mathcal Q_M)$ is the Toeplitz operator generated by the symbol $\mathcal Q_M\kappa_\alpha \mathcal Q_M$. 	Since $\mathcal Q_M\kappa_\alpha(\theta)\mathcal Q_M$ vanishes on $\ker(\mathcal Q_M)$ and leaves $\operatorname{Im}(\mathcal Q_M)$ invariant, the latter operator decomposes into the matrix-valued Toeplitz operator with symbol $\lambda I_M - \kappa^{(M)}_\alpha(\theta)$ and the invertible operator $\lambda I$ on the complementary sequence space. Hence, by the classical index formula for matrix-valued symbols
	\cite{GohbergKrein1958},
	\[
	\operatorname{ind}\bigl(\lambda I-\widetilde{\mathbb K}^*_\alpha\bigr)
	=-\operatorname{wind}\det\bigl(\lambda I_M-\kappa^{(M)}_\alpha(\theta)\bigr),
	\]
	the winding number being taken as $\theta$ traverses $[-\pi,\pi]$. Writing it as
	the integral of the logarithmic derivative yields the stated formula.
\end{proof}

\begin{remark}
	The truncation size $M$ may depend on $\lambda$. In particular, no uniform choice of $M$ over the whole Fredholm set is asserted.	
\end{remark}

\begin{corollary}[Index as eigenvalue winding]\label{cor:index-winding}
	Under the hypotheses of Theorem~\ref{thm:fredholm-index}, let $\Gamma_1,\dots,\Gamma_L$ be the closed contours traced out in $\mathbb{C}$ by the eigenvalues of $\kappa^{(M)}_\alpha(\theta)$, counted with multiplicity, as $\theta$ traverses $[-\pi,\pi]$. Then
	\[
	\operatorname{ind}(\lambda I-\widetilde K^*_\Gamma)
	=-\sum_{k=1}^{L}\operatorname{wind}_\lambda\Gamma_k ,
	\]
	the winding numbers being taken about $\lambda$.
\end{corollary}

\subsection{Main result 2: the spectrum in $H^{-1/2}_0(\Gamma)$}\label{sec:physical space}
The physical energy space corresponds to $s = -1/2$ and $p = 2$, yielding the exact scaling exponent $\alpha = d/2$. At this scaling, the spectral picture simplifies considerably. 

\begin{theorem}[The spectrum in the energy space]\label{thm:physical-spectrum}
	Let $\widetilde{K}^*_\Gamma$ act on the mean-zero physical energy space $H^{-1/2}_0(\Gamma)$. Then
	$$
	\sigma(\widetilde{K}^*_\Gamma) = \sigma_{\mathrm{ess}}(\widetilde{K}^*_\Gamma) \cup \{\lambda_k\}_{k\ge 1} \subset \mathbb{R},
	$$
	where
	$$
	\sigma_{\mathrm{ess}}(\widetilde{K}^*_\Gamma) = \bigcup_{\theta\in[-\pi,\pi]} \sigma\bigl(\kappa_{d/2}(\theta)\bigr) \subset \mathbb{R},
	$$
	and $\{\lambda_k\}_{k\ge 1}$ denotes the possibly finite or empty collection of isolated eigenvalues outside the essential spectrum. Each $\lambda_k$ has finite multiplicity, and these eigenvalues may accumulate only at the essential spectrum.
\end{theorem} 

\begin{proof}
	By Lemma~\ref{lem: Symmetrisation of K toeplitz} the operator
	$\widetilde{\mathbb K}^*_{d/2}$ is self-adjoint for the inner product
	$\langle\cdot,\cdot\rangle_{\mathbb J}$, which is equivalent to the standard
	one. Hence
	$\sigma(\tilde K^*_\Gamma)=\sigma(\widetilde{\mathbb K}^*_{d/2})\subset\mathbb R$.
	
	The characterisation of the essential spectrum is
	Theorem~\ref{thm:essential-spectrum}. Since
	$\sigma_{\mathrm{ess}}(\tilde K^*_\Gamma)$ is a compact subset of $\mathbb R$,
	its complement in $\mathbb C$ is connected. Because
	$\lambda I-\tilde K^*_\Gamma$ is invertible for every non-real $\lambda$, and
	the Fredholm index is constant on connected components of the Fredholm set,
	\[
	\operatorname{ind}(\lambda I-\tilde K^*_\Gamma)=0,
	\qquad \lambda\notin\sigma_{\mathrm{ess}}(\tilde K^*_\Gamma).
	\]
	Consequently every spectral point outside the essential spectrum is an isolated
	real eigenvalue of finite multiplicity, with possible accumulation only at the
	essential spectrum.
\end{proof}

One can further characterise the eigenvalues $\lambda_k$ through an
operator-valued Wiener--Hopf factorisation. The relevant theory, for symbols that
are compact perturbations of the identity and holomorphic on an annulus containing
$\mathbb T$, goes back to Gohberg and Leiterer
\cite{GohbergLeiterer1973}; see
Appendix~\ref{app:wiener-hopf} for the statements used here.
For $\lambda \notin \sigma_{\mathrm{ess}}(\widetilde{K}^*_\Gamma)$, define the operator-valued function
$$
a_\lambda(z) = I - \lambda^{-1} \sum_{m\in\mathbb{Z}} r^{md/2} \widetilde{K}^*_m z^m.
$$	
By the discussion following Definition~\ref{def:symbol}, this Laurent series converges in the operator norm of $\mathcal L(H^{-1/2}_0(\Gamma_0))$ on an annulus containing the unit circle $\mathbb T$. On $\mathbb{T}$, with $z = e^{i\theta}$,
$$
a_\lambda(e^{i\theta}) = I - \lambda^{-1}\kappa_{d/2}(\theta).
$$
\begin{theorem}\label{thm:wiener-hopf}
	For every $\lambda\notin\sigma_{\mathrm{ess}}(\widetilde K^*_\Gamma)$, the function
	$a_\lambda$ admits a right Wiener--Hopf factorisation with respect to
	$\mathbb T$ in the sense of Definition~\ref{def:app-factorisation}, with fibre
	space $H^{-1/2}_0(\Gamma_0)$,
	\[
	a_\lambda(z)=a_{\lambda,-}(z)\,\Lambda_\lambda(z)\,a_{\lambda,+}(z),
	\qquad
	\Lambda_\lambda(z)=
	\begin{pmatrix}
		z^{\nu_1} & & &\\
		& \ddots & &\\
		& & z^{\nu_J} &\\
		& & & I
	\end{pmatrix},
	\]
	where the middle factor acts on the direct sum
	$H^{-1/2}_0(\Gamma_0)=\operatorname{Im}(\Pi_1)\oplus\cdots\oplus
	\operatorname{Im}(\Pi_J)\oplus\operatorname{Im}(I-\Pi)$, where $\Pi_1,\dots,\Pi_J$ are rank-one idempotents with $\Pi_i\Pi_j=0$ for $i\neq j$, and
	$\nu_1,\dots,\nu_J\in\mathbb Z\setminus\{0\}$ are the non-zero partial indices.
	Since the Fredholm index vanishes outside the essential spectrum, these satisfy $\sum_{j=1}^{J}\nu_j=0$. Moreover
	\[
	\lambda \text{ is an eigenvalue of } \widetilde K^*_\Gamma
	\iff
	\text{the factorisation of $a_\lambda$ is non-canonical, i.e. }J\ge1,
	\]
	and in that case
	\[
	\dim\ker(\lambda I-\widetilde K^*_\Gamma)=\sum_{\nu_j<0}(-\nu_j).
	\]
	The corresponding eigenfunctions are square summable across the components of the chain.
\end{theorem}
\begin{proof}
	Since $0\in\sigma_{\mathrm{ess}}(\widetilde K^*_\Gamma)$, every
	$\lambda\notin\sigma_{\mathrm{ess}}(\widetilde K^*_\Gamma)$ satisfies $\lambda\neq0$.
	Because the blocks $\widetilde K^*_m$ are compact, $a_\lambda-I$ takes values in
	$\mathcal K\bigl(H^{-1/2}_0(\Gamma_0)\bigr)$ throughout the annulus of
	convergence, and by Theorem~\ref{thm:essential-spectrum} the operator
	$a_\lambda(e^{i\theta})$ is boundedly invertible for every
	$\theta$. Theorem~\ref{thm:app-factorisation} therefore yields the stated
	factorisation, with finitely many non-zero partial indices.
	
	By construction $\lambda I-\widetilde{\mathbb K}^*_{d/2}=\lambda\,T(a_\lambda)$, where $T(a_\lambda)$ is the Toeplitz operator generated by the symbol $a_\lambda$ (as defined in Appendix~\ref{app:toeplitz}). Hence, Theorem~\ref{thm:app-partial-indices} gives 
\[
\dim\ker(\lambda I-\widetilde{\mathbb K}^*_{d/2})=\sum_{\nu_j<0}(-\nu_j),
\qquad
\dim\operatorname{coker}(\lambda I-\widetilde{\mathbb K}^*_{d/2})
=\sum_{\nu_j>0}\nu_j,
\]
and $\sum_{j}\nu_j=0$ because the index vanishes outside the essential spectrum.
In particular the kernel is non-trivial precisely when some $\nu_j$ is negative,
which by $\sum_j\nu_j=0$ happens exactly when $J\ge1$. Since the eigenvectors lie in
$\ell^2\bigl(\mathbb N,H^{-1/2}_0(\Gamma_0)\bigr)$ they are square summable across
the components of the chain.
\end{proof}

\section{Asymptotics for large finite chains (main result 3)}\label{sec:asymptotics}
We now turn to the truncated system consisting of a large but finite number of
inclusions, $N\gg1$. Throughout this section we work on $H^{-1/2}_0(\Gamma)$, at the energy scaling
$\alpha=d/2$. 

Let $Q_N$ denote the projection onto the first $N$ components of
the sequence space. The projected Neumann--Poincar\'e operator for the finite
chain of $N$ inclusions is represented by the truncated block-Toeplitz operator
\[
\mathbb A_N:=Q_N\widetilde{\mathbb K}^*_{d/2}Q_N =T_N\bigl(\kappa_{d/2}\bigr),
\]
where $T_N(a)$ is the $N\times N$ block-Toeplitz operator generated by a given symbol $a$ (as defined in Appendix~\ref{app:szego}). We are interested in the asymptotic distribution of the eigenvalues of $\mathbb{A}_N$ filling the essential spectrum bands of $\widetilde{\mathbb{K}}^*_{d/2}$. These eigenvalues correspond exactly to the eigenvalues of the projected Neumann--Poincar\'e operator on the finite chain $\widetilde{K}^*_{\Gamma^{(N)}}$.

The natural tool for this analysis is the Szeg\H o limit trace theorem for block-Toeplitz operators which is classical
for matrix-valued self-adjoint symbols \cite{MirandaTilli2000,Widom1974,BottcherSilbermann2006}. Here the symbol takes values in an infinite-dimensional
space and is not self-adjoint. The symbol is however symmetrisable by a $\theta$-dependent
inner product, so we prove a version of Szeg\H o's theorem adapted to our context. This is done in Appendix~\ref{app:szego},
Theorem~\ref{thm:app-szego}.

The following three lemmas gives conditions to apply Theorem~\ref{thm:app-szego} to $\mathbb A_N$. First, recall
$\mathbb J =-D^{-1/2}\widetilde{\mathbb S}_{\frac{d-1}{2}}D^{-1/2}$ of
Lemma~\ref{lem: Symmetrisation of K toeplitz}, which is block-Toeplitz because
diagonal conjugation preserves the Toeplitz structure. We denote by $\tau$ the symbol of $\mathbb J$.

\begin{lemma}[Exact symmetry of the finite chain]\label{lem:finite-symmetry}
	Write $\ell^2_N:=\ell^2\bigl(\{0,\dots,N-1\},H^{-1/2}_0(\Gamma_0)\bigr)$. For
	every $N\ge1$,
	\[
	T_N(\tau)\,\mathbb A_N=\mathbb A_N^*\,T_N(\tau),
	\]
	and there are constants $0<c\le C$ independent of $N$ with
	\[
	c\,\|u\|^2_{\ell^2_N}\le\langle T_N(\tau)u,u\rangle\le C\,\|u\|^2_{\ell^2_N},
	\qquad u\in\ell^2_N .
	\]
	Consequently $\mathbb A_N$ is compact and self-adjoint for the inner product
	$\langle T_N(\tau)\cdot,\cdot\rangle$ on $\ell^2_N$, with real eigenvalues.
\end{lemma}

\begin{proof}
	The truncated boundary $\Gamma^{(N)}$ is a finite union of disjoint smooth
components, so the classical symmetrisation identity applies:
\[
\widetilde S_{\Gamma^{(N)}}\widetilde K^*_{\Gamma^{(N)}}
=\widetilde K_{\Gamma^{(N)}}\widetilde S_{\Gamma^{(N)}} .
\]
Since the kernels of the layer potentials do not see the discarded components,
$\widetilde S_{\Gamma^{(N)}}=Q_N\widetilde S_\Gamma Q_N$,
$\widetilde K^*_{\Gamma^{(N)}}=Q_N\widetilde K^*_\Gamma Q_N$ and
$\widetilde K_{\Gamma^{(N)}}=Q_N\widetilde K_\Gamma Q_N$. Every factor of
$\mathcal M_\alpha=D^\alpha\mathcal P\mathcal T$ acts diagonally in the component
index, as does $D^{\pm1/2}$, so all of them commute with $Q_N$; conjugating the
identity above by $\mathcal M_\alpha$ with $\alpha=\frac{d-1}{2}$ therefore gives
\[
Q_N\widetilde{\mathbb S}_{\alpha}Q_N\;Q_N\widetilde{\mathbb K}^*_{\alpha}Q_N
=Q_N\widetilde{\mathbb K}_{\alpha}Q_N\;Q_N\widetilde{\mathbb S}_{\alpha}Q_N .
\]
Multiplying by $-D^{-1/2}$ on both sides, inserting $D^{-1/2}D^{1/2}=I$ between
the factors and using
$\widetilde{\mathbb K}^*_{d/2}=D^{1/2}\widetilde{\mathbb K}^*_{\alpha}D^{-1/2}$
together with $\widetilde{\mathbb K}_{\alpha}
=\bigl(\widetilde{\mathbb K}^*_{\alpha}\bigr)^{*}$, as in the proof of
Lemma~\ref{lem: Symmetrisation of K toeplitz}, yields
\[
T_N(\tau)\,\mathbb A_N=\mathbb A_N^{*}\,T_N(\tau).
\]
	
	For the two-sided bound, note that $T_N(\tau)=Q_N\mathbb JQ_N$, so that
	$\langle T_N(\tau)u,u\rangle=\langle\mathbb Ju,u\rangle$ for every $u$ supported
	on the first $N$ components. The bound for $\mathbb J$ established in
	Lemma~\ref{lem: Symmetrisation of K toeplitz} therefore applies with the same
	constants, which are in particular independent of $N$.
\end{proof}

We will also need a symmetrisation result at the level of symbols.
\begin{lemma}[Energy symmetry of the fibre]\label{lem:fibre-symmetry}
	For every $\theta\in[-\pi,\pi]$ the operator
	$\tau(\theta):H^{-1/2}_0(\Gamma_0)\to H^{1/2}_0(\Gamma_0)$ is symmetric for the
	$L^2(\Gamma_0)$ pairing and uniformly coercive:
	\begin{equation}\label{eq:symbol-coercive}
		c\,\|v\|^2_{H^{-1/2}(\Gamma_0)}\;\le\;\langle\tau(\theta)v,v\rangle\;\le\;
		C\,\|v\|^2_{H^{-1/2}(\Gamma_0)},
		\qquad v\in H^{-1/2}_0(\Gamma_0),
	\end{equation}
	with $c,C>0$ independent of $\theta$. Moreover
	\begin{equation}\label{eq:symbol_identity}
		\bigl\langle\tau(\theta)\kappa_{d/2}(\theta)u,v\bigr\rangle
		=\bigl\langle\tau(\theta)u,\kappa_{d/2}(\theta)v\bigr\rangle,
		\qquad u,v\in H^{-1/2}_0(\Gamma_0),
	\end{equation}
	so $\kappa_{d/2}(\theta)$ is self-adjoint on $H^{-1/2}_0(\Gamma_0)$ for the
	energy inner product $\langle\tau(\theta)\cdot,\cdot\rangle$, which is uniformly
	equivalent to the $H^{-1/2}$ one.
\end{lemma}
\begin{proof}
	Let $t_m:=-r^{m(d-2)/2}\widetilde S_m$, so that
	$\tau(\theta)=\sum_{m\in\mathbb Z}t_me^{im\theta}$. The decay bounds of
	Lemma~\ref{lem:decay-SL-blocks} give $\|t_m\|\lesssim r^{|m|d/2}$ for
	$m$ of either sign, so the series converges absolutely in
	$\mathcal L\bigl(H^{-1/2}_0(\Gamma_0),H^{1/2}_0(\Gamma_0)\bigr)$, and the
	adjoint relation $\widetilde S^*_m=r^{-m(d-2)}\widetilde S_{-m}$ gives
	$t^*_m=t_{-m}$, whence $\tau(\theta)^*=\tau(\theta)$.
	
	For \eqref{eq:symbol-coercive}, fix $\theta$ and $v$, and set
	$u^{(N)}_j:=N^{-1/2}e^{-ij\theta}v$ for $0\le j<N$ and $u^{(N)}_j:=0$ otherwise.
	Then $\|u^{(N)}\|_{\ell^2(\mathbb N,H^{-1/2}_0(\Gamma_0))}
	=\|v\|_{H^{-1/2}(\Gamma_0)}$ for every $N$ and every $\theta$, and
	\[
	\bigl\langle\mathbb Ju^{(N)},u^{(N)}\bigr\rangle
	=\sum_{|m|<N}\Bigl(1-\frac{|m|}N\Bigr)e^{im\theta}\langle t_mv,v\rangle
	\;\longrightarrow\;\langle\tau(\theta)v,v\rangle,
	\]
	where the limit follows by dominated convergence from the absolute summability of
	the coefficients. The two-sided bound for $\mathbb J$ established in
	Lemma~\ref{lem: Symmetrisation of K toeplitz} therefore passes to the limit with
	the same constants, which are independent of $\theta$ because the test sequence
	has the same norm for every $\theta$.
	
	For \eqref{eq:symbol_identity}, write $\kappa_m:=r^{md/2}\widetilde K^*_m$ for
the coefficients of $\kappa_{d/2}$. Equating the $(j,k)$ entries on both sides of
the Plemelj identity
$\mathbb J\widetilde{\mathbb K}^*_{d/2}
=\bigl(\widetilde{\mathbb K}^*_{d/2}\bigr)^*\mathbb J$ of
Lemma~\ref{lem: Symmetrisation of K toeplitz} gives
\[
\sum_{i\ge0}t_{j-i}\,\kappa_{i-k}
=\sum_{i\ge0}\kappa^*_{i-j}\,t_{i-k},
\qquad j,k\ge0 .
\]
Both sides are absolutely convergent, and letting $j,k\to\infty$ with $j-k=m$
fixed extends the summation to $i\in\mathbb Z$, so that
$\sum_{i+i'=m}t_i\kappa_{i'}=\sum_{i+i'=m}\kappa^*_it_{i'}$ for every
$m\in\mathbb Z$. Resumming against $e^{im\theta}$ yields
$\tau(\theta)\kappa_{d/2}(\theta)=\kappa_{d/2}(\theta)^*\tau(\theta)$, which is
\eqref{eq:symbol_identity}. Together with \eqref{eq:symbol-coercive} this states
precisely that $\kappa_{d/2}(\theta)$ is symmetric, hence self-adjoint, for the
inner product $\langle\tau(\theta)\cdot,\cdot\rangle$.
\end{proof}

\begin{remark}
	The content of Lemma~\ref{lem:fibre-symmetry} is that $\kappa_{d/2}(\theta)$ is
	the Neumann--Poincar\'e operator of the fibre written in its natural energy
	pairing; its apparent non-symmetry is an artefact of the ambient $H^{-1/2}$ norm,
	exactly as in the classical case of a single smooth boundary.
\end{remark}

It remains to check the regularity hypothesis of Theorem~\ref{thm:app-szego},
namely membership in the two Schatten--Wiener classes of
Definition~\ref{def:app-schattenwiener}. The exponent $q>d-1$ is forced by the
diagonal coefficient.

\begin{lemma}[Schatten--Wiener class of the symbol]\label{lem:schatten-wiener}
	For every $q>d-1$,
	\[
	\kappa_{d/2}\in \mathcal W^1\bigl(\mathfrak S_q(H^{-1/2}_0(\Gamma_0))\bigr)
	\cap \mathcal W^1_\flat\bigl(\mathfrak S_2(H^{-1/2}_0(\Gamma_0))\bigr).
	\]
\end{lemma}

\begin{proof}
	By Lemma~\ref{lem:schatten-blocks} the weighted coefficients
	$r^{md/2}\tilde K^*_m$ have $\mathfrak S_q$ norm $O(r^{|m|d/2})$ for $m$ of either
	sign, which is summable against the weight $(1+|m|)$; this gives the first
	membership. The second follows in the same way at $q=2$, the bounds for $m\neq0$
	holding for every $q\ge1$ and the coefficient of order zero being exempt by
	definition of the flat class.
\end{proof}

\begin{theorem}[Spectral asymptotics]\label{thm:szego-chain}
	Let $\{\lambda^{(N)}_j\}_{j\ge1}$ be the non-zero eigenvalues of $\widetilde{K}^*_{\Gamma^{(N)}}$, repeated according to multiplicity, and let
	$\{\mu_j(\theta)\}_{j\ge1}$ be the non-zero eigenvalues of
	$\kappa_{d/2}(\theta)$. Then for every $f\in\mathcal F$ (given by Definition~\ref{def:app-testfunctions}),
	\begin{equation}\label{eq:szego}
		\lim_{N\to\infty}\frac1N\sum_{j\ge1}f\bigl(\lambda^{(N)}_j\bigr)
		=\frac1{2\pi}\int_{-\pi}^{\pi}\sum_{j\ge1}f\bigl(\mu_j(\theta)\bigr)\,d\theta .
	\end{equation}
\end{theorem}

\begin{proof}
	$\{\lambda^{(N)}_j\}_{j\ge1}$ correspond to the non-zero eigenvalues of
	$\mathbb A_N$. We have $\mathbb A_N=T_N(\kappa_{d/2})$, and from the three
	lemmas above, with $q>d-1$ in Lemma~\ref{lem:schatten-wiener}, the symbol
	$\kappa_{d/2}$ satisfies all the requirements of
	Theorem~\ref{thm:app-szego}. This gives directly \eqref{eq:szego}.
\end{proof}

\begin{corollary}[Eigenvalue counting function]
	\label{cor:density-of-states}
	Let
	$
	[a,b]\subset\mathbb{R}\setminus\{0\}
	$
	be a compact interval such that
	\[
	\left| \left\{ \theta\in[-\pi,\pi] : a\in\sigma\bigl(\kappa_{d/2}(\theta)\bigr) \text{ or } b\in\sigma\bigl(\kappa_{d/2}(\theta)\bigr) \right\} \right| = 0.
	\]
	Then
	\begin{equation}
		\label{eq:density-of-states}
		\lim_{N\to\infty} \frac{1}{N} \#\left\{ j:\lambda_j^{(N)}\in 	[a,b] \right\} = \frac{1}{2\pi} \int_{-\pi}^{\pi} N_{	[a,b]}\bigl(\kappa_{d/2}(\theta)\bigr) \,d\theta,
	\end{equation}
	where $N_{	[a,b]}(\kappa_{d/2}(\theta))$ denotes the number of eigenvalues of $\kappa_{d/2}(\theta)$ lying in $[a,b]$, counted with multiplicity.
\end{corollary}

\begin{proof}
	Since $[a,b]$ is separated from the accumulation point $0$, its characteristic function can be approximated from above and below by functions in $\mathcal{F}$ with continuous transitions from zero to $a$ and $b$. Applying Theorem~\ref{thm:szego-chain} to these approximations and letting the transition regions shrink to the endpoints $a$ and $b$ gives \eqref{eq:density-of-states}. The hypothesis on the endpoints ensures that the upper and lower limits coincide.
\end{proof}

\section{A computable example: concentric 2D annuli}\label{sec:annuli}

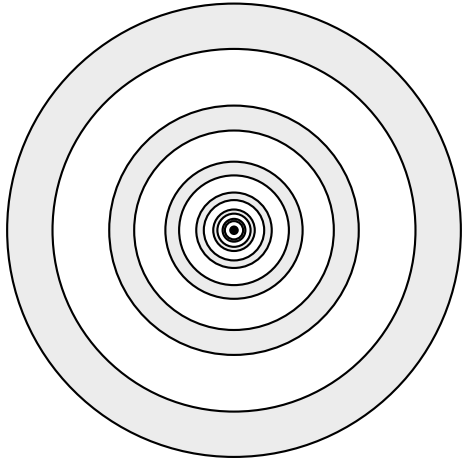
\begin{figure}[!htbp]
	\centering
	\begin{tikzpicture}[scale=3]
		
		\def\r{0.55}       
		\def\q{0.8}        
		\def\R{1.0}        
		\def\nAnnuli{5}    
		
		\foreach \j in {0,...,\nAnnuli} {
			\pgfmathsetmacro{\Rout}{\R*(\r)^(\j)}
			\pgfmathsetmacro{\Rin}{\q*\R*(\r)^(\j)}
			
			\filldraw[
			thick,
			fill=gray!15,
			even odd rule
			]
			(0,0) circle (\Rout)
			(0,0) circle (\Rin);
		}
		
		\pgfmathsetmacro{\rlast}{\R*(\r)^(\nAnnuli)}
		
		\fill (0,0) circle (0.6pt);
		
		
	\end{tikzpicture}
	\caption{A self-similar chain of concentric annuli. The reference annulus $\Gamma_0$ has outer radius $R_1$ and inner radius $R_2$. The non-overlap condition $r R_1 < R_2$ ensures that consecutive scaled components remain strictly disjoint.}
	\label{fig:annuli}
\end{figure}

We illustrate our main results by considering a self-similar chain of two-dimensional rings (as in Figure~\ref{fig:annuli}), for which the relevant operator-valued symbols can be computed explicitly. For simplicity, we focus on the Sobolev spaces $H_0^s(\Gamma)$.

\subsection{Geometry and function spaces}
Let the fundamental domain $\Omega_0\subset\mathbb R^2$ be an annulus centred at the origin, with outer radius $R_1$ and inner radius $R_2$, where $R_1>R_2>0$. Its boundary
\[
\Gamma_0=\partial\Omega_0
\]
consists of two disjoint concentric circles,
\begin{equation*}
	\Gamma_0
	=
	\Gamma_{\mathrm{out}}\cup\Gamma_{\mathrm{in}},
	\qquad
	\Gamma_{\mathrm{out}}=\{|x|=R_1\},
	\qquad
	\Gamma_{\mathrm{in}}=\{|x|=R_2\}.
\end{equation*}
The outward normal to the annulus points radially outward on $\Gamma_{\mathrm{out}}$ and radially inward on $\Gamma_{\mathrm{in}}$.

Recall that the self-similar geometry is generated by the pure scaling maps $\eta_j(x)=r^jx$, $j\in\mathbb N$, and that consecutive rings remain strictly disjoint provided $rR_1<R_2$.

We parametrise the two components of $\Gamma_0$ by their common polar angle $\phi\in[0,2\pi)$. Any distribution $v$ on $\Gamma_0$ may then be decomposed into its Fourier series
\begin{equation*}
	v(\phi)=\sum_{n\in\mathbb Z}\mathbf c_ne^{in\phi},
	\qquad
	\mathbf c_n=\begin{pmatrix}a_n\\ b_n\end{pmatrix}\in\mathbb C^2,
\end{equation*}
where $a_n$ and $b_n$ are the Fourier coefficients on $\Gamma_{\mathrm{out}}$ and $\Gamma_{\mathrm{in}}$, and, up to constants depending on $R_1$, $R_2$ and $s$ but not on $n$,
\begin{equation*}
	\|v\|_{H^s(\Gamma_0)}^2\asymp\sum_{n\in\mathbb Z}(1+|n|^2)^s\|\mathbf c_n\|_{\mathbb C^2}^2 .
\end{equation*}

The mean-zero condition should be imposed on the whole boundary of each annular cell, rather than separately on its two connected components. On the reference boundary it reads $\int_{\Gamma_0}v\,d\sigma=0$. Since $d\sigma=R_1\,d\phi$ on $\Gamma_{\mathrm{out}}$ and $d\sigma=R_2\,d\phi$ on $\Gamma_{\mathrm{in}}$, and all non-zero Fourier modes integrate to zero, this constrains only the constant mode, $R_1a_0+R_2b_0=0$. Introducing the aspect ratio 
$$\rho:=\frac{R_2}{R_1}\in(0,1),$$
the admissible zero-mode subspace is
\begin{equation*}
	L_0:=\left\{\begin{pmatrix}a_0\\ b_0\end{pmatrix}\in\mathbb C^2: a_0+\rho b_0=0\right\}
	=\operatorname{span}\left\{\begin{pmatrix}-\rho\\1\end{pmatrix}\right\},
\end{equation*}
so that
\begin{equation*}
	H^s_0(\Gamma_0)=\left\{v=\sum_{n\in\mathbb Z}\mathbf c_ne^{in\phi}\;:\;\mathbf c_0\in L_0,\ \sum_{n\in\mathbb Z}(1+|n|^2)^s\|\mathbf c_n\|_{\mathbb C^2}^2<\infty\right\}.
\end{equation*}

\subsection{The block-Toeplitz operator $\widetilde{\mathbb K}^*_\alpha$ for the annuli chain}

We consider the scaled Neumann--Poincar\'e operator $\widetilde{\mathbb K}^*_\alpha$ acting on $\ell^2\bigl(\mathbb N,H^s_0(\Gamma_0)\bigr)$ where $s=\frac12-\alpha$.
As established in Section~\ref{sec:block-Toeplitz}, it has the block-Toeplitz representation
\[
\widetilde{\mathbb K}^*_\alpha
=
\begin{pmatrix}
	\tilde K^*_0
	&
	r^{-\alpha}\tilde K^*_{-1}
	&
	r^{-2\alpha}\tilde K^*_{-2}
	&
	\cdots
	\\
	r^\alpha\tilde K^*_1
	&
	\tilde K^*_0
	&
	r^{-\alpha}\tilde K^*_{-1}
	&
	\cdots
	\\
	r^{2\alpha}\tilde K^*_2
	&
	r^\alpha\tilde K^*_1
	&
	\tilde K^*_0
	&
	\cdots
	\\
	\vdots&\vdots&\vdots&\ddots
\end{pmatrix},
\]
where the local blocks $\tilde K_m^*=P_{0,0}K^*_m$ map onto the mean-zero subspace $H^s_0(\Gamma_0)$, with the unprojected action on the reference cell given by
\[
(K_m^*\psi)(\tilde x)
=
\int_{\Gamma_0}
\frac{\partial G}
{\partial\nu_{\eta_m(\tilde x)}}
\bigl(\eta_m(\tilde x),\tilde y\bigr)
\psi(\tilde y)\,d\sigma_0(\tilde y).
\]

Because the concentric geometry and the Neumann--Poincar\'e kernel are invariant under rotations, each $\tilde K_m^*$ commutes with the rotation group. Consequently, every angular Fourier subspace is invariant, and the operator diagonalises with respect to the Fourier index.

For every $n\neq0$, the action on the $n$-th Fourier sector is represented by a matrix
\[
A_{n,m}\in\mathbb C^{2\times2},
\]
so that
\[
\tilde K_m^*
\bigl(\mathbf c_ne^{in\phi}\bigr)
=
A_{n,m}\mathbf c_ne^{in\phi}.
\]
The constant Fourier sector is treated separately: after projection onto $H^s_0(\Gamma_0)$, it gives a one-dimensional operator on $L_0$. Accordingly, the Fourier decomposition of $\tilde K_m^*$ takes the form
\[
\tilde K_m^*
=
A_{0,m}^{(0)}
\oplus
\bigoplus_{n\neq0}A_{n,m},
\]
where $A_{0,m}^{(0)}$ denotes the restriction of the projected zero-mode action to $L_0$.

For the non-zero modes, no Fourier frequencies are mixed, and on each sector the Sobolev weight $(1+|n|^2)^s$ acts as a positive scalar multiple of the identity, so it commutes with $A_{n,m}$ and cancels between input and output. The $2\times2$ multipliers $A_{n,m}$, and their norms, are therefore (and as expected) independent of the Sobolev exponent $s$.

We now compute these matrices explicitly. The Neumann--Poincar\'e kernel in two dimensions is
\[
\frac{\partial G}{\partial \nu_x}(x, y)
=
\frac{1}{2\pi}
\frac{\langle x-y,\nu_x\rangle}{|x-y|^2}.
\]

For the intra-ring interaction $m=0$, the self-interaction of either circle vanishes on every non-constant Fourier mode. For the interaction from the inner circle to the outer circle, one obtains, for $n\neq0$,
\begin{align*}
	\frac{1}{2\pi}
	\int_0^{2\pi}
	\frac{R_1-R_2\cos\theta}
	{R_1^2+R_2^2-2R_1R_2\cos\theta}
	e^{-in\theta}R_2\,d\theta
	=
	\frac12\rho^{|n|+1}.
\end{align*}
Taking into account the opposite orientation of the normal on the inner boundary, the four pairings give
\begin{equation*}
	A_{n,0}
	=
	\frac12
	\begin{pmatrix}
		0&\rho^{|n|+1}\\[1mm]
		\rho^{|n|-1}&0
	\end{pmatrix},
	\qquad n\neq0.
\end{equation*}
In particular, the two eigenvalues of the isolated-annulus block are
$
\pm\frac12\rho^{|n|}.
$

For $m>0$, the target copy $\eta_m(\Gamma_0)$ lies strictly inside the source copy $\Gamma_0$. Evaluating the four circle-to-circle interactions gives
\begin{equation*}
	A_{n,m}
	=
	\frac12r^{m(|n|-1)}
	\begin{pmatrix}
		-1&-\rho^{1-|n|}\\[1mm]
		\rho^{|n|-1}&1
	\end{pmatrix},
	\qquad m>0,\quad n\neq0.
\end{equation*}

For $m<0$, the target copy lies outside the source copy, and the corresponding calculation gives
\begin{equation*}
	A_{n,m}
	=
	\frac12r^{|m|(|n|+1)}
	\begin{pmatrix}
		1&\rho^{|n|+1}\\[1mm]
		-\rho^{-|n|-1}&-1
	\end{pmatrix},
	\qquad m<0,\quad n\neq0.
\end{equation*}

Consequently, for each fixed non-zero Fourier mode $n$, the $(j,k)$-th block of the fully scaled self-similar operator acts through
\[
r^{(j-k)\alpha}A_{n,j-k}.
\]
Thus, apart from the one-dimensional zero-mode sector, the infinite-dimensional problem decomposes exactly into independent $2\times2$ block-Toeplitz problems, one for each non-zero angular Fourier mode. The zero mode gives an additional scalar block-Toeplitz problem on the spaces $L_0$, which can be computed separately.

\subsection{The spectrum}
We begin by computing the symbol of the block-Toeplitz operator. Recall that
\[
\kappa_\alpha(\theta)
=
\sum_{m\in\mathbb Z}
r^{m\alpha}\tilde K_m^*e^{im\theta},
\qquad
\theta\in[-\pi,\pi].
\]
Since each block $\tilde K_m^*$ diagonalises in the angular Fourier basis, so does the symbol:
\begin{equation*}
	\kappa_\alpha(\theta)
	=
	\kappa_{\alpha,0}(\theta)
	\oplus
	\bigoplus_{n\neq0}\kappa_{\alpha,n}(\theta).
\end{equation*}
The zero Fourier sector is one-dimensional. On this sector the density is radially symmetric with vanishing total charge, so its potential vanishes outside the annulus and the interactions between distinct cells vanish, while $\tilde K_0^*$ acts as multiplication by $-1/2$. Hence
\begin{equation*}
	\kappa_{\alpha,0}(\theta)=-\frac12.
\end{equation*}
This flat branch is the classical eigenvalue $-\tfrac12$ associated with the holes of $\Omega$, here of infinite multiplicity because the chain has infinitely many components, and therefore belongs to the essential spectrum.

For $n\neq0$, using the matrices obtained in the previous subsection,
\begin{equation*}
	\kappa_{\alpha,n}(\theta)
	=
	A_{n,0}
	+
	\sum_{m=1}^{\infty}
	r^{m\alpha}A_{n,m}e^{im\theta}
	+
	\sum_{m=1}^{\infty}
	r^{-m\alpha}A_{n,-m}e^{-im\theta}.
\end{equation*}
Both series converge geometrically, with ratios $r^{|n|-1+\alpha}$ and $r^{|n|+1-\alpha}$, the exponents being positive for $|n|\ge1$ precisely because $0<\alpha<2$. Summing gives
\begin{align}
	\kappa_{\alpha,n}(\theta)
	=
	\frac12\Bigg[
	&
	\begin{pmatrix}
		0&\rho^{|n|+1}\\
		\rho^{|n|-1}&0
	\end{pmatrix}
	+
	\frac{r^{|n|-1+\alpha}e^{i\theta}}
	{1-r^{|n|-1+\alpha}e^{i\theta}}
	\begin{pmatrix}
		-1&-\rho^{1-|n|}\\
		\rho^{|n|-1}&1
	\end{pmatrix}
	\nonumber\\
	&+
	\frac{r^{|n|+1-\alpha}e^{-i\theta}}
	{1-r^{|n|+1-\alpha}e^{-i\theta}}
	\begin{pmatrix}
		1&\rho^{|n|+1}\\
		-\rho^{-|n|-1}&-1
	\end{pmatrix}
	\Bigg].
	\label{eq:ring-symbol-general}
\end{align}
Since the matrices depend on $n$ only through $|n|$,
\[
\kappa_{\alpha,-n}(\theta)
=
\kappa_{\alpha,n}(\theta).
\]

\paragraph{The essential spectrum in $H^s_0(\Gamma)$.}
In dimension $d=2$, the scaling exponent associated with $H^s_0(\Gamma)$ is
\[
\alpha=\frac12-s.
\]
For each Fourier mode, the corresponding $2\times2$ matrix has zero trace, and hence its eigenvalues occur in opposite pairs. We denote them by
\[
\mu^\pm_{\alpha,n}(\theta)
=
\pm\mu_{\alpha,n}(\theta),
\qquad
\theta\in[-\pi,\pi].
\]
The zero mode contributes the constant flat branch 
$
\mu_{\alpha,0}(\theta)=-\frac12. 
$
Defining the spectral bands
\[
\mathcal B^\pm_{\alpha,n}
:=
\left\{
\mu^\pm_{\alpha,n}(\theta):
\theta\in[-\pi,\pi]
\right\},
\qquad n\geq1,
\]
and noting that $0\in\sigma\bigl(\kappa_\alpha(\theta)\bigr)$ for every $\theta$, each fibre operator being compact, the general essential-spectrum characterisation gives
\begin{equation*}
	\sigma_{\mathrm{ess}}
	\bigl(\widetilde{K}_\Gamma^*\bigr)
	=
	\left\{-\frac12\right\}
	\cup
	\{0\}
	\cup
	\bigcup_{n\geq1}
	\left(
	\mathcal B^+_{\alpha,n}
	\cup
	\mathcal B^-_{\alpha,n}
	\right),
\end{equation*}
the bands accumulating at $0$ as $n\to\infty$.

\begin{figure}[!htbp]
	\centering
			\includegraphics[width=\textwidth]{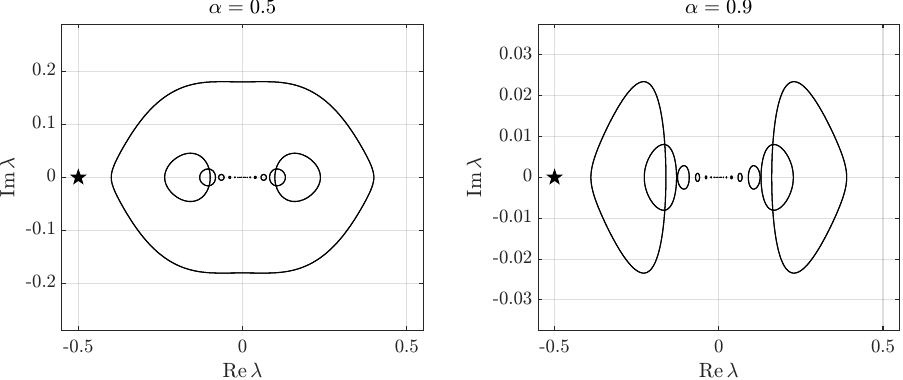}%
\caption{Analytical essential spectrum for the infinite chain of concentric annuli across different functional-space scalings $\alpha$. The horizontal and vertical axes denote the real and imaginary parts of the spectral values, respectively. The spectrum consists of a constant flat branch at $\lambda = -1/2$ (arising from the zero Fourier mode) and continuous symmetric bands $B_{\alpha,n}^\pm$ generated by the non-zero angular modes. Geometric parameters: $r = 0.2$ and $\rho = 0.6$.}
	\label{fig:essential-spectrum-annuli}
\end{figure}

These can be obtained explicitly from the eigenvalues of the $2 \times 2$ matrices $\kappa_{\alpha,n}(\theta)$. Their computation is illustrated in Figure~\ref{fig:essential-spectrum-annuli}. Notably, the closed loops formed by the continuous spectrum visually partition the complex plane. In accordance with Corollary~\ref{cor:index-winding}, the Fredholm index at a point $\lambda$ outside the bands is minus the total winding number of these curves about $\lambda$, so the components on which the winding fails to cancel carry a non-zero index.

\begin{figure}[!htbp]
	\centering
	\includegraphics[width=0.7\textwidth] {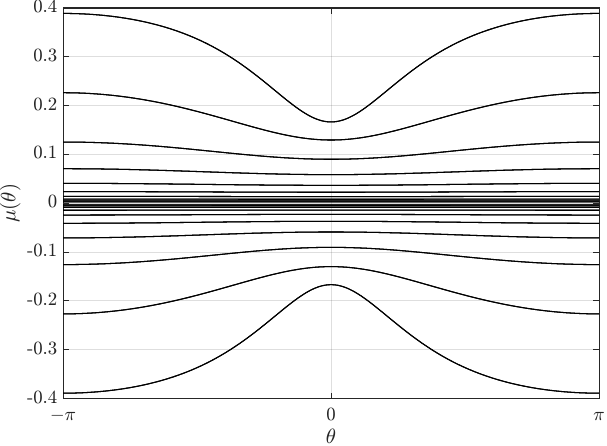}%
	\caption{Exact analytical spectral bands for the chain of annuli in the physical energy space ($\alpha = 1$). The vertical axis represents the eigenvalues of $\kappa_{1,n}(\theta)$ as a function of the Floquet phase $\theta \in [-\pi, \pi]$. Geometric parameters: $r = 0.2$ and $\rho = 0.6$.}
	\label{fig:dispersion-annuli}
\end{figure}

\paragraph{The spectrum in $H^{-1/2}_0(\Gamma)$.}
We now specialise to the physical energy space, for which $\alpha=d/2=1$. The zero Fourier mode gives the flat branch $-1/2$, so it remains to consider $n\neq0$. For a fixed angular mode, define
\[
a_{\lambda,n}(z)
:=
\lambda I-
\sum_{m\in\mathbb Z}r^m A_{n,m}z^m ,
\]
which is $\lambda$ times the restriction of $a_\lambda$ of Section~\ref{sec:physical space} to the $n$-th Fourier sector; the factor is immaterial for the factorisation, since $\lambda\neq0$. On the unit circle, writing $z = e^{i\theta}$, we have
\[
a_{\lambda,n}(e^{i\theta}) = \lambda I - \kappa_{1,n}(\theta).
\]
The resulting spectral bands for these physical energy states are plotted in Figure~\ref{fig:dispersion-annuli}.

\paragraph{Edge modes.} For $\lambda$ outside the $n$-th bulk bands, $a_{\lambda,n}$ is invertible on the unit circle and admits a Wiener--Hopf factorisation. A direct calculation gives
\begin{equation*}
	\det a_{\lambda,n}(z)
	=
	\frac{P_{\lambda,n}(z)}
	{4\rho^{2|n|}(z-r^{|n|})(r^{|n|}z-1)},
\end{equation*}
where
\begin{align}
	P_{\lambda,n}(z)
	={}&
	\rho^{2|n|}r^{|n|}(4\lambda^2-1)(z^2+1)+
	\left[
	\rho^{4|n|}
	+r^{2|n|}
	-4\rho^{2|n|}\lambda^2(1+r^{2|n|})
	\right]z .
	\label{eq:ring-characteristic-polynomial}
\end{align}
The coefficients of $z^2$ and $z^0$ in $P_{\lambda,n}$ coincide, so its two roots have product one. Denoting by $\zeta_{\lambda,n}$ the root lying inside the unit disk, the determinant therefore has a canonical scalar Wiener--Hopf factorisation.

The matrix factorisation can also be constructed explicitly. The pole of $a_{\lambda,n}$ inside the unit disk is at $z=r^{|n|}$; its residue there has rank one, with range spanned by
\[
u_n=
\begin{pmatrix}
	1\\
	-\rho^{-|n|-1}
\end{pmatrix},
\]
so a rank-one middle factor suffices. Let $v_{\lambda,n}^{\top}$ be a left null vector of $a_{\lambda,n}(\zeta_{\lambda,n})$. Provided $v_{\lambda,n}^{\top}u_n\neq0$, which as we show below is the only degenerate case, the rank-one projection
\[
\Pi_{\lambda,n}
=
\frac{u_nv_{\lambda,n}^{\top}}
{v_{\lambda,n}^{\top}u_n}
\]
gives the canonical factorisation
\[
a_{\lambda,n}(z)
=
a_{\lambda,n,-}(z)a_{\lambda,n,+}(z),
\]
with
\[
a_{\lambda,n,-}(z)
=
I+
\frac{r^{|n|}-\zeta_{\lambda,n}}
{z-r^{|n|}}\Pi_{\lambda,n},
\qquad
a_{\lambda,n,+}(z)
=
a_{\lambda,n,-}(z)^{-1}a_{\lambda,n}(z).
\]
The only possible obstruction to this canonical factorisation is
\[
\operatorname{adj}
a_{\lambda,n}(\zeta_{\lambda,n})\,u_n=0,
\]
where $\operatorname{adj}
a_{\lambda,n}(\zeta_{\lambda,n})$ denotes the adjugate of $a_{\lambda,n}(\zeta_{\lambda,n})$. Using \eqref{eq:ring-characteristic-polynomial}, this condition reduces to
\[
\lambda=0,
\qquad
\zeta_{\lambda,n}
=
\left(\frac{\rho^2}{r}\right)^{|n|}.
\]
Consequently, for $\lambda\neq0$ the factorisation is canonical throughout every gap of the $n$-th bulk spectrum. At $\lambda=0$ there are two different regimes:
\[
\begin{cases}
	r<\rho^2, & \text{the factorisation is canonical},\\[1mm]
	r>\rho^2, & \text{the factorisation is non-canonical with partial indices }
	(-1,1).
\end{cases}
\]
Both regimes occur, the disjointness condition being $r<\rho$. At the transition $r=\rho^2$, the corresponding zero reaches the unit circle and $0$ belongs to the bulk band.

In the non-canonical regime $r>\rho^2$, each non-zero Fourier mode therefore supports a zero-energy edge state. In fact, it can be written explicitly as
\[
\mathbf c_j
=
\left(\frac{\rho^2}{r}\right)^{j|n|}
\begin{pmatrix}
	-\rho^{1-|n|}\\
	1
\end{pmatrix},
\qquad j\geq0,
\]
which is square summable precisely when $r>\rho^2$. Since every non-zero mode contributes such a state, $\widetilde K^*_\Gamma$ has infinite-dimensional kernel when $r>\rho^2$, in contrast with a single annulus, whose mean-zero spectrum $\{\pm\tfrac12\rho^{|n|}\}\cup\{-\tfrac12\}$ omits the origin.
\begin{remark}\label{rem:modewise-edge-states}
	These states are not the isolated eigenvalues of
	Theorem~\ref{thm:wiener-hopf}, which requires
	$\lambda\notin\sigma_{\mathrm{ess}}$, whereas $0$ always in the essential spectrum. The
	mechanism of Theorem~\ref{thm:wiener-hopf} is, however, encountered modewise: within the $n$-th angular sector the symbol is a
	$2\times2$ matrix, $\lambda=0$ sits in a gap of its bulk spectrum, and the
	factorisation produces a genuine gap eigenvalue. The bands of the remaining
	modes then fill the origin, leaving the state embedded in the essential spectrum.
\end{remark}

\begin{figure}[!htbp]
	\centering
	\includegraphics[width=0.7\textwidth] {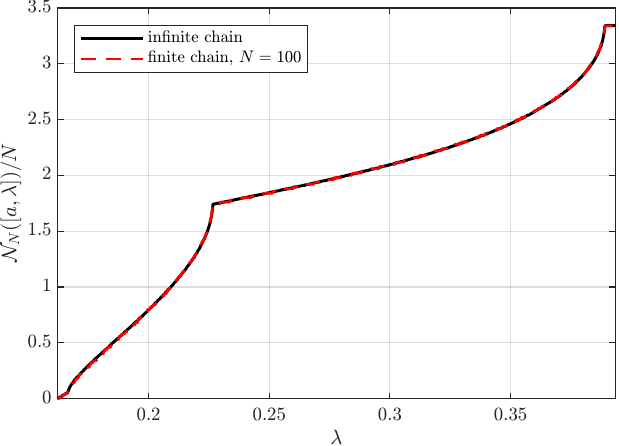}%
	\caption{Normalised eigenvalue counting function $\mathcal N_N([a,\lambda])/N$, defined in \eqref{eq:counting function I_lambda}, for a truncated chain of $N$ annuli, compared against its limit value computed from the infinite-chain symbol (Corollary~\ref{cor:density-of-states}). We show the results for $a = 0.16$ and $\lambda\in [0.16,0.4]$. Geometric parameters: $r = 0.2$ and $\rho = 0.6$.}
	\label{fig:counting-annuli}
\end{figure}

\paragraph{Asymptotics for the finite chain.}
We finally illustrate the large-chain asymptotics through the eigenvalue counting function. Let $\{\lambda_j^{(N)}\}$ denote the eigenvalues of the projected Neumann--Poincar\'e operator on a chain of $N$ annuli, and for a compact interval
$
[a,\lambda]\subset\mathbb R\setminus\{0\},
$
define
\begin{equation}\label{eq:counting function I_lambda}
	\mathcal N_N([a,\lambda])
	:=
	\#\left\{
	j:\lambda_j^{(N)}\in [a,\lambda]
	\right\}.
\end{equation}
Provided that the endpoints $a$ and $\lambda$ do not lie on a spectral band for a set of phases of positive measure, the eigenvalue counting formula gives
\begin{equation*}
	\lim_{N\to\infty}
	\frac{\mathcal N_N([a,\lambda])}{N}
	=
	\mathbf 1_{[a,\lambda]}\!\left(-\tfrac12\right)
	+
	\frac{1}{2\pi}
	\sum_{n\neq0}
	\int_{-\pi}^{\pi}
	\#\left\{
	\mu_{\alpha,n}^{\pm}(\theta)\in [a,\lambda]
	\right\}
	\,d\theta ,
\end{equation*}
the first term being the contribution of the flat branch of the zero mode, which supplies one eigenvalue per component. Here the notation inside the counting function includes both branches $\mu_{\alpha,n}^{+}$ and $\mu_{\alpha,n}^{-}$. See Figure~\ref{fig:counting-annuli}.

\section{A numerical example}\label{sec:disks}

We now illustrate our results for a self-similar chain of disks (as in Figure~\ref{fig:disks}), for which the operator-valued symbol is not available in closed form. We numerically approximate both the symbol $\kappa_\alpha(\theta)$ of the infinite chain and the spectrum of large finite truncations of the scaled Neumann--Poincar\'e operator $\widetilde{\mathbb{K}}^*_1$ (whose spectrum is real).

\subsection{Assembly of the scaled operator and symbol}

Let $\Gamma_0$ be the boundary of the reference disk and $\Gamma_j = r^j \Gamma_0$. We discretise $\Gamma_0$ using the equidistant Nystr\"om mesh
\[
t_\ell = \frac{2\pi\ell}{M}, \qquad \ell = 0, \ldots, M - 1.
\]
A density on $\Gamma_0$ is therefore represented by a vector in $\mathbb{C}^M$. Since the theoretical construction is carried out exclusively on the mean-zero space, we introduce the discrete mean-zero projection
\[
P_M = I - \frac{1}{M} \mathbf{1}\mathbf{1}^\top.
\] 

For each Toeplitz shift $m \in \mathbb{Z}$, let $\tilde{C}^{(M)}_m$ denote the Nystr\"om discretization of the reference block $\tilde{K}^*_m$. Before projection, its entries are given by
\[
[B^{(M)}_m]_{k,\ell} = \frac{2\pi}{M} \frac{\partial G}{\partial \nu_{\tilde{x}}} \bigl( \eta_m(x(t_k)), x(t_\ell) \bigr) J(t_\ell),
\]
where $J$ is the Jacobian of the parametrisation of $\Gamma_0$. We then set
\[
\tilde{C}^{(M)}_m = P_M B^{(M)}_m P_M.
\]
For $m \neq 0$ the two boundaries are disjoint and the kernel is smooth, so the periodic trapezoidal rule can be applied directly. The diagonal block $m = 0$ corresponds to the classical Neumann-Poincaré operator, projected to mean zero functions. Because the reference domain is a circle, $\tilde{C}^{(M)}_0 = 0$ on $\operatorname{Im}(P_M)$. For other domains, this can be computed using standard singularity removal for discretising integral operators with a singular kernel \cite{Ammari_book:18}.

For a finite chain of $N$ disks, the resulting discretisation of the scaled operator is the $NM \times NM$ block-Toeplitz matrix
\[
\widetilde{\mathbb{K}}^{*(N)}_{\alpha,M} =
\begin{pmatrix}
	\tilde{C}^{(M)}_0 & r^{-\alpha}\tilde{C}^{(M)}_{-1} & r^{-2\alpha}\tilde{C}^{(M)}_{-2} & \cdots \\
	r^\alpha\tilde{C}^{(M)}_1 & \tilde{C}^{(M)}_0 & r^{-\alpha}\tilde{C}^{(M)}_{-1} & \cdots \\
	r^{2\alpha}\tilde{C}^{(M)}_2 & r^\alpha\tilde{C}^{(M)}_1 & \tilde{C}^{(M)}_0 & \cdots \\
	\vdots & \vdots & \vdots & \ddots
\end{pmatrix}.
\]
Equivalently, its $(j, k)$ block is $r^{(j-k)\alpha}\tilde{C}^{(M)}_{j-k}$. Note that the spectrum of this matrix is independent of $\alpha$. Assembling $\widetilde{\mathbb{K}}^{*(N)}_{\alpha,M}$ serves to build the corresponding $M \times M$ approximation of the infinite-chain symbol:
\[
\kappa^{(M)}_\alpha(\theta) = \sum_{m\in\mathbb{Z}} r^{m\alpha} \tilde{C}^{(M)}_m e^{im\theta}, \qquad \theta \in [-\pi, \pi].
\]
The exponential decay of the Toeplitz blocks allows this series to be truncated once the remaining terms fall below a prescribed numerical tolerance. The essential spectrum is then approximated by computing the eigenvalues of $\kappa^{(M)}_\alpha(\theta)$ along a sufficiently fine discretization of the phase variable $\theta$.

\subsection{Numerical results}

\begin{figure}[!htbp]
	\centering
	\includegraphics[width=\textwidth]{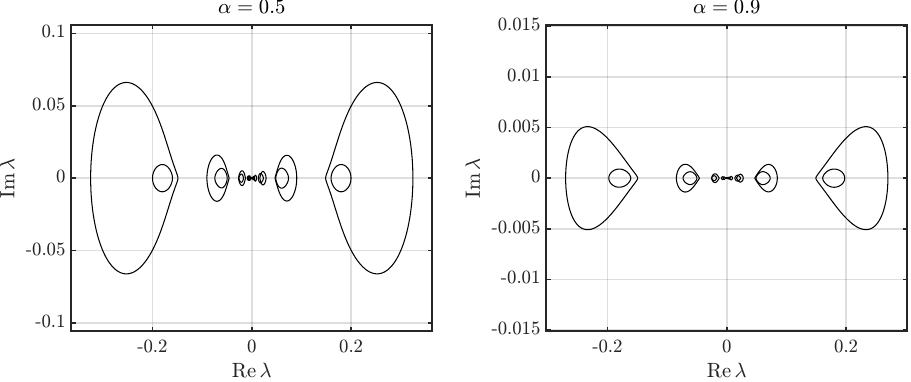}%
	\caption{Numerical approximation of the essential spectrum for a self-similar chain of disks. The horizontal and vertical axes represent the real and imaginary parts of the spectrum, respectively. The spectral bands are computed via a Nystr\"om discretization of the Floquet--Toeplitz symbol $\kappa_\alpha^{(M)}(\theta)$ across different Sobolev space scaling $\alpha$.  Geometric parameters: $r = 0.4$ with reference disk $\Gamma_0$ centred at $(1,0)$ and radius $0.4$. Each disk is discretised with $M = 48$.}
	\label{fig:essential-spectrum-disks}
\end{figure}

We use the discretised symbol $\kappa^{(M)}_\alpha(\theta)$ to approximate the essential spectrum of the infinite chain. For a sufficiently fine discretization of $\theta \in [-\pi, \pi]$, we compute the eigenvalues of $\kappa^{(M)}_\alpha(\theta)$ and take their union over the phase variable. According to Theorem~\ref{thm:essential-spectrum}, these spectral curves approximate
\[
\sigma_{\mathrm{ess}}(\widetilde{K}^*_\Gamma) = \bigcup_{\theta\in[-\pi,\pi]} \sigma\bigl(\kappa_\alpha(\theta)\bigr).
\]
This allows us to visualise directly how the essential spectrum depends on the functional space through the parameter $\alpha = (d - 1)/p - s$, as shown in Figure~\ref{fig:essential-spectrum-disks}. Furthermore, the plotted curves provide a visual estimate for the Fredholm regions based on Corollary~\ref{cor:index-winding}, by tracking the winding number of the eigenvalue loops.

\begin{figure}[!htbp]
	\centering
	\includegraphics[width=0.7\textwidth] {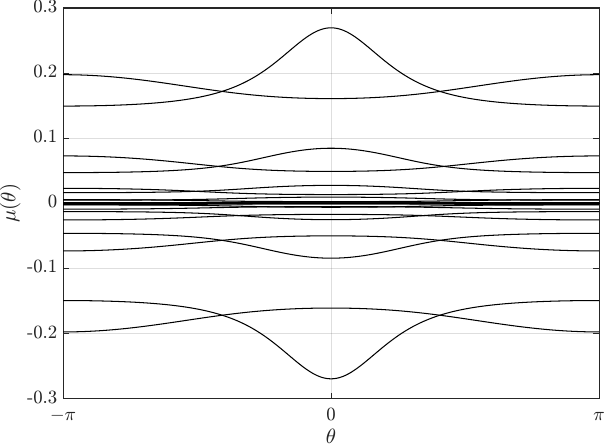}%
\caption{Numerical spectral bands for the chain of disks in the physical energy space ($\alpha = 1$). The vertical axis represents the continuous, real spectral values obtained by tracing the eigenvalues of the discretised symbol $\kappa_{1}^{(M)}(\theta)$ over the phase variable $\theta \in [-\pi, \pi]$.  Geometric parameters: $r = 0.4$ with reference disk $\Gamma_0$ centred at $(1,0)$ and radius $0.4$. Each disk is discretised with $M = 48$.}
	\label{fig:dispersion-disks}
\end{figure}

In the physical space ($\alpha = 1$), we correspondingly obtain the continuous spectral bands plotted in Figure~\ref{fig:dispersion-disks}.

\begin{figure}[!htbp]
	\centering
	\includegraphics[width=0.7\textwidth] {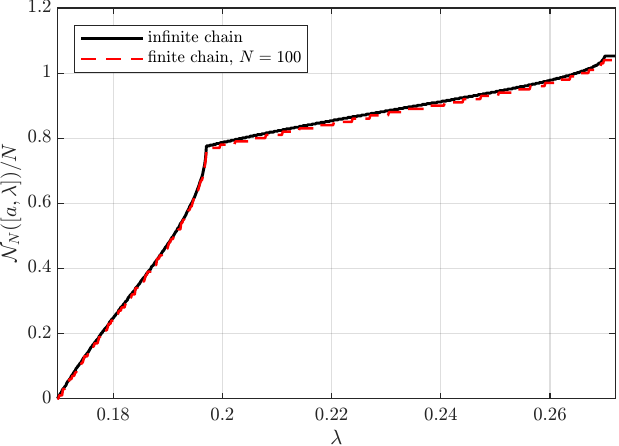}%
\caption{Normalised eigenvalue counting function $\mathcal N_N([a,\lambda])/N$ for a truncated chain of $N=100$ disks, compared against its limit value computed from the numerical chain symbol (Corollary~\ref{cor:density-of-states}). We show the result for $a = 0.17$ and $\lambda\in[0.17,0.27]$. Geometric parameters: $r = 0.4$ with reference disk $\Gamma_0$ centred at $(1,0)$ and radius $0.4$. Each disk is discretised with $M = 48$.}
	\label{fig:density_of_states-disks}
\end{figure}

Finally, at the physical scaling $\alpha = d/2 = 1$, we compare the finite-chain eigenvalue counting function with the symbol-based computation prediction of Corollary~\ref{cor:density-of-states}. We use the notation introduced in \eqref{eq:counting function I_lambda}: for a compact interval $[a,\lambda]  \subset \mathbb{R} \setminus \{0\}$, we compute
\[
\mathcal{N}_N([a,\lambda]) = \#\left\{ j : \lambda^{(N)}_j \in [a,\lambda] \right\}
\]
and compare $\mathcal{N}_N([a,\lambda])/N$ with the corresponding phase average obtained from the eigenvalues of $\kappa^{(M)}_{d/2}(\theta)$. This provides a direct numerical verification of the asymptotic counting formula derived in Corollary~\ref{cor:density-of-states}, shown in Figure~\ref{fig:density_of_states-disks}.

\section{Concluding remarks}\label{sec:conclusion}
In this work, we have established a spectral framework for the Neumann--Poincar\'e operator on self-similar chains of disjoint domains. By constructing a discrete geometric scaling isomorphism, we mapped the different layer potential operators onto block-Toeplitz operator with operator-valued symbols in the Wiener algebra. This allowed us to characterise the essential spectrum in $W_0^{s,p}(\Gamma)$ through the scaling parameter $\alpha = (d - 1)/p - s$, and to derive an exact asymptotic expression for the eigenvalue counting function for large finite truncations.

The operator-valued Wiener--Hopf factorisation presented in Section~\ref{sec:physical space} demonstrate that edge states are governed by a bulk quantity, namely the partial indices of the symbol, which is reminiscent of bulk-boundary correspondences in solid state physics. The explicit geometry of Section~\ref{sec:annuli} sharply illustrates this dependence. Because the symbol within each angular Fourier sector reduces to a $2 \times 2$ matrix with the origin in its bulk spectral gap, the factorisation is rendered canonical or non-canonical strictly by the geometric ratio $r/\rho^{2}$. At the critical value $r = \rho^{2}$, the relevant root reaches the unit circle and the modewise gap closes. Beyond this threshold, every mode simultaneously acquires partial indices $(-1,1)$ and supports a zero-energy state localized on the outermost rings. This mechanism is similar to the localized zero modes found in discrete one-dimensional systems such as the Su--Schrieffer--Heeger model \cite{su1979solitons, ammari2020topological}.

A further question we haven't investigated here concerns the spectral radius. On the energy
space the symmetrisation confines the spectrum to $[-\frac{1}{2},\frac{1}{2}]$,
but on other function spaces neither the bound nor the shape of the essential
spectrum is inherited, and determining the radius has proved challenging for
Lipschitz boundaries \cite{ChandlerWildeHaggerPerfektVirtanen2023}. The same
question has been studied for the layer potentials of elastostatics and
hydrostatics on nonsmooth domains \cite{Mitrea1999, mitrea2002spectra}. In this respect, our setting could be used to study the spectral radius for a non-Lipschitz domain in a range of Sobolev spaces, since the essential spectrum is given
exactly by the fibre spectra of $\kappa_\alpha$ and depends on the functional
setting only through $\alpha$. For the chain of annuli it is available in closed
form, so the spectral radius can be computed as a function of $\alpha$ and the
geometric parameters rather than merely estimated. 

Our approach could be translated to study the spectrum of the Neumann-Poincaré operator in other dilation-invariant configurations, such as corners and cones, where the geometric singularity can be localised and studied as a Wiener-Hopf operator in the half line. By employing analogous asymptotic trace formulas for these continuous operators, one could extract the density of states formulas for smooth geometries approaching a scale-invariant singular limit such as smooth corners with high curvature. For example, evaluating the spectrum of a rounded corner as its radius of curvature tends to zero would establish a rigorous, continuous counterpart to the finite-to-infinite chain asymptotics developed in this work.

Finally, we hope that the techniques put forward in this work could also be translated to other scenarios in wave-physics research where periodic configurations can also be studied from a block-Toeplitz perspective \cite{ammari2020topological, alase2023wiener}.

\appendix
\numberwithin{theorem}{section}
\numberwithin{lemma}{section}
\numberwithin{remark}{section}
\numberwithin{corollary}{section}
\numberwithin{definition}{section}
\numberwithin{equation}{section}

\section{Operator-valued Toeplitz operators}\label{app:toeplitz}

We collect the results from operator-valued Toeplitz theory used in Sections~\ref{sec:spectrum} and~\ref{sec:asymptotics}, in the form required for this work. For the general theory see \cite{BottcherSilbermann1996, BottcherSilbermann2006}.

Throughout, $X$ is a complex Banach space, $\mathcal{L}(X)$ and $\mathcal{K}(X)$ denote the bounded and the compact operators on $X$, and $\mathbb{T}$ is the unit circle. The operator-valued Wiener algebra $\mathcal W(\mathcal{L}(X))$ consists of the functions
\[
a(z)=\sum_{m\in\mathbb{Z}}a_mz^m,
\qquad
a_m\in\mathcal{L}(X),
\qquad
\|a\|_W:=\sum_{m\in\mathbb{Z}}\|a_m\|_{\mathcal{L}(X)}<\infty,
\]
and we write $\mathcal W(\mathcal{K}(X))$ when in addition every $a_m$ is compact. For $a\in \mathcal W(\mathcal{L}(X))$, the Toeplitz operator
\[
\bigl(T(a)\bigr)_{j,k}:=a_{j-k},\qquad j,k\ge0,
\]
is bounded on $\ell^p(\mathbb{N},X)$ for every $p\in(1,\infty)$, with norm at most $\|a\|_{W}$, by the Schur test in the proof of Theorem~\ref{thm:boundedness}.

The symbols occurring in this paper are scalar multiples of the identity plus a compact-valued term, which is the situation in which Fredholmness admits the following characterisation.
\begin{theorem}[Fredholm criterion]\label{thm:app-fredholm}
	Let $k\in \mathcal W(\mathcal{K}(X))$ and $\lambda\in\mathbb{C}$. Then $\lambda I-T(k)=T(\lambda I-k)$ is Fredholm on $\ell^p(\mathbb{N},X)$ if and only if $\lambda I-k(z)$ is boundedly invertible on $X$ for every $z\in\mathbb{T}$.
\end{theorem}
\begin{proof}
	For $\lambda\neq 0$ this is the identity-plus-compact case of the Fredholm theory for Toeplitz operators with operator-valued symbols in the Wiener algebra \cite{BottcherSilbermann1996}. Independence of $p$ holds by \cite{Lindner2008}. For $\lambda=0$, none of the operators in the theorem statement can be Fredholm: $k(z)$ is compact on an infinite-dimensional space, hence not invertible, and $T(k)$ is not Fredholm, since a normalised weakly null sequence in $X$, placed in a single component, is mapped to norm zero by each block and hence by $T(k)$.
\end{proof}

\begin{remark}\label{rem:app-pindependence}
	The independence of $p$ is worth isolating. In this paper the exponent $p$ enters both through the fibre space $X=W^{s,p}_0(\Gamma_0)$ and through the sequence space $\ell^p$. Theorem~\ref{thm:app-fredholm} shows that the second route contributes nothing, so the whole dependence of the spectrum on the functional setting is carried by the fibre and by the scaling exponent $\alpha$.
\end{remark}

\section{Operator-valued Wiener--Hopf factorisation}\label{app:wiener-hopf}

The factorisation theory described here and used in Section~\ref{sec:physical space} can be found in \cite{GohbergLeiterer2009}. As above $X$ denotes a Banach space.

\begin{definition}\label{def:app-factorisation}
	Let $a:\mathbb{T}\to\mathcal{L}(X)$. A \emph{right Wiener--Hopf factorisation} of $a$ with respect to the unit circle $\mathbb{T}$ is a representation
	\[
	a(z)=a_-(z)\Lambda(z)a_+(z),
	\qquad z\in\mathbb{T},
	\]
	in which $a_+$ and $a_+^{-1}$ extend holomorphically to the interior unit disk $\mathbb{D}$, $a_-$ and $a_-^{-1}$ extend holomorphically to the exterior domain $\mathbb{E} = \{z : |z| > 1\} \cup \{\infty\}$, all four with continuous boundary values on $\mathbb{T}$, and
	\[
	\Lambda(z)=I-\Pi+\sum_{j=1}^{J}z^{\nu_j}\Pi_j,
	\qquad
	\Pi:=\sum_{j=1}^{J}\Pi_j,
	\]
	where $\Pi_1,\dots,\Pi_J$ are rank-one idempotents with $\Pi_i\Pi_j=0$ for $i\neq j$ and $\nu_1,\dots,\nu_J\in\mathbb{Z}\setminus\{0\}$ are the \emph{non-zero partial indices}, not necessarily distinct. The factorisation is \emph{canonical} if $J=0$, that is $\Lambda\equiv I$.
\end{definition}

\begin{theorem}[Factorisation of identity plus compact]\label{thm:app-factorisation}
	Let $k(z)$ be compact-valued and holomorphic on an annulus containing $\mathbb{T}$, and set $a(z):=I+k(z)$. If $a(z)$ is boundedly invertible for every $z\in\mathbb{T}$, then $a$ admits a right Wiener--Hopf factorisation with finitely many non-zero partial indices.
\end{theorem}

\begin{proof}
	See \cite{GohbergLeiterer2009}. Holomorphy on an annulus is what makes the family finitely meromorphic and the partial indices finite in number.
\end{proof}

\begin{theorem}[Partial indices]\label{thm:app-partial-indices}
	Let $a=a_-\Lambda a_+$ be a right Wiener--Hopf factorisation as above. Then
	\[
	\dim\ker T(a)=\sum_{\nu_j<0}(-\nu_j),
	\qquad
	\dim\operatorname{coker}T(a)=\sum_{\nu_j>0}\nu_j,
	\qquad
	\operatorname{ind}T(a)=-\sum_{j=1}^{J}\nu_j,
	\]
	and $T(a)$ is invertible if and only if the factorisation is canonical.
\end{theorem}

\begin{proof}
	The analytic factors induce invertible Toeplitz operators, so the defect numbers are those of $T(\Lambda)$, which splits into scalar shifts $T(z^{\nu_j})$ on the ranges of the $\Pi_j$ and the identity elsewhere \cite{GohbergLeiterer2009}.
\end{proof}

\begin{remark}\label{rem:app-eigenvalue-criterion}
	Applied to $a_\lambda:=I-\lambda^{-1}k$ with $\lambda\neq0$, and combined with $\lambda I-T(k)=\lambda T(a_\lambda)$, this is the exact criterion used in Section~\ref{sec:physical space} for $\lambda$ outside the essential spectrum, $\lambda$ is an eigenvalue of $T(k)$ if and only if the factorisation of $a_\lambda$ is non-canonical. In that case, $\dim\ker(\lambda I-T(k))=\sum_{\nu_j<0}(-\nu_j)$. When the Fredholm index vanishes, the partial indices sum to zero, so non-canonicity forces at least one negative partial index.
\end{remark}

\section{The Szeg\H o limit for block-Toeplitz operators}\label{app:szego}

This appendix develops the Szeg\H{o}-type trace asymptotic used as input to
Section~\ref{sec:asymptotics}. For self-adjoint matrix-valued symbols the result is
classical \cite{MirandaTilli2000,Widom1974,BottcherSilbermann2006}. We adapt it here
to an operator-valued symbol that is not self-adjoint, but for which each fibre and each finite
section is self-adjoint in an equivalent inner product.

Throughout, $X$ and $Y$ are
separable Hilbert spaces and $\mathfrak S_q(X,Y)$, $q\in[1,\infty)$, is the
Schatten class of order $q$, with
$\mathfrak S_q(X):=\mathfrak S_q(X,X)$. We use the ideal property in its two-space
form \cite{GohbergKrein1969}: if $A\in\mathfrak S_q(X_1,X_2)$, $B\in\mathcal L(X_0,X_1)$ and
$C\in\mathcal L(X_2,X_3)$, then $CAB\in\mathfrak S_q(X_0,X_3)$ with
\begin{equation}\label{eq:app-ideal}
	\|CAB\|_{\mathfrak S_q}\le\|C\|\,\|A\|_{\mathfrak S_q}\|B\| ;
\end{equation}
only the middle factor need be Schatten. Renorming a Hilbert space by an
equivalent inner product changes Schatten norms by at most a fixed factor and
leaves traces unchanged; both facts are used repeatedly below.

We write $\bigl(H(a)\bigr)_{j,k}:=a_{j+k+1}$, $j,k\ge0$, for the Hankel operator
generated by $a\in \mathcal W(\mathcal{L}(X))$, which is bounded with norm at
most $\|a\|_{W}$, and we use the Widom product formula \cite{Widom1974}
\begin{equation}\label{eq:widomproduct}
	T(ab)=T(a)T(b)+H(a)H(\tilde b),
	\qquad
	\tilde b(z):=b(1/z).
\end{equation}

For $N\ge1$ let $Q_N$ denote the projection of $\ell^p(\mathbb N,X)$ onto its first
$N$ components. The finite sections of $T(a)$ are
\[
T_N(a):=Q_N\,T(a)\,Q_N ,
\qquad
\bigl(T_N(a)\bigr)_{j,k}=a_{j-k},\quad 0\le j,k<N,
\]
regarded as operators on $\ell^p(\{0,\dots,N-1\},X)$.

\begin{definition}[Schatten--Wiener classes]\label{def:app-schattenwiener}
	For $a(z)=\sum_{m\in\mathbb Z}a_mz^m$ with $a_m\in\mathcal L(X,Y)$ put
	\begin{align*}
		a\in \mathcal W^1\bigl(\mathfrak S_q(X,Y)\bigr)
		&\iff a_m\in\mathfrak S_q(X,Y)\ \text{for all }m,\quad
		\textstyle\sum_m(1+|m|)\|a_m\|_{\mathfrak S_q}<\infty,\\[2pt]
		a\in \mathcal W^1_\flat\bigl(\mathfrak S_q(X,Y)\bigr)
		&\iff a_m\in\mathfrak S_q(X,Y)\ \text{for }m\neq0,\quad
		\textstyle\|a_0\|+\sum_{m\neq0}(1+|m|)\|a_m\|_{\mathfrak S_q}<\infty .
	\end{align*}
\end{definition}

\begin{remark}
	The second class in Definition~\ref{def:app-schattenwiener} motivated by the structure of the Hankel operators, since
	$\bigl(H(a)\bigr)_{j,k}=a_{j+k+1}$ never involves $a_0$. This is necessary since the diagonal coefficient of the
	symbol studied in this work is the Neumann--Poincar\'e operator of the reference component, which is
	Schatten only for $q>d-1$ and in particular not Hilbert--Schmidt in dimension three.
\end{remark}

\begin{lemma}[Hankel and finite-section bounds]\label{lem:app-hankel}
	Let $a\in \mathcal W^1_\flat(\mathfrak S_q(X,Y))$. Then $H(a)\in\mathfrak S_q$
	with
	\[
	\|H(a)\|_{\mathfrak S_q}\le\sum_{m\ge0}(1+m)^{1/q}\|a_{m+1}\|_{\mathfrak S_q} .
	\]
	If moreover $a\in \mathcal W^1(\mathfrak S_q(X,Y))$, then
	\[
	\|T_N(a)\|_{\mathfrak S_q}\le N^{1/q}\sum_{m\in\mathbb Z}\|a_m\|_{\mathfrak S_q},
	\qquad N\ge1 .
	\]
\end{lemma}

\begin{proof}
	Decompose $H(a)$ into block anti-diagonals: the $m$-th carries $a_{m+1}$
	repeated $m+1$ times, so its Schatten norm is
	$(m+1)^{1/q}\|a_{m+1}\|_{\mathfrak S_q}$; the coefficient $a_0$ does not occur.
	Likewise $T_N(a)$ decomposes into block diagonals, the $m$-th carrying $a_m$
	exactly $N-|m|$ times, and here $a_0$ occurs $N$ times, whence the stronger
	hypothesis. The triangle inequality gives both bounds.
\end{proof}

\begin{corollary}[Trace-class corrections]\label{cor:app-corrections}
	Let $a\in \mathcal W^1_\flat(\mathfrak S_2(X,Y))$ and
	$b\in \mathcal W^1_\flat(\mathfrak S_2(Y,X))$. Then
	\[
	\sup_{N\ge1}\bigl\|T_N(ab)-T_N(a)T_N(b)\bigr\|_{\mathfrak S_1}<\infty .
	\]
\end{corollary}

\begin{proof}
	Since $T_N(\cdot)=Q_N T(\cdot)Q_N$, the Widom product formula
	\eqref{eq:widomproduct} gives
	\[
	T_N(ab)-T_N(a)T_N(b)
	=Q_NH(a)H(\tilde b)Q_N+Q_NT(a)(I-Q_N)T(b)Q_N .
	\]
	The first term is a product of two Hankel operators, each in $\mathfrak S_2$ by
	Lemma~\ref{lem:app-hankel} with bounds independent of $N$. In the second, the
	entries of $Q_NT(a)(I-Q_N)$ are $a_{j-k}$ with $0\le j<N\le k$, so each
	coefficient $a_m$ with $m\le-1$ occurs at most $|m|$ times along a partial
	isometry; decomposing as in the proof of Lemma~\ref{lem:app-hankel} gives
	$\|Q_NT(a)(I-Q_N)\|_{\mathfrak S_2}\le\sum_{m\le-1}|m|^{1/2}\|a_m\|_{\mathfrak S_2}$,
	and likewise for $(I-Q_N)T(b)Q_N$. Neither term involves $a_0$ or $b_0$, so both
	are products of $\mathfrak S_2$ operators with bounds independent of $N$, hence
	lie in $\mathfrak S_1$ uniformly in $N$.
\end{proof}

\begin{remark}\label{rem:app-dimension}
	No relation between $q$ and the dimension is required by
	Corollary~\ref{cor:app-corrections}, whose hypothesis constrains only the
	off-diagonal coefficients; in the application these are integral operators with
	smooth kernels between disjoint surfaces, hence of every Schatten order in every
	dimension. The condition $q>d-1$ enters only through the diagonal coefficient,
	in the finite-section bound of Lemma~\ref{lem:app-hankel} and in the moment
	estimate of Theorem~\ref{thm:app-szego}.
\end{remark}

\begin{definition}[Admissible test functions]\label{def:app-testfunctions}
	Let $\mathcal F$ be the set of continuous $f:\mathbb R\to\mathbb R$ with
	$0\notin\operatorname{supp}f$.
\end{definition}

Note that if $A$ is compact with $\|A\|\le A_0$ and real spectrum, and
$f\in\mathcal F$ vanishes on $(-\delta,\delta)$, then $f(A)$ is supported on the
spectral subspace associated with $\{\delta\le|x|\le A_0\}$, whose dimension is at
most $\delta^{-q}\|A\|^q_{\mathfrak S_q}$; in particular $f(A)$ is finite rank, and
$f$ may be approximated by polynomials divisible by $x^{\lceil q\rceil}$.

\begin{theorem}[First Szeg\H o limit theorem, operator-valued symbols]\label{thm:app-szego}
	Let $a\in \mathcal W^1(\mathfrak S_q(X))\cap \mathcal W^1_\flat\bigl(\mathfrak S_2(X)\bigr)$ and write
	$\ell^2_N:=\ell^2(\{0,\dots,N-1\},X)$. Assume:
	\begin{enumerate}
		\item[(i)] for every $\theta$, the operator $a(\theta)$ is self-adjoint for
		some inner product on $X$ equivalent to the one in $X$;
		\item[(ii)] for every $N$, the operator $T_N(a)$ is self-adjoint for some
		inner product on $\ell^2_N$ equivalent to the standard one.
	\end{enumerate}
	Then for every $f\in\mathcal F$,
	\[
	\lim_{N\to\infty}\frac1N\operatorname{Tr}f\bigl(T_N(a)\bigr)
	=\frac1{2\pi}\int_{-\pi}^{\pi}\operatorname{Tr}f\bigl(a(\theta)\bigr)\,d\theta .
	\]
\end{theorem}

\begin{proof}
	Put $A_0:=\sum_m\|a_m\|_{\mathfrak S_q}$.
	
	\smallskip
	\noindent\textbf{Step 0: both sides are defined.} Since
	$\|\cdot\|\le\|\cdot\|_{\mathfrak S_q}$ we have $\|T_N(a)\|\le A_0$ and
	$\|a(\theta)\|\le A_0$, while Lemma~\ref{lem:app-hankel} gives
	$\|T_N(a)\|_{\mathfrak S_q}\le N^{1/q}A_0$ and
	$\|a(\theta)\|_{\mathfrak S_q}\le A_0$. Both operators are compact, and by (i)
	and (ii) each is self-adjoint in an equivalent Hilbert structure. Their spectra
	are therefore real, they admit spectral decompositions
	\[
	T_N(a)=\sum_j\lambda^{(N)}_jP^{(N)}_j,
	\qquad
	a(\theta)=\sum_j\mu_j(\theta)P_j(\theta),
	\]
	and for any continuous $f$ we set $f(\cdot):=\sum_jf(\lambda_j)P_j$ as usual.
	Eigenvalues, spectral projections and traces are unchanged by the renorming, so
	all quantities below may be computed in the original structure. In particular
	Weyl's inequality applies there \cite{GohbergKrein1969} and gives
	\begin{equation}\label{eq:app-moment}
		\frac1N\sum_j\bigl|\lambda^{(N)}_j\bigr|^{q}\le A_0^{q},
		\qquad
		\sum_j\bigl|\mu_j(\theta)\bigr|^{q}\le A_0^{q}.
	\end{equation}
	If $f\in\mathcal F$ vanishes on $(-\delta,\delta)$, then only the eigenvalues
	with $|\lambda_j|\ge\delta$ contribute, and by \eqref{eq:app-moment} there are at
	most $\delta^{-q}NA_0^q$ of them for $T_N(a)$ and $\delta^{-q}A_0^q$ for
	$a(\theta)$. Hence $f(T_N(a))$ and $f(a(\theta))$ are finite rank, and
	\begin{equation}\label{eq:app-trace-eigen}
		\operatorname{Tr}f\bigl(T_N(a)\bigr)=\sum_jf\bigl(\lambda^{(N)}_j\bigr),
		\qquad
		\operatorname{Tr}f\bigl(a(\theta)\bigr)=\sum_jf\bigl(\mu_j(\theta)\bigr).
	\end{equation}
	
	\smallskip
	\noindent\textbf{Step 1: the finite sections almost respect products.} We claim
	that for every $k\ge1$
	\begin{equation}\label{eq:app-power-correction}
		\sup_{N\ge1}\bigl\|T_N(a)^k-T_N(a^k)\bigr\|_{\mathfrak S_1}<\infty .
	\end{equation}
	Both classes of Definition~\ref{def:app-schattenwiener} are stable under
	products. Hence $a^k$ lies in both
	classes for every $k$. Now \eqref{eq:app-power-correction} is trivial for $k=1$,
	and assuming it for $k$,
	\[
	T_N(a)^{k+1}-T_N(a^{k+1})
	=T_N(a)\bigl[T_N(a)^{k}-T_N(a^{k})\bigr]
	+\bigl[T_N(a)T_N(a^{k})-T_N(a\,a^{k})\bigr],
	\]
	where the first bracket is uniformly bounded in $\mathfrak S_1$ by hypothesis and
	multiplied by an operator of norm at most $A_0$, and the second is uniformly
	bounded in $\mathfrak S_1$ by Corollary~\ref{cor:app-corrections}, applicable
	because $a$ and $a^k$ both lie in $\mathcal W^1_\flat(\mathfrak S_2(X))$.
	
	\smallskip
	\noindent\textbf{Step 2: the trace of a finite section.} This is where the
	theorem comes from. Let $b$ be a symbol with
	$\sum_m\|b_m\|_{\mathfrak S_1}<\infty$. All $N$ diagonal blocks of $T_N(b)$ are
	equal to $b_0$, and $b_0=\frac1{2\pi}\int_{-\pi}^{\pi}b(\theta)\,d\theta$. Since $\operatorname{Tr}$ is bounded on $\mathfrak S_1$ it passes through
	the integral, so the normalised trace of a finite section is exactly the fibre
	average:
	\begin{equation}\label{eq:app-section-trace}
		\frac1N\operatorname{Tr}T_N(b)=\operatorname{Tr}b_0
		=\frac1{2\pi}\int_{-\pi}^{\pi}\operatorname{Tr}b(\theta)\,d\theta .
	\end{equation}
	For $k>q$ the Fourier coefficients $(a^k)_m$ of $a^k$ are trace class, by H\"older's inequality for
	Schatten norms applied to each term of the convolution, so
	\eqref{eq:app-section-trace} applies to $b=a^k$. Combining it with
	\eqref{eq:app-power-correction},
	\begin{equation}\label{eq:app-monomial}
		\frac1N\operatorname{Tr}T_N(a)^k
		=\frac1{2\pi}\int_{-\pi}^{\pi}\operatorname{Tr}a(\theta)^k\,d\theta
		+O(N^{-1}),\qquad k>q ,
	\end{equation}
	and by linearity the theorem holds for every polynomial all of whose monomials
	have degree exceeding $q$.
	
	\smallskip
	\noindent\textbf{Step 3: from polynomials to $\mathcal F$.} Let $f\in\mathcal F$
	vanish on $(-\delta,\delta)$ and let $\varepsilon>0$. Fix an integer $\kappa>q$
	and set $g(x):=x^{-\kappa}f(x)$, which is continuous on $[-A_0,A_0]$ because $f$
	vanishes near the origin. Choose a polynomial $h$ with
	$\sup_{[-A_0,A_0]}|g-h|\le\varepsilon$ and put $p:=x^{\kappa}h$, so that
	\eqref{eq:app-monomial} applies to $p$ termwise and
	\[
	|f-p|\le\varepsilon|x|^{\kappa}\le\varepsilon A_0^{\kappa-q}|x|^{q}
	\qquad\text{on }[-A_0,A_0].
	\]
	By \eqref{eq:app-trace-eigen} and the moment bounds \eqref{eq:app-moment},
	\[
	\frac1N\bigl|\operatorname{Tr}f(T_N(a))-\operatorname{Tr}p(T_N(a))\bigr|
	\le\varepsilon A_0^{\kappa},
	\qquad
	\bigl|\operatorname{Tr}f(a(\theta))-\operatorname{Tr}p(a(\theta))\bigr|
	\le\varepsilon A_0^{\kappa},
	\]
	the second uniformly in $\theta$. Applying \eqref{eq:app-monomial} to $p$ gives
	\[
	\limsup_{N\to\infty}
	\left|\frac1N\operatorname{Tr}f(T_N(a))
	-\frac1{2\pi}\int_{-\pi}^{\pi}\operatorname{Tr}f(a(\theta))\,d\theta\right|
	\le2\varepsilon A_0^{\kappa},
	\]
	and letting $\varepsilon\to0$ completes the proof.
\end{proof}

\printbibliography
\end{document}